\documentclass[a4paper,10pt]{article}
\usepackage[margin=0.80in,bottom=0.90in,top=0.85in]{geometry}
\usepackage[T1]{fontenc}
\usepackage[utf8]{inputenc}
\usepackage[english]{babel}
\usepackage[dvipsnames]{xcolor}

\usepackage{amsmath,amssymb,amsfonts}%
\usepackage{amsthm}%
\usepackage{mathrsfs}%
\usepackage[title]{appendix}%
\usepackage[colorlinks=true,    
            linkcolor=blue,    
            citecolor=blue,     
            urlcolor=blue]{hyperref}
\usepackage{cleveref}
\usepackage{upgreek}
\usepackage{enumitem}
\usepackage{graphicx}
\DeclareMathOperator{\sech}{sech}

\newcommand{\overbar}[1]{\mkern 1.5mu\overline{\mkern-1.5mu#1\mkern-1.5mu}\mkern 1.5mu}

\theoremstyle{plain}
\newtheorem{theorem}{Theorem}[section]
 
\newtheorem{lemma}[theorem]{Lemma}
\newtheorem{proposition}[theorem]{Proposition}

\theoremstyle{definition}
\newtheorem{definition}{Definition}[section]

\theoremstyle{remark}
\newtheorem{remark}{Remark}

\numberwithin{equation}{section}

\begin{document}

\title{Unstable Stokes waves \protect\\
in constant vorticity flows of finite depth}

\author{
Ting-Yang Hsiao%
\thanks{International School for Advanced Studies (SISSA),
Via Bonomea 265, 34136 Trieste, Italy.
\newline
\textit{Email:} \texttt{thsiao@sissa.it}}
,\quad
Vera Mikyoung Hur%
\thanks{Department of Mathematics,
University of Illinois Urbana--Champaign,
Urbana, IL 61801, USA.
\newline
\textit{Email:} \texttt{verahur@illinois.edu}}
,\quad
Zhao Yang%
\thanks{State Key Laboratory of Mathematical Sciences,
Academy of Mathematics and Systems Science,
Chinese Academy of Sciences,
Beijing 100190, China.
\newline
\textit{Email:} \texttt{yangzhao@amss.ac.cn}}
}
\date{}

\maketitle

\begin{abstract}
We study the spectral instability of small-amplitude Stokes waves in finite-depth flows with constant vorticity.  A center-manifold reduction and the periodic Evans function give two explicit indices, one for long-wave instability near the origin and the other for an order-two resonance away from the origin.  We prove instability when the low-frequency index is positive or the high-frequency index is nonzero, and obtain the leading expansions of the unstable spectral points.  At low frequency, first- and second-order monodromy coefficients determine the leading figure-eight curve. At high frequency, we prove that the resonance exists uniquely and derive a limiting ellipse.  An algebraic identity expresses the high-frequency index as a positive factor times the square of a real coupling coefficient.  Thus every parameter pair with nonzero coupling gives high-frequency instability.  The explicit formulas also permit a comparison of the two criteria in the parameter plane.
\end{abstract}

\tableofcontents

\section{Introduction}\label{sec:intro}

We study the spectral instability of small periodic gravity waves in water of finite depth with constant vorticity.  The fluid is incompressible and inviscid, and the motion is two dimensional.  We consider perturbations with the same constant vorticity as the underlying flow.  Our purpose is to determine how unstable spectrum arises near the origin and at a nonzero resonance, and to calculate the corresponding instability criteria in terms of the vorticity and wave number.

Periodic traveling water waves are commonly called Stokes waves, following the work of Stokes \cite{Sto47}.  Their existence at small amplitude was proved by Nekrasov \cite{Nek21}, Levi-Civita \cite{LC25}, and Struik \cite{Str26} in the irrotational setting.  For rotational flows, we refer to \cite{CS04,Wah09,CV11} and the references therein.

The Benjamin--Feir instability concerns long-wave perturbations and the spectrum near the origin of the complex plane.  For irrotational gravity waves, the formal calculations of \cite{Ben67,Whi67,ZK70} predict a change of modulational stability at the dimensionless depth $1.3627827\ldots$ in the zero-amplitude limit.  Bridges and Mielke \cite{BM95} proved instability in finite depth above this threshold.  Nguyen and Strauss \cite{NS23} subsequently proved the deep-water instability.  Berti, Maspero, and Ventura \cite{BMV22,BMV23} described the local spectrum in deep and finite depth, including its closed figure-eight shape in the unstable regime.  The same authors treated the degenerate case at the critical depth in \cite{BMV24b}.

Instability can also occur away from the origin.  Numerical studies by McLean \cite{McL82} and Deconinck and Oliveras \cite{DO11} revealed unstable spectral bands arising from collisions of Floquet branches.  These bands form small isolated curves, often called bubbles or isolas.  Creedon, Deconinck, and Trichtchenko \cite{CDT22} obtained their asymptotic expansions.  Hur and Yang \cite{HY23b} developed a periodic Evans function approach and proved instability at an order-two resonance for finite-depth irrotational waves under an explicit index condition.  More recently, Berti, Maspero, and Ventura \cite{BMV25} treated the first deep-water isola, and Berti, Corsi, Maspero, and Ventura \cite{BCMV25} proved the existence of isolas of arbitrarily high order in finite depth, outside the exceptional depth sets specified in their theorem.  We refer to Berti \cite{Ber24} for a survey of longitudinal instability, and to Creedon, Nguyen, and Strauss \cite{CNS25} for a survey that also discusses transverse perturbations.

Constant vorticity changes both the dispersion relation and the modulational instability criterion.  Thomas, Kharif, and Manna \cite{TKM12} derived a nonlinear Schr\"odinger equation for finite-depth waves on a uniform shear flow and studied the dependence of its instability criterion on the vorticity.  Related modulational questions for dispersive shallow-water models were studied by Hur \cite{Hur19}.  In infinite depth, Guyenne, Kairzhan, and Sulem \cite{GKS22} derived a Hamiltonian Dysthe equation and compared its predictions with numerical solutions of the Euler equations.

For the full finite-depth Euler equations with nonzero constant vorticity, Hsiao \cite{Hsi24} established a modulational instability criterion in his PhD thesis. Using Kato's similarity transform, Fan and Sun \cite{FS26} rigorously recover the modulational stability criterion predicted by \cite{TKM12}. Building on Hsiao's analysis, we refine the low-frequency instability result, describe the leading unstable spectrum near the origin, and establish spectral instability at an order-two resonance away from the origin.

In infinite depth, Chen and Nguyen \cite{CN26} rigorously establish the Benjamin--Feir stability transition away from the critical vorticity. Independently, our companion paper \cite{HHY26} obtains the same criterion, determines the first unstable isola, and proves spectral instability of sufficiently small nonzero Stokes waves for every fixed constant vorticity.

We follow the periodic Evans function approach of \cite{HY23b}.  After fixing the fluid domain, we formulate the spectral problem as an evolution equation in the horizontal variable and reduce it to a finite-dimensional system.  The relevant center space is four dimensional at the origin and two dimensional at the nonzero resonance.  Expanding the corresponding monodromy matrices gives two real indices, $\mathrm{ind}_{\mathrm{low}}(\omega,\kappa)$ and $\mathrm{ind}_{\mathrm{high}}(\omega,\kappa)$. Here $\omega$ is the nondimensional vorticity, and the fluid depth and wave period are normalized to $1$ and $2\pi/\kappa$, respectively.

\begin{theorem}[Instability criteria]\label{thm:main}
Fix $\omega\in\mathbb R\setminus\{0\}$ and $\kappa>0$, and let the $2\pi/\kappa$-periodic Stokes branch be normalized as in Proposition~\ref{prop:Stokes}.  For all sufficiently small nonzero amplitudes $\epsilon$, the following statements hold.
\begin{enumerate}[label=\textup{(\roman*)}]
\item If $\mathrm{ind}_{\mathrm{low}}(\omega,\kappa)>0$, then the linearized operator has spectral points with nonzero real part at Bloch exponents of order $\epsilon$.  These points tend to the origin as $\epsilon\to0$, with the expansion in Theorem~\ref{thm:BF}.
\item There is a unique order-two resonance $k_2(\sigma_H)-k_4(\sigma_H)=2\kappa$.  It is simple and satisfies $\sigma_H>\sigma_c$.  If $a_{12}^{(0,2)}\ne0$, equivalently if $\mathrm{ind}_{\mathrm{high}}(\omega,\kappa)\ne0$, then the linearized operator has spectral points with nonzero real part near $\lambda=i\sigma_H$ and $\lambda=-i\sigma_H$, with the expansion in Theorem~\ref{thm:high}.
\end{enumerate}
In either case the Stokes wave is spectrally unstable.
\end{theorem}

The low-frequency index follows from the first- and second-order monodromy coefficients and determines the leading splitting of the carrier modes.  At high frequency, we prove that the required resonance and slope conditions hold throughout the parameter range. An algebraic relation between the two coupling coefficients then expresses the high-frequency index as a positive factor times a square.  Thus a nonzero coupling suffices for instability.  The explicit coefficients and the algebraic verification are given in the appendices.  Figure~\ref{fig:leading-spectra} illustrates the leading figure-eight curve and the limiting ellipse.  Their parameter dependence is compared in Figure~\ref{fig:low-high-stability}.

Section~\ref{sec:formulation} formulates the water wave problem, and Section~\ref{sec:Stokes} constructs the Stokes branch. Section~\ref{sec:spec} develops the spectral reduction and the periodic Evans function.  The low- and high-frequency calculations are carried out in Sections~\ref{sec:BF} and \ref{sec:high-freq}, respectively. Section~\ref{sec:comparison-indices} illustrates the resulting instability criteria.  The appendices contain the auxiliary proofs and explicit coefficient formulas.

\section{The water wave problem}\label{sec:formulation}

We begin with the governing equations.  Let \((X,y)\) be inertial coordinates, let \(\tau\) denote dimensional time, and let \(x=X-c\tau\) be the horizontal coordinate in a frame moving to the right with speed \(c>0\).  The fluid of unit density occupies
\[
 \mathcal D(\tau)=\{(x,y)\in\mathbb R^2:0<y<h+H(x,\tau)\}, \qquad h+H(x,\tau)>0.
\]
Gravity is \(-g\mathbf e_y\) and the pressure on the free surface is taken to be zero.  We use the vorticity convention
\[
 \partial_x U_2-\partial_y U_1=\Omega .
\]
A constant-vorticity velocity field can then be written
\begin{equation}\label{eqn:velocity}
 \mathbf U(x,y,\tau) =\big(-\Omega(y-h)+\Phi_x,\,\Phi_y\big), \qquad \mathbf U-c\mathbf e_x =\big(-\Omega(y-h)-c+\Phi_x,\,\Phi_y\big),
\end{equation}
where \(\mathbf U\) is the physical velocity and \(\mathbf U-c\mathbf e_x\) is the velocity relative to the moving frame.  Incompressibility gives
\begin{equation}\label{eqn:Laplace}
 \Delta\Phi=0\qquad\text{in }\mathcal D(\tau).
\end{equation}
Since the domain is simply connected, \(\Phi\) has a harmonic conjugate \(\Psi\), normalized by
\begin{equation}\label{eqn:CR}
 \Phi_x=\Psi_y,\qquad \Phi_y=-\Psi_x,\qquad \Psi=0\quad\text{on }y=0.
\end{equation}
The last normalization is compatible with the impermeability condition
\begin{equation}\label{eqn:Neumann}
 \Phi_y=0\quad\text{on }y=0,
\end{equation}
because \(\Psi_x=-\Phi_y\) vanishes there.  We retain the separate label
\begin{equation}\label{eqn:Dirichlet}
 \Psi=0\quad\text{on }y=0
\end{equation}
for later reference.

In moving coordinates, the Euler equations are
\begin{equation*} \mathbf U_\tau+ \big((\mathbf U-c\mathbf e_x)\cdot\nabla\big)\mathbf U =-\nabla p-g\mathbf e_y,\qquad \nabla\cdot\mathbf U=0,
\end{equation*}
where \(p\) is the fluid pressure.  Using \((\mathbf U\cdot\nabla)\mathbf U =(\nabla\times\mathbf U)\times\mathbf U +\nabla|\mathbf U|^2/2\), together with \eqref{eqn:velocity} and \eqref{eqn:CR}, one obtains
\[
 \nabla\!\left( \Phi_\tau+\frac12|\nabla\Phi|^2 -\big(\Omega(y-h)+c\big)\Phi_x +\Omega\Psi+p+g(y-h)\right)=0.
\]
A time-dependent additive constant in \(\Phi\) may be used to fix the Bernoulli gauge, so the spatial constant on the right may be denoted by \(B\). Evaluating at \(y=h+H(x,\tau)\), where \(p=0\), yields
\begin{equation*} \Phi_\tau+\frac12\Phi_x^2+\frac12\Phi_y^2 -\Omega H\Phi_x-c\Phi_x+gH+\Omega\Psi=B.
\end{equation*}
The condition that the free surface moves with the fluid is
\begin{equation}\label{eqn:kinematic}
 H_\tau=\Phi_y-\big(\Phi_x-\Omega H-c\big)H_x \qquad\text{on }y=h+H(x,\tau).
\end{equation}
Equations \eqref{eqn:Laplace}--\eqref{eqn:kinematic}, with \eqref{eqn:Neumann} and \eqref{eqn:Dirichlet}, are the dimensional free-boundary problem in the moving frame.

Introduce
\[
 \bar x=\frac{x}{h},\qquad \bar y=\frac{y}{h},\qquad \bar t=\frac{c\tau}{h},\qquad \eta=\frac{H}{h},\qquad \phi=\frac{\Phi}{ch},\qquad \psi=\frac{\Psi}{ch},
\]
and the dimensionless parameters
\begin{equation*} b=\frac{B}{c^2},\qquad \omega=\Omega\sqrt{\frac{h}{g}},\qquad \mu=\frac{\sqrt{gh}}{c}>0.
\end{equation*}
Thus \(1/\mu=c/\sqrt{gh}\) is the Froude number, and \(\omega\mu=\Omega h/c\).  Dropping bars, the equations on the moving domain \(0<y<1+\eta(x,t)\) are
\begin{equation}\label{eqn:surface}
\left\{
\begin{aligned}
 &\phi_x-\psi_y=0,\qquad \phi_y+\psi_x=0
 &&\text{for }0<y<1+\eta,\\
 &\phi_y=0,\qquad \psi=0
 &&\text{on }y=0,\\
 &\eta_t+(\psi_y-\omega\mu\eta-1)\eta_x=\phi_y
 &&\text{on }y=1+\eta,\\
 &\phi_t+\frac12\psi_y^2+\frac12\phi_y^2 -(\omega\mu\eta+1)\psi_y+\mu^2\eta+\omega\mu\psi=b
 &&\text{on }y=1+\eta .
\end{aligned}
\right.
\end{equation}

To fix the domain, set
\begin{equation*} \bar t=t,\qquad \bar x=x,\qquad \bar y=\frac{y}{1+\eta(x,t)}.
\end{equation*}
This is a diffeomorphism whenever \(1+\eta\) is bounded away from zero. After the change of variables, and again dropping bars,
\[
 \partial_x^{\rm old} =\partial_x-\frac{y\eta_x}{1+\eta}\partial_y,\qquad \partial_y^{\rm old}=\frac{1}{1+\eta}\partial_y,\qquad \partial_t^{\rm old} =\partial_t-\frac{y\eta_t}{1+\eta}\partial_y.
\]
Substituting into \eqref{eqn:surface} and then eliminating \(\eta_t\) from the Bernoulli condition by use of the kinematic condition, we arrive at the equivalent fixed-strip system
\begin{equation}\label{eqn:flattened}
\left\{
\begin{aligned}
 &\phi_x-\frac{y\eta_x}{1+\eta}\phi_y -\frac{1}{1+\eta}\psi_y=0
 &&\text{for }0<y<1,\\
 &\psi_x-\frac{y\eta_x}{1+\eta}\psi_y +\frac{1}{1+\eta}\phi_y=0
 &&\text{for }0<y<1,\\
 &\phi_y=0,\qquad\psi=0
 &&\text{on }y=0,\\
 &\eta_t+ \left(\frac{\psi_y}{1+\eta}-\omega\mu\eta-1\right)\eta_x -\frac{\phi_y}{1+\eta}=0
 &&\text{on }y=1,\\
 &\phi_t+ \frac{\left(\dfrac{\psi_y}{1+\eta}-\omega\mu\eta-1\right) \phi_y\eta_x}{1+\eta} -\frac{\phi_y^2}{2(1+\eta)^2} +\frac{\psi_y^2}{2(1+\eta)^2}\\
 &\hspace{4em} +\mu^2\eta+\omega\mu\psi -\frac{(\omega\mu\eta+1)\psi_y}{1+\eta}=b
 &&\text{on }y=1 .
\end{aligned}
\right.
\end{equation}

\section{Stokes waves of sufficiently small amplitude}\label{sec:Stokes}

A Stokes wave is a time-independent solution of \eqref{eqn:flattened}.  It therefore satisfies
\begin{equation}\label{eqn:Stokes}
\left\{
\begin{aligned}
 &\phi_x-\frac{y\eta_x}{1+\eta}\phi_y -\frac{1}{1+\eta}\psi_y=0
 &&\text{for }0<y<1,\\
 &\psi_x-\frac{y\eta_x}{1+\eta}\psi_y +\frac{1}{1+\eta}\phi_y=0
 &&\text{for }0<y<1,\\
 &\phi_y=0,\qquad\psi=0
 &&\text{on }y=0,\\
 &\left(\frac{\psi_y}{1+\eta}-\omega\mu\eta-1\right)\eta_x -\frac{\phi_y}{1+\eta}=0
 &&\text{on }y=1,\\
 &\frac{\phi_y^2+\psi_y^2}{2(1+\eta)^2} +\mu^2\eta+\omega\mu\psi -\frac{(\omega\mu\eta+1)\psi_y}{1+\eta}=b
 &&\text{on }y=1 .
\end{aligned}
\right.
\end{equation}
We impose period \(T=2\pi/\kappa\), with wave number \(\kappa>0\), and fix the translation and Bernoulli gauges by requiring
\begin{equation}\label{eqn:Stokes-normalization}
 \phi(-x,y)=-\phi(x,y),\qquad \psi(-x,y)=\psi(x,y),\qquad \eta(-x)=\eta(x),\qquad \frac1T\int_0^T\eta(x)\,dx=0.
\end{equation}

The Crandall--Rabinowitz theorem determines the local solution set, but its curve parameter is not canonical.  For this reason, before calculating the expansion, we fix the Stokes amplitude by an exact condition. On the odd potential component define
\begin{equation*} \mathfrak a(\phi) :=\frac{2}{T}\int_0^T\phi(x,0)\sin(\kappa x)\,dx .
\end{equation*}
The trace map from \(C_{\rm per}^{m+1,\alpha}(\overline S_T)\) to \(C_{\rm per}^{m+1,\alpha}(\mathbb R/T\mathbb Z)\) is bounded, so \(\mathfrak a\) is a bounded linear functional on the potential space used below.  For
\[
 \phi_1(x,y)=\sin(\kappa x)\cosh(\kappa y)
\]
one has \(\mathfrak a(\phi_1)=1\).  The signed parameter \(\epsilon\) will be chosen so that the complete nonlinear branch, not only its tangent at the origin, satisfies
\begin{equation}\label{eqn:amplitude}
\mathfrak a(\phi(\epsilon))=\epsilon.
\end{equation}

We first establish the local branch which will be used throughout the spectral analysis.  Let \(S_T=(\mathbb R/T\mathbb Z)\times(0,1)\), fix an integer \(m\ge2\) and \(0<\alpha<1\), and write subscripts \({\rm e}\) and \({\rm o}\) for even and odd functions of \(x\). Set
\begin{align*}
 \mathcal X^{m+1,\alpha} :=\big\{(\phi,\psi,\eta,b)\in{}& C_{\rm per,o}^{m+1,\alpha}(\overline S_T) \times C_{\rm per,e}^{m+1,\alpha}(\overline S_T) \times C_{\rm per,e}^{m+1,\alpha}(\mathbb R/T\mathbb Z) \times\mathbb R:\\
 &\psi|_{y=0}=0,\quad \textstyle\int_0^T\eta(x)\,dx=0\big\}, \\
 \mathcal Y^{m,\alpha} :={}&C_{\rm per,e}^{m,\alpha}(\overline S_T) \times C_{\rm per,o}^{m,\alpha}(\overline S_T) \times C_{\rm per,o}^{m,\alpha}(\mathbb R/T\mathbb Z) \times C_{\rm per,e}^{m,\alpha}(\mathbb R/T\mathbb Z).
\end{align*}
For solutions, the second interior equation at \(y=0\), together with \(\psi_x(x,0)=0\), implies \(\phi_y(x,0)=0\).  Thus the remaining bottom condition follows from the equations. On the open set \(\inf_x(1+\eta(x))>0\), the two interior residuals and the negatives of the two surface residuals in \eqref{eqn:Stokes} define a real-analytic mapping
\[
 \mathscr F:\mathcal X^{m+1,\alpha}\times(0,\infty) \longrightarrow\mathcal Y^{m,\alpha}.
\]
The laminar states form the trivial sheet \(\mathscr F(0,0,0,0;\mu)=0\).

\begin{proposition}[Analytic Stokes branch]\label{prop:Stokes}
For every \(\omega\in\mathbb R\) and \(\kappa>0\), define
\begin{equation}\label{eqn:dispersion}
 \mu_0=\frac{-\omega+ \sqrt{\omega^2+4\kappa\coth\kappa}}{2}>0.
\end{equation}
There exist \(\epsilon_0>0\) and a real-analytic curve
\[
 (-\epsilon_0,\epsilon_0)\ni\epsilon\longmapsto (\phi(\epsilon),\psi(\epsilon),\eta(\epsilon),b(\epsilon),\mu(\epsilon))
\]
of solutions of \eqref{eqn:Stokes} such that the value at \(\epsilon=0\) is \((0,0,0,0,\mu_0)\), the exact amplitude identity \eqref{eqn:amplitude} holds, and
\begin{equation}\label{eqn:Stokes-expansion}
 \partial_\epsilon\phi(0)=\sin(\kappa x)\cosh(\kappa y),\quad
 \partial_\epsilon\psi(0)=\cos(\kappa x)\sinh(\kappa y),\quad
 \partial_\epsilon\eta(0)=\sinh(\kappa)\cos(\kappa x),\quad
 \partial_\epsilon b(0)=0.
\end{equation}
In the symmetric class \eqref{eqn:Stokes-normalization}, every sufficiently small non-laminar solution lies on this curve.  Without imposing the parity conditions, the same statement holds up to horizontal translation.
\end{proposition}

\begin{proof}
At a laminar state the linearized equations are
\begin{equation*}\left\{
\begin{aligned}
 &\phi_x-\psi_y=0,\qquad \psi_x+\phi_y=0
 &&\text{for }0<y<1,\\
 &\phi_y=0,\qquad\psi=0
 &&\text{on }y=0,\\
 &\eta_x+\phi_y=0
 &&\text{on }y=1,\\
 &\psi_y-\omega\mu\psi-\mu^2\eta+b=0
 &&\text{on }y=1.
\end{aligned}
\right.
\end{equation*}
For the \(n\)-th Fourier mode, \(q_n=n\kappa\), \(n\ge1\), a kernel element is a multiple of
\[
 \phi=\sin(q_nx)\cosh(q_ny),\qquad \psi=\cos(q_nx)\sinh(q_ny),\qquad \eta=\sinh(q_n)\cos(q_nx),
\]
and the remaining boundary condition is
\begin{equation}\label{eqn:Fourier-symbol}
 \Delta_n(\mu):= q_n\cosh(q_n)-(\mu^2+\omega\mu)\sinh(q_n)=0.
\end{equation}
Equation \(\Delta_1(\mu_0)=0\) is precisely \(\mu_0^2+\omega\mu_0=\kappa\coth\kappa\), which gives \eqref{eqn:dispersion}.  Moreover,
\[
 \frac{d}{dr}(r\coth r) =\frac{\sinh r\cosh r-r}{\sinh^2r}>0 \qquad(r>0).
\]
Consequently \(\Delta_n(\mu_0)\ne0\) for every \(n\ge2\).  The zero Fourier block is also invertible: the mean-zero constraint forces the zero mode of \(\eta\) to vanish, the Cauchy--Riemann equations and the bottom conditions then force the zero modes of \(\phi\) and \(\psi\) to vanish in the chosen parity class, and the last boundary equation determines \(b\).

We include the Fredholm argument, since it is needed in the bifurcation theorem.  Put
\[
 A_\mu=D_{(\phi,\psi,\eta,b)} \mathscr F(0,0,0,0;\mu): \mathcal X^{m+1,\alpha}\longrightarrow\mathcal Y^{m,\alpha}.
\]
After Fourier expansion in \(x\), the interior Cauchy--Riemann system with \(\psi|_{y=0}=0\) is an ordinary boundary-value problem in \(y\).  On the \(n\)-th nonzero block its homogeneous boundary matrix has determinant \(\Delta_n(\mu)\).  At \(\mu=\mu_0\),
\begin{equation*} |\Delta_n(\mu_0)| =\sinh(q_n)\big(q_n\coth q_n-\kappa\coth\kappa\big) \ge c(1+q_n)\sinh(q_n),\qquad n\ge2,
\end{equation*}
for some \(c=c(\kappa)>0\).  Indeed the factor in parentheses is positive for every \(n\ge2\), is asymptotic to \(q_n\), and only finitely many modes remain before the asymptotic estimate applies. To justify the estimate in H\"older spaces, consider \(A_\mu U=(f_1,f_2,g_1,g_2)\).  Write \(\langle\cdot\rangle\) for the horizontal mean, \(r^\perp=r-\langle r\rangle\), and let \(K=\partial_x^{-1}\) denote the periodic primitive with zero mean. Setting \(h=K(g_1-f_2(\cdot,1))\), elimination of \(\phi\) and \(\eta\) gives
\[
 \begin{cases}
 \Delta\psi^\perp=\partial_xf_2-\partial_yf_1^\perp
     &\text{in }S_T,\\
 \psi^\perp(\cdot,0)=0,\\
 \partial_y\psi^\perp(\cdot,1) -(\mu^2+\omega\mu)\psi^\perp(\cdot,1) =g_2^\perp+\mu^2h.
 \end{cases}
\]
On even functions with zero horizontal mean, this scalar operator is a compact Robin perturbation of the Dirichlet--Neumann isomorphism from \(C^{m+1,\alpha}\) to \(C^{m-1,\alpha}\times C^{m,\alpha}\), and is therefore Fredholm of index zero.  At \(\mu=\mu_0\), its kernel is one dimensional by \eqref{eqn:Fourier-symbol}.  The remaining variables are recovered from
\[
 \langle\psi(\cdot,y)\rangle =-\int_0^y\langle f_1(\cdot,z)\rangle\,dz,\qquad \phi=K(f_1+\psi_y),\qquad \eta=\psi^\perp(\cdot,1)+h,
\]
with \(b\) determined by the mean of the last boundary equation. Here \(\phi_y=f_2-\psi_x\), which also gives the required vertical regularity of \(\phi\).  The scalar Schauder estimate and these reconstruction formulas yield
\begin{equation*} \|U\|_{\mathcal X^{m+1,\alpha}} \le C\big(\|A_{\mu_0}U\|_{\mathcal Y^{m,\alpha}} +\|\mathcal P_1U\|_{\mathcal X^{m+1,\alpha}}\big),
\end{equation*}
where \(U=(\phi,\psi,\eta,b)\) and \(\mathcal P_1\) is the projection onto the first Fourier block.  The scalar Fredholm alternative and the reconstruction show that \(A_{\mu_0}\) is Fredholm of index zero, its kernel is one dimensional and is spanned by the vector in \eqref{eqn:Stokes-expansion}, and its range has codimension one.

Finally,
\[
 \partial_\mu\Delta_1(\mu_0) =-(2\mu_0+\omega)\sinh\kappa =-\sqrt{\omega^2+4\kappa\coth\kappa}\,\sinh\kappa\ne0.
\]
Let \(\ell\) be the cokernel functional obtained by solving the interior first-harmonic equations and taking the remaining surface residual.  With the normalization in \eqref{eqn:Stokes-expansion},
\[
 \ell\big(\partial_\mu A_{\mu_0} \partial_\epsilon U(0)\big) =\partial_\mu\Delta_1(\mu_0)\ne0.
\]
Thus \(\partial_\mu A_{\mu_0}\partial_\epsilon U(0)\) does not belong to \(\operatorname{Ran}A_{\mu_0}\), which is precisely the Crandall--Rabinowitz transversality condition.  The analytic Crandall--Rabinowitz theorem \cite{CR71} first gives an analytic curve
\[
 s\longmapsto (\widetilde U(s),\widetilde\mu(s))
\]
in the parity-fixed space, oriented so that
\[
 \widetilde U(0)=0,\qquad \widetilde\mu(0)=\mu_0,\qquad \partial_s\widetilde U(0)=U_1,
\]
where \(U_1\) is the vector displayed in \eqref{eqn:Stokes-expansion}.  Put
\[
 A(s):=\mathfrak a(\widetilde\phi(s)).
\]
Then
\[
 A(0)=0,\qquad A'(0)=\mathfrak a(\partial_s\widetilde\phi(0)) =\mathfrak a(\phi_1)=1.
\]
The analytic inverse-function theorem therefore gives a unique analytic inverse \(s=s(\epsilon)\) near the origin.  Defining
\[
 U(\epsilon):=\widetilde U(s(\epsilon)),\qquad \mu(\epsilon):=\widetilde\mu(s(\epsilon))
\]
proves the exact identity \eqref{eqn:amplitude} without changing the first derivatives in \eqref{eqn:Stokes-expansion}.  The local uniqueness statement in the Crandall--Rabinowitz theorem, together with the local injectivity of \(A\), shows that this is the unique local parameterization of the nontrivial branch satisfying \eqref{eqn:amplitude}.  In particular, the condition eliminates the freedom to replace the branch parameter by \(\epsilon+c_2\epsilon^2+c_3\epsilon^3+\cdots\).

In the space without parity restrictions the first-harmonic kernel is the two-dimensional span of the sine and cosine modes. Every nonzero vector in that plane is a horizontal translate of a sine mode.  Imposing the phase condition that the first cosine coefficient of \(\phi\) vanish gives a local slice transverse to the translation orbit; the preceding uniqueness on that slice proves local uniqueness modulo translation.
\end{proof}

We henceforth use the exact amplitude coordinate \eqref{eqn:amplitude}.  Analyticity gives the convergent expansions
\begin{equation}\label{eqn:branch-exp}
\begin{aligned}
 \phi(\epsilon)&=\epsilon\phi_1+\epsilon^2\phi_2 +\epsilon^3\phi_3+O(\epsilon^4),
 &\psi(\epsilon)&=\epsilon\psi_1+\epsilon^2\psi_2 +\epsilon^3\psi_3+O(\epsilon^4),\\
 \eta(\epsilon)&=\epsilon\eta_1+\epsilon^2\eta_2 +\epsilon^3\eta_3+O(\epsilon^4),
 &b(\epsilon)&=\epsilon b_1+\epsilon^2b_2 +\epsilon^3b_3+O(\epsilon^4),\\
 \mu(\epsilon)&=\mu_0+\epsilon\mu_1+\epsilon^2\mu_2 +\epsilon^3\mu_3+O(\epsilon^4),&&
\end{aligned}
\end{equation}
in the spaces used above.  If \((\tau f)(x,y)=f(x+T/2,y)\), then periodicity and \(\kappa T/2=\pi\) give the exact identity
\[
 \mathfrak a(\tau\phi)=-\mathfrak a(\phi).
\]
The traveling-wave equations are translation invariant, and the half-period translate of an odd, respectively even, periodic function is again odd, respectively even.  Thus \(\tau U(\epsilon)\) satisfies the same parity and mean conditions as \(U(\epsilon)\), while its amplitude is \(-\epsilon\).  Exact amplitude normalization and local uniqueness therefore give
\[
 (\phi,\psi,\eta)(x,y;-\epsilon) =(\phi,\psi,\eta)(x+T/2,y;\epsilon),\qquad b(-\epsilon)=b(\epsilon),\qquad \mu(-\epsilon)=\mu(\epsilon).
\]
In particular,
\begin{equation*} b_1=b_3=\mu_1=\mu_3=0.
\end{equation*}
Moreover, applying the bounded linear functional \(\mathfrak a\) coefficientwise to \eqref{eqn:branch-exp} and using \eqref{eqn:amplitude} yields
\begin{equation*} \mathfrak a(\phi_1)=1,\qquad \mathfrak a(\phi_j)=0\quad(j\ge2).
\end{equation*}
These are the normalization conditions used in the order-by-order Stokes calculation.

Collecting the coefficient of \(\epsilon^2\) in \eqref{eqn:Stokes}, a straightforward calculation reveals
\begin{equation*}\left\{
\begin{aligned}
 {\phi_2}_x-{\psi_2}_y
 &=y{\eta_1}_x{\phi_1}_y-\eta_1{\psi_1}_y,\\
 {\psi_2}_x+{\phi_2}_y
 &=y{\eta_1}_x{\psi_1}_y+\eta_1{\phi_1}_y,
 &&0<y<1,\\
 {\phi_2}_y&=0,\qquad\psi_2=0,
 &&y=0,\\
 {\phi_2}_y+{\eta_2}_x
 &={\phi_1}_x{\eta_1}_x-\omega\mu_0\eta_1{\eta_1}_x +\eta_1{\phi_1}_y,
 &&y=1,\\
 {\psi_2}_y-\omega\mu_0\psi_2-\mu_0^2\eta_2
 &=\frac12\big({\psi_1}_y^2+{\phi_1}_y^2\big) +(1-\omega\mu_0)\eta_1{\psi_1}_y-b_2,
 &&y=1 .
\end{aligned}
\right.
\end{equation*}
The first two equations in turn imply
\begin{align*}
 {\phi_2}_{xx}+{\phi_2}_{yy}&=\kappa^2\sinh\kappa\,\sin(2\kappa x)\cosh(\kappa y)-\frac32\kappa^3\sinh\kappa\, y\sin(2\kappa x)\sinh(\kappa y),
 &&0<y<1,\\
 {\phi_2}_y&=0,&&y=0.
\end{align*}
Solving the Fourier equations and imposing zero mean on \(\eta_2\) gives
\begin{align*}
 \phi_2(x,y)
 &=A\sin(2\kappa x)\cosh(2\kappa y) +\frac{\kappa\sinh\kappa}{2} \sin(2\kappa x)y\sinh(\kappa y),\\
 \psi_2(x,y)
 &=A\cos(2\kappa x)\sinh(2\kappa y) +\frac{\kappa\sinh\kappa}{2} \big(\cos(2\kappa x)+1\big)y\cosh(\kappa y),\\
 \eta_2(x)&=E\cos(2\kappa x),
\end{align*}
where
\begin{align}
 A&=\frac{3\kappa}{8\sinh^2\kappa} -\frac{\omega\big(\mu_0^3\sinh\kappa +2\kappa\mu_0\cosh\kappa\big)} {8\kappa\sinh\kappa},\notag\\
 E&=A\sinh(2\kappa) +\frac{\kappa\sinh\kappa\cosh\kappa}{2} -\frac{\omega\mu_0\sinh^2\kappa}{4},\label{eqn:coefE}\\
 b_2&=\frac34\kappa^2\cosh^2\kappa -\frac14\kappa^2\sinh^2\kappa -\frac12\kappa\mu_0^2\sinh\kappa\cosh\kappa -\frac14\kappa\mu_0\omega\sinh(2\kappa) =\frac{\kappa^2}{4}.\label{eqn:coefb2}
\end{align}
The last equality follows from \(\mu_0^2+\omega\mu_0=\kappa\coth\kappa\).  The third term in \eqref{eqn:coefb2} is the contribution of the constant Fourier component of the order-two Bernoulli condition.

We take matters further and record the order-three solvability condition, which will be needed later through \(\mu_2\).  For compactness put \(s=\sinh\kappa\), \(c=\cosh\kappa\). The polynomials $N_1$ and $N_2$ are given in Appendix~\ref{app:Stokes-order3}. Projection of the order-three boundary equations onto the first harmonic yields the scalar Fredholm solvability condition
\begin{equation}\label{eqn:coefmu2}
 \mu_2=-\frac{N_1}{\kappa(2\mu_0+\omega)s^3} +\frac{N_2}{\kappa(2\mu_0+\omega)s^2}.
\end{equation}
Its denominator never vanishes because \(2\mu_0+\omega=\sqrt{\omega^2+4\kappa\coth\kappa}>0\).

To distinguish solvability from normalization, write the third-order equation abstractly as
\[
 A_{\mu_0}U_3=F_3+\mu_2G_1.
\]
The cokernel condition determines \(\mu_2\) uniquely because the Crandall--Rabinowitz transversality coefficient is nonzero.  Once that condition is satisfied, \(U_3\) is determined only up to the addition of \(C U_1\), since \(\ker A_{\mu_0}=\operatorname{span}\{U_1\}\).  In the Green-function representation, this freedom is precisely the first-harmonic homogeneous term
\[
 C\sin(\kappa x)\cosh(\kappa y)
\]
in \(\phi_3\).  Because \(\mathfrak a(\phi_1)=1\), the coefficient condition \(\mathfrak a(\phi_3)=0\) forces \(C=0\).  Thus the Fredholm condition determines \(\mu_2\), while the exact amplitude normalization fixes the kernel component of \(U_3\); the two roles are separate and together make the third-order expansion unambiguous.

In the irrotational limit \(\omega=0\),
\begin{equation*} \mu_2=-\frac{\mu_0^5 (8\sinh^4\kappa+8\sinh^2\kappa+9)} {16\cosh^2\kappa}.
\end{equation*}
The coefficient printed in \cite{HY23b} is expressed in a different gauge: its potential contains a mean term \(\overbar\phi_2x\), whereas \eqref{eqn:branch-exp} requires \(\phi\) to be periodic and lets \(b\) vary.  In the irrotational case,
\[
 \overbar\phi_2=\frac14\mu_0^4\tanh^2\kappa,\qquad \mu_2^{\rm HY}=\mu_2-\mu_0\overbar\phi_2 =-\frac{\mu_0^5 (8\sinh^4\kappa+12\sinh^2\kappa+9)} {16\cosh^2\kappa}.
\]
Thus the two expressions describe the same physical branch after the gauge conversion.

\section{The spectral problem}
\label{sec:spec}

Throughout this section, \(\kappa>0\) is fixed and
\begin{equation}\label{eqn:mu0}
 \mu _0=\frac{-\omega+\sqrt{\omega^2+4\kappa\coth\kappa}}2>0, \qquad \mu_0^2+\omega\mu_0=\kappa\coth\kappa.
\end{equation}
We turn the attention to the spectral problem.  The new constant-vorticity results of this paper are stated for \(\omega\ne0\).  The operator formulation and the zero-amplitude dispersion calculation below, however, are valid for every \(\omega\in\mathbb R\), including \(\omega=0\).

\subsection{The linearized problem}
\label{sec:linearize}

Let \((\phi_\epsilon,\psi_\epsilon,\eta_\epsilon, \mu_\epsilon,b_\epsilon)\) be the real-analytic family of \(T=2\pi/\kappa\)-periodic Stokes waves constructed in Section~\ref{sec:Stokes}, with the crest placed at \(x=0\).  We assume the regularity furnished there, enlarged if necessary by elliptic bootstrapping, so that all coefficient traces occurring below belong to \(C^2_{\rm per}(\mathbb R)\).

Linearizing the flattened equations at this wave and introducing the Alinhac good unknowns
\begin{equation*} \phi\mapsto\phi- \frac{y\phi_{\epsilon,y}}{1+\eta_\epsilon}\eta, \qquad \psi\mapsto\psi- \frac{y\psi_{\epsilon,y}}{1+\eta_\epsilon}\eta,
\end{equation*}
gives, after inserting the normal-mode ansatz \(e^{\lambda t}\),
\begin{equation}\label{eqn:linearized}
\left\{
\begin{aligned}
&\phi_x-\frac{y\eta_{\epsilon,x}}{1+\eta_\epsilon}\phi_y -\frac1{1+\eta_\epsilon}\psi_y=0,
&&0<y<1,\\
&\psi_x-\frac{y\eta_{\epsilon,x}}{1+\eta_\epsilon}\psi_y +\frac1{1+\eta_\epsilon}\phi_y=0,
&&0<y<1,\\
&A_1\eta_x+\frac{\eta_{\epsilon,x}}{(1+\eta_\epsilon)^2}\psi_y -\frac1{(1+\eta_\epsilon)^2}\phi_y+A_2\eta=0,
&&y=1,\\
&\lambda\phi+A_1\psi_y+A_3\eta+A_4\phi_y +\omega\mu_\epsilon\psi=0,
&&y=1,\\
&\phi_y=0,\qquad \psi=0,&&y=0,
\end{aligned}
\right.
\end{equation}
where all unmarked surface traces are evaluated at \(y=1\), and
\begin{align}
A_1={}&\frac{\psi_{\epsilon,y}}{(1+\eta_\epsilon)^2} -\frac{\omega\mu_\epsilon\eta_\epsilon+1}{1+\eta_\epsilon}, \notag\\
A_2={}&\frac{\lambda-\omega\mu_\epsilon\eta_{\epsilon,x}} {1+\eta_\epsilon} +\frac{\psi_{\epsilon,yy}\eta_{\epsilon,x}-\phi_{\epsilon,yy}} {(1+\eta_\epsilon)^3}, \notag\\
A_3={}&\mu_\epsilon^2+ \frac{\phi_{\epsilon,y}\phi_{\epsilon,yy} +\psi_{\epsilon,y}\psi_{\epsilon,yy}}{(1+\eta_\epsilon)^3} -\frac{(\omega\mu_\epsilon\eta_\epsilon+1)\psi_{\epsilon,yy}} {(1+\eta_\epsilon)^2}, \notag\\
A_4={}&\frac{\phi_{\epsilon,y}}{(1+\eta_\epsilon)^2}.
\label{coefficients:A}
\end{align}
The calculation is in Appendix~\ref{A:good-unknown-proof}.  Notice, in particular, that no division by \(\omega\) is involved.

Since \(A_1=-1+O(\epsilon)\) and \(A_3=\mu_0^2+O(\epsilon)\) in \(C^1_{\rm per}\), there is \(\epsilon_0>0\) such that
\begin{equation}\label{eqn:A13}
 \inf_x|A_1(x)|\ge\frac12, \qquad \inf_x|A_3(x)|\ge\frac{\mu_0^2}{2} \quad (|\epsilon|<\epsilon_0).
\end{equation}
Define \(\chi\) by
\begin{align}\label{eqn:barpsi}
\chi(x,y)=\frac{\mu_0^2}{A_3(x)}\bigg\{
&\lambda\int_0^y\phi(x,s)\,ds +\omega\mu_\epsilon\int_0^y\psi(x,s)\,ds\notag\\
&+A_1(x)\psi(x,y) +A_4(x)(\phi(x,y)-\phi(x,0))\bigg\}.
\end{align}
Then the two surface conditions in \eqref{eqn:linearized} become
\begin{equation*} \chi_y(x,1)+\mu_0^2\eta(x)=0, \qquad \chi(x,0)=0.
\end{equation*}
By \eqref{eqn:A13} it defines a bounded isomorphism on each Sobolev space used below.  Its inverse is
\begin{align}
\psi(x,y)={}&e^{-\frac{\omega\mu_\epsilon}{A_1(x)}y} \int_0^y e^{\frac{\omega\mu_\epsilon}{A_1(x)}s} \left[-\frac{\omega\mu_\epsilon A_3(x)}{\mu_0^2A_1(x)^2}\chi(x,s) +\frac{\omega\mu_\epsilon A_4(x)-\lambda A_1(x)}{A_1(x)^2}
\phi(x,s)\right]ds \notag\\
&+\frac{A_3(x)}{\mu_0^2A_1(x)}\chi(x,y) -\frac{A_4(x)}{A_1(x)} \left(\phi(x,y)-e^{-\frac{\omega\mu_\epsilon}{A_1(x)}y}\phi(x,0)\right).
\label{Psi:recover}
\end{align}
Notice that \eqref{Psi:recover} has the regular limit obtained by putting \(\omega=0\) directly into the formula.

Set
\begin{equation}\label{def:V}
 \mathcal V_\lambda\binom{\phi}{\chi}(y) :=e^{\omega\mu_0y}\int_0^y e^{-\omega\mu_0s} \big(-\omega\mu_0\chi(s)+\lambda\phi(s)\big)\,ds.
\end{equation}
For later reference we denote this function by \(\upsilon\).  Let
\begin{align*}
Y&=H^1(0,1)\times H^1_D(0,1)\times\mathbb C,\\
H^s_D(0,1)&:=\{f\in H^s(0,1):f(0)=0\} \quad(s>1/2),\qquad H^0_D(0,1):=L^2(0,1),\\
\mathcal D&=\left\{(\phi,\chi,\eta)\in H^2(0,1)\times H^2_D(0,1)\times\mathbb C: \phi_y(0)=0,\ \chi_y(1)+\mu_0^2\eta=0\right\}.
\end{align*}
Both spaces are complex Hilbert spaces, and \(\mathcal D\) is equipped with the product \(H^2\times H^2\times\mathbb C\) norm. Define \(\mathbf L(\lambda):\mathcal D\subset Y\to Y\) by
\begin{equation}\label{eqn:operatorL}
\mathbf L(\lambda)
\begin{pmatrix}\phi\\ \chi\\ \eta\end{pmatrix}
=
\begin{pmatrix}
\omega\mu_0\upsilon-\omega\mu_0\chi-\chi_y+\lambda\phi\\
\lambda\upsilon-\lambda\chi- \omega\mu_0(\phi-\phi(0))+\phi_y\\
\lambda\eta-\phi_y(1)
\end{pmatrix}.
\end{equation}
Putting together the preceding transformations, we arrive at
\begin{equation}\label{eqn:spatial-system}
 \mathbf u_x= \big(\mathbf L(\lambda)+\mathbf B(x;\lambda,\epsilon)\big)\mathbf u, \qquad \mathbf u=(\phi,\chi,\eta)^T,
\end{equation}
where
\begin{equation}\label{eqn:B}
 \mathbf B\in C^1_{\rm per} \big(\mathbb R;\mathcal L(\mathcal D,Y)\big), \qquad \mathbf B(x;\lambda,0)=0.
\end{equation}
For every compact \(K\subset\mathbb C\), it is analytic in \((\lambda,\epsilon)\in K\times(-\epsilon_0,\epsilon_0)\), locally uniformly in \(x\). Appendix~\ref{A:operator-properties} proves that \(\mathbf L(\lambda)\) is closed on the common dense domain \(\mathcal D\), has compact resolvent, and that the graph norm of \(\mathbf L(\lambda)\) is locally uniformly equivalent to the \(\mathcal D\)-norm.

\subsection{Spectral stability and instability}
\label{sec:L}

Notice that the spectral problem is an analytic operator pencil in \(\lambda\), rather than the spectrum of a single operator denoted by \(\mathcal L(\epsilon)\).  We make the following definition. Let
\begin{align*}
X_{\rm per}&=L^2_{\rm per}(0,T;Y),\\
\mathcal D_{\rm per}&=H^1_{\rm per}(0,T;Y) \cap L^2_{\rm per}(0,T;\mathcal D),
\end{align*}
and, for \(\theta\in[-\kappa/2,\kappa/2)\), define
\begin{equation}\label{eqn:Bloch}
 \mathcal L_\theta(\lambda,\epsilon) :=\partial_x+i\theta-\mathbf L(\lambda) -\mathbf B(x;\lambda,\epsilon): \mathcal D_{\rm per}\subset X_{\rm per}\longrightarrow X_{\rm per}.
\end{equation}

\begin{definition}[The spectrum]\label{def:temporal-spectrum}
The spectrum of the Stokes wave is
\begin{equation}\label{eqn:spec}
 \Sigma(\epsilon)= \left\{\lambda\in\mathbb C: \mathcal L_\theta(\lambda,\epsilon) \text{ is not boundedly invertible for some } \theta\in[-\kappa/2,\kappa/2)\right\}.
\end{equation}
The wave is spectrally stable if \(\Sigma(\epsilon)\subset\{\operatorname{Re}\lambda\le0\}\), and spectrally unstable otherwise.
\end{definition}

This definition is equivalent to the whole-line definition obtained from
\begin{equation*}
 \mathcal L(\lambda,\epsilon)= \partial_x-\mathbf L(\lambda)-\mathbf B(x;\lambda,\epsilon)
\end{equation*}
on \(L^2(\mathbb R;Y)\), with domain \(H^1(\mathbb R;Y)\cap L^2(\mathbb R;\mathcal D)\), because the Bloch transform unitarily decomposes the latter pencil as the direct integral of \eqref{eqn:Bloch}.  More precisely, the whole-line pencil is boundedly invertible if and only if every fiber is invertible and the fiber inverses are uniformly bounded in \(\theta\).  In the spectral neighborhoods considered below the fibers are Fredholm of index zero and depend continuously on \(\theta\); hence failure of invertibility of one fiber is equivalent to failure of invertibility of the direct integral.  The proof is given in Appendix~\ref{A:bloch-proof}.

\begin{lemma}[Spectral symmetries]\label{lem:sym}
For real \(\epsilon\),
\begin{equation}\label{eqn:symm}
 \lambda\in\Sigma(\epsilon) \quad\Longrightarrow\quad \overline\lambda,-\lambda,-\overline\lambda \in\Sigma(\epsilon).
\end{equation}
Consequently spectral stability is equivalent to \(\Sigma(\epsilon)\subset i\mathbb R\).
\end{lemma}

The proof in Appendix~\ref{A:symm} keeps track of the Bloch exponent: complex conjugation and reversibility both replace \(\theta\) by \(-\theta\).  This point is necessary for a correct whole-line statement.

\subsection{The spectrum at zero amplitude}
\label{sec:eps=0}

We next work out the spectrum at zero amplitude.  For \(\epsilon=0\), a Fourier mode \(e^{ikx}\mathbf U(y)\), \(k\in\mathbb R\), satisfies
\begin{equation*} \mathbf L(\lambda)\mathbf U=ik\mathbf U.
\end{equation*}
For \(k\ne0\), Appendix~\ref{A:eps=0} shows directly that a nonzero solution exists precisely when
\begin{equation}\label{eqn:disp}
 d(k,\lambda):= \lambda^2+i\big(\omega\mu_0\tanh k-2k\big)\lambda -k^2+(\mu_0^2+\omega\mu_0)k\tanh k=0.
\end{equation}
Thus \(\lambda=i\sigma\), where
\begin{equation}\label{eqn:dispersion-curve}
 \left(\sigma+\frac{\omega\mu_0\tanh k}{2}-k\right)^2 =\mu_0^2k\tanh k+ \frac{\omega^2\mu_0^2\tanh^2k}{4},
\end{equation}
or equivalently
\begin{equation}\label{eqn:sigma-pm}
 \sigma_\pm(k)=k-\frac{\omega\mu_0\tanh k}{2} \pm\sqrt{\mu_0^2k\tanh k+ \frac{\omega^2\mu_0^2\tanh^2k}{4}}.
\end{equation}
The square root is real because \(k\tanh k\ge0\).  The formula extends continuously through \(k=0\), and
\begin{equation}\label{eqn:sigma-sym}
 \sigma_+(-k)=-\sigma_-(k).
\end{equation}

\begin{lemma}\label{lem:flat-spectrum}
For every \(\omega\in\mathbb R\), the flat wave is spectrally stable and
\begin{equation}\label{eqn:spec0}
 \Sigma(0)=i\mathbb R.
\end{equation}
\end{lemma}

\begin{proof}
For each Bloch exponent \(\theta\), the admissible Fourier numbers are \(k=\theta+n\kappa\), \(n\in\mathbb Z\).  Equation \eqref{eqn:disp} shows that every corresponding temporal eigenvalue is purely imaginary.  Conversely, \(\sigma_+\) is strictly increasing from \(0\) to \(+\infty\) on \((0,\infty)\), by Lemma~\ref{lem:sigma-geometry} below.  Hence every \(i\sigma\), \(\sigma\ge0\), occurs, and \eqref{eqn:sigma-sym} gives the negative half-axis.  The Bloch direct-integral characterization then yields \eqref{eqn:spec0}.
\end{proof}

\begin{lemma}[Geometry of the dispersion curves]
\label{lem:sigma-geometry}
The following statements hold.
\begin{enumerate}[label=(\roman*)]
\item \(\sigma_+\) is strictly increasing on \((0,\infty)\).
\item \(\sigma_-\) is strictly convex on \((0,\infty)\), satisfies \(\sigma_-(0)=\sigma_-(\kappa)=0\), is negative on \((0,\kappa)\), and is strictly increasing on \((\kappa,\infty)\).
\item There is a unique \(r_c\in(0,\kappa)\) at which \(\sigma_-\) attains its minimum.  Put
\begin{equation*} k_c=-r_c, \qquad \sigma_c=-\sigma_-(r_c)=\sigma_+(k_c)>0.
\end{equation*}
Then \(k_c\) is the unique critical point of \(\sigma_+\) on the negative half-line.
\end{enumerate}
\end{lemma}

The proof is in Appendix~\ref{A:disp}; its final identity is a strictly positive decomposition.

For \(\sigma\ge0\), let \(k_j(\sigma)\) denote the real roots of \(\sigma_\pm(k)=\sigma\), ordered as follows:
\begin{enumerate}[label=(\roman*)]
\item If \(\sigma=0\), then \(k_1=-\kappa\), \(k_2=\kappa\), and \(k_3=k_4=0\).
\item If \(0<\sigma<\sigma_c\), there are four simple roots \(k_1<k_3<0<k_4<k_2\), with \(\sigma_+(k_j)=\sigma\) for \(j=1,3,4\) and \(\sigma_-(k_2)=\sigma\).
\item If \(\sigma=\sigma_c\), then \(k_1=k_3=k_c<0<k_4<k_2\); the root \(k_c\) is double and the other two roots are simple.
\item If \(\sigma>\sigma_c\), there are exactly two simple roots \(0<k_4<k_2\).
\end{enumerate}

\subsection{The reduced spaces and biorthonormal projections}
\label{sec:proj}

For a nonzero root \(k\) of \eqref{eqn:dispersion-curve}, set
\begin{equation*}\Upphi(k,\sigma)=
\begin{pmatrix}
\cosh(ky)\\[2mm]
\dfrac1{ik}\big((k-\sigma)\sinh(ky) -\omega\mu_0\cosh(ky)+\omega\mu_0\big)\\[2mm]
\dfrac{i}{\mu_0^2}\big((k-\sigma)\cosh k -\omega\mu_0\sinh k\big)
\end{pmatrix}.
\end{equation*}
Write \(\Upphi_j(\sigma)=\Upphi(k_j(\sigma),\sigma)\) whenever \(k_j(\sigma)\ne0\).  At \(\sigma=0\), define
\begin{equation*} \Upphi_j(0)=\Upphi(k_j(0),0),\qquad j=1,2,
\end{equation*}
and
\begin{equation*} \Upphi_3(0)=
 \begin{pmatrix}0\\ \omega\mu_0y^2/2-y\\
 (1-\omega\mu_0)/\mu_0^2\end{pmatrix},
 \qquad
 \Upphi_4(0)=\begin{pmatrix}1\\0\\0\end{pmatrix}.
\end{equation*}

\begin{lemma}[Spectrum of \(\mathbf L(i\sigma)\)]
\label{lem:flat-eigenvalues}
The spatial eigenvalues on the imaginary axis have precisely the following algebraic structure.
\begin{enumerate}[label=(\roman*)]
\item At \(\sigma=0\), \(-i\kappa\) and \(i\kappa\) are simple, with eigenvectors \(\Upphi_1(0)\) and \(\Upphi_2(0)\), respectively. The eigenvalue \(0\) has algebraic multiplicity two and geometric multiplicity one, and
\begin{equation}\label{eqn:Jordan0}
 \mathbf L(0)\Upphi_4(0)=0, \qquad \mathbf L(0)\Upphi_3(0)=\Upphi_4(0).
\end{equation}
\item If \(0<\sigma<\sigma_c\), the four numbers \(ik_j(\sigma)\), \(j=1,2,3,4\), are simple spatial eigenvalues with eigenvectors \(\Upphi_j(\sigma)\).
\item If \(\sigma>\sigma_c\), the only imaginary spatial eigenvalues are the two simple eigenvalues \(ik_2(\sigma)\) and \(ik_4(\sigma)\).
\item At \(\sigma=\sigma_c\), \(ik_2\) and \(ik_4\) are simple. The eigenvalue \(ik_c\) has algebraic multiplicity two and geometric multiplicity one.  A Jordan chain is
\begin{equation}\label{eqn:Jordanc}
 \Upphi_1(\sigma_c)=\Upphi(k_c,\sigma_c), \qquad \Upphi_3(\sigma_c)=-i\,\partial_k \Upphi(k,\sigma_c)\big|_{k=k_c},
\end{equation}
so that
\begin{equation*}
 (\mathbf L(i\sigma_c)-ik_c)\Upphi_1(\sigma_c)=0, \qquad (\mathbf L(i\sigma_c)-ik_c)\Upphi_3(\sigma_c) =\Upphi_1(\sigma_c).
\end{equation*}
\end{enumerate}
No longer Jordan chains occur in (i) or (iv).
\end{lemma}

The proof in Appendix~\ref{A:Jordan} uses the order of the zeros of the scalar determinant \eqref{eqn:disp}; it does not infer algebraic multiplicity from a formal calculation.

Only two spectral anchors are needed below.  The low-frequency anchor is \(\sigma_0=0\).  For the high-frequency calculation, let \(\sigma_H>\sigma_c\) satisfy
\begin{equation}\label{eqn:resonance}
 k_2(\sigma_H)-k_4(\sigma_H)=2\kappa.
\end{equation}
Put
\begin{equation}\label{eqn:center-spaces}
 \begin{aligned}
 \mathbf Y(0)&=\operatorname{span}\{\Upphi_1(0),\Upphi_2(0), \Upphi_3(0),\Upphi_4(0)\},\\
 \mathbf Y(\sigma_H)&=\operatorname{span}\{\Upphi_2(\sigma_H), \Upphi_4(\sigma_H)\}.
 \end{aligned}
\end{equation}

We use the inner product, linear in its first entry,
\begin{equation}\label{eqn:Y-product}
\langle(\phi_1,\chi_1,\eta_1),(\phi_2,\chi_2,\eta_2)\rangle_Y =\int_0^1(\phi_1\overline{\phi_2} +\phi_{1,y}\overline{\phi_{2,y}})\,dy +\int_0^1(\chi_1\overline{\chi_2} +\chi_{1,y}\overline{\chi_{2,y}})\,dy +\eta_1\overline{\eta_2}.
\end{equation}
At \(\sigma=0\), choose \(\uppsi_j(0)\in\mathcal D(\mathbf L(0)^*)\), \(1\le j\le4\), such that
\begin{align}
 \mathbf L(0)^*\uppsi_j(0)&=-ik_j(0)\uppsi_j(0),
 &&j=1,2,\notag\\
 \mathbf L(0)^*\uppsi_3(0)&=0, \qquad\mathbf L(0)^*\uppsi_4(0)=\uppsi_3(0),
 \label{eqn:adjoint-chain-zero-b}\\
 \langle\Upphi_\ell(0),\uppsi_j(0)\rangle_Y&=\delta_{\ell j},
 &&1\le \ell,j\le4.
\notag
\end{align}
At \(\sigma=\sigma_H\), choose \(\uppsi_j(\sigma_H)\in\mathcal D(\mathbf L(i\sigma_H)^*)\), \(j=2,4\), such that
\begin{align*}
 \mathbf L(i\sigma_H)^*\uppsi_j(\sigma_H)
 &=-ik_j(\sigma_H)\uppsi_j(\sigma_H),
 &&j=2,4,\\
 \langle\Upphi_\ell(\sigma_H),\uppsi_j(\sigma_H)\rangle_Y
 &=\delta_{\ell j},
 &&\ell,j\in\{2,4\}.
\end{align*}
The construction and the coefficient formulas are given in Appendix~\ref{A:adjoint-coordinate-formulas}.

We now define the projections directly from these biorthonormal families:
\begin{align}
 P_0u&=\sum_{j=1}^4 \langle u,\uppsi_j(0)\rangle_Y\Upphi_j(0),
 \label{eqn:biorthogonal-projection-zero}\\
 P_Hu&=\sum_{j\in\{2,4\}} \langle u,\uppsi_j(\sigma_H)\rangle_Y\Upphi_j(\sigma_H).
 \label{eqn:biorthogonal-projection-high}
\end{align}

\begin{lemma}\label{lem:biorthogonal-projections}
Each of \(P_0\) and \(P_H\) is a bounded finite-rank projection onto \(\mathbf Y(0)\) and \(\mathbf Y(\sigma_H)\), respectively, and maps \(Y\) into \(\mathcal D\).  Moreover,
\begin{equation*}
 P_0\mathbf L(0)=\mathbf L(0)P_0, \qquad P_H\mathbf L(i\sigma_H)=\mathbf L(i\sigma_H)P_H \quad\hbox{on }\mathcal D.
\end{equation*}
Consequently, both projections give invariant splittings of \(Y\).
\end{lemma}

\begin{proof}
For either projection, biorthonormality gives \(P^2=P\) and identifies its range with the corresponding space in \eqref{eqn:center-spaces}. Boundedness follows because the sums in \eqref{eqn:biorthogonal-projection-zero}-- \eqref{eqn:biorthogonal-projection-high} are finite, and their ranges are contained in \(\mathcal D\).  Finally, the eigenvector relations and the two Jordan-chain relations \eqref{eqn:Jordan0} and \eqref{eqn:adjoint-chain-zero-b} give, for \(u\in\mathcal D\), the two commutation identities above. This also proves the invariant splitting.
\end{proof}

The oscillatory part of the low-frequency ordering is fixed by
\begin{equation*} k_1(0)=-\kappa,\qquad k_2(0)=\kappa, \qquad \mathbf L(0)\Upphi_1=-i\kappa\Upphi_1, \qquad \mathbf L(0)\Upphi_2=i\kappa\Upphi_2.
\end{equation*}
The dual vectors are constructed directly in this order.  Namely,
\begin{equation}\label{eqn:dual-order}
 (\uppsi_1(0),\uppsi_2(0),\uppsi_3(0),\uppsi_4(0)) \longleftrightarrow (\Upphi_{-\kappa},\Upphi_{+\kappa}, \Upphi_0^{\rm gen},\Upphi_0), \qquad \langle\Upphi_\ell(0),\uppsi_j(0)\rangle_Y=\delta_{\ell j}.
\end{equation}
Thus the coefficient column defining \(\uppsi_j(0)\) pairs with \(\Upphi_j\) from the outset.  No intermediate ordering or subsequent permutation is used, and every row and column of the reduced system has the ordering in \eqref{eqn:dual-order}.
\subsection{Reduction of the spectral problem. The periodic Evans function}
\label{sec:reduction}

Fix one of the two spectral anchors \(\sigma_*\in\{0,\sigma_H\}\).  Put \(P=P_0\) if \(\sigma_*=0\), and \(P=P_H\) if \(\sigma_*=\sigma_H\), and let \(Q=I-P\). Write
\begin{equation*} \lambda=i\sigma_*+\delta, \qquad |\delta|+|\epsilon|\ll1,
\end{equation*}
and
\begin{equation}\label{eqn:anchored-system}
 \mathbf u_x=\mathbf L(i\sigma_*)\mathbf u +\widetilde{\mathbf B}(x;\sigma_*,\delta,\epsilon)\mathbf u,
\end{equation}
where
\begin{equation*}
 \widetilde{\mathbf B}= \mathbf L(i\sigma_*+\delta)-\mathbf L(i\sigma_*) +\mathbf B(x;i\sigma_*+\delta,\epsilon).
\end{equation*}
Let \(\mathbf L_h=\mathbf L(i\sigma_*)|_{QY}\).  The key estimate is the following.
\begin{align}
 \mathbf v_x&=\mathbf L(i\sigma_*)\mathbf v +P\widetilde{\mathbf B}(x;\sigma_*,\delta,\epsilon) (\mathbf v+\mathbf w),\label{eqn:center-evolution}\\
 \mathbf w_x&=\mathbf L_h\mathbf w +Q\widetilde{\mathbf B}(x;\sigma_*,\delta,\epsilon) (\mathbf v+\mathbf w).\label{eqn:complement-evolution}
\end{align}

\begin{lemma}[Resolvent estimate]
\label{lem:resolvent}
There is \(C>0\) such that
\begin{equation}\label{eqn:resolvent-estimate}
 \|(\mathbf L_h-ia)^{-1}\|_{\mathcal L(QY)} \le\frac{C}{1+|a|}, \qquad a\in\mathbb R.
\end{equation}
\end{lemma}

The proof is in Appendix~\ref{A:resolvent-estimate}.  In particular, the compact-resolvent property alone is not used as a substitute for the high-frequency estimate.

We shall use the following precise form of the reduction theorem.  It is the linear, nonautonomous specialization of Mielke's reduction principle \cite[Theorem~1]{Mie88}; see also \cite{Mie86,Mie90}.

\begin{theorem}[Periodic linear reduction principle]
\label{thm:Mielke-specialized}
Let \(A:\mathcal D\subset Y\to Y\) be closed and densely defined, and let \(P\) be a bounded finite-rank projection satisfying \(P\mathcal D\subset\mathcal D\) and \(APu=PAu\) for \(u\in\mathcal D\).  Put \(Q=I-P\), \(Y_c=PY\), and \(A_h=A|_{QY}\).  Assume
\begin{equation}\label{eqn:resolvent}
 i\mathbb R\subset\rho(A_h),\qquad \sup_{a\in\mathbb R}(1+|a|) \|(A_h-ia)^{-1}\|_{\mathcal L(QY)}<\infty.
\end{equation}
Let \(R(\cdot;p)\in C^1_{\rm per} (\mathbb R;\mathcal L(\mathcal D,Y))\) depend analytically on the finite-dimensional parameter \(p\), with \(R(x;0)=0\).  Then, for \(|p|\) sufficiently small, there is a unique
\begin{equation*}
 W(\cdot;p)\in C^1_{\rm per} (\mathbb R;\mathcal L(Y_c,Q\mathcal D)), \qquad \|W(\cdot;p)\|\le C|p|,
\end{equation*}
analytic in \(p\), whose graph is invariant under \(u_x=(A+R(x;p))u\).  Every bounded classical solution lies on this graph, and its center component satisfies
\begin{equation*}
 v_x=Av+PR(x;p)(I+W(x;p))v.
\end{equation*}
Conversely every bounded solution of the reduced equation lifts uniquely to a bounded solution of the full equation.
\end{theorem}

We verify the hypotheses.  By Lemma~\ref{lem:L-closed-compact}, the operator \(A=\mathbf L(i\sigma_*)\) is closed and densely defined. Lemma~\ref{lem:biorthogonal-projections} shows that \(P\) is a bounded finite-rank projection which reduces \(A\).  Estimate \eqref{eqn:resolvent-estimate} is exactly \eqref{eqn:resolvent}.  Finally, \(R=\widetilde{\mathbf B}\) is \(T\)-periodic, belongs to \(C^1_{\rm per}(\mathbb R;\mathcal L(\mathcal D,Y))\), is analytic in \(p=(\delta,\epsilon)\), and vanishes at \(p=0\), by \eqref{eqn:B}.  Therefore Theorem~\ref{thm:Mielke-specialized} gives a \(T\)-periodic operator
\begin{equation*} \mathcal W(x;\sigma_*,\delta,\epsilon) \in\mathcal L(\mathbf Y(\sigma_*),Q\mathcal D), \qquad \mathcal W(x;\sigma_*,0,0)=0,
\end{equation*}
analytic in \((\delta,\epsilon)\), such that every bounded solution of \eqref{eqn:anchored-system} is uniquely of the form
\begin{equation*}
 \mathbf u(x)=\mathbf v(x)+ \mathcal W(x;\sigma_*,\delta,\epsilon)\mathbf v(x), \qquad \mathbf v(x)\in\mathbf Y(\sigma_*).
\end{equation*}
The center component satisfies
\begin{equation}\label{eqn:reduced-evolution}
 \mathbf v_x=\mathbf L(i\sigma_*)\mathbf v +P\widetilde{\mathbf B}(x;\sigma_*,\delta,\epsilon) \big(I+\mathcal W(x;\sigma_*,\delta,\epsilon)\big)\mathbf v.
\end{equation}
Since the original problem is linear, \(\mathcal W\) is linear in \(\mathbf v\), and \(\|\mathcal W\|_{\mathcal L(\mathbf Y(\sigma_*),Q\mathcal D)} =O(|\delta|+|\epsilon|)\).  This is the precise form of the reduction used in the subsequent expansions.

Using the ordered basis in \eqref{eqn:center-spaces}, we write \(\mathbf v=\sum_j a_j\Upphi_j(\sigma_*)\).  Then \eqref{eqn:reduced-evolution} becomes
\begin{equation}\label{eqn:center-system}
 \mathbf a_x=\mathbf A(x;\sigma_*,\delta,\epsilon)\mathbf a,
\end{equation}
where \(\mathbf A\) is \(T\)-periodic in \(x\), analytic in \((\delta,\epsilon)\), and has order four at the low-frequency anchor \(\sigma_*=0\).  Let \(\mathbf X\) be the principal fundamental matrix,
\begin{equation*}
 \mathbf X_x=\mathbf A\mathbf X, \qquad \mathbf X(0)=I.
\end{equation*}

\begin{definition}[The periodic Evans function]
\label{def:Delta}
For \(k\in\mathbb C\), define
\begin{equation}\label{eqn:Evans}
 \Delta(\lambda,k;\epsilon) :=\det\big(e^{ikT}I- \mathbf X(T;\sigma_*,\lambda-i\sigma_*,\epsilon)\big).
\end{equation}
\end{definition}

For fixed \(\sigma_*\), the function in \eqref{eqn:Evans} is analytic in \((\lambda,k,\epsilon)\) near \((i\sigma_*,k_*,0)\).  A change of analytic basis multiplies it by a nonvanishing analytic factor and therefore does not change its zero set or zero multiplicities.

\begin{lemma}\label{lem:Evans-spectrum}
There is a neighborhood \(U\) of \(i\sigma_*\) such that, for \(|\epsilon|\) sufficiently small and \(\lambda\in U\),
\begin{equation}\label{eqn:Evans-spectrum}
 \lambda\in\Sigma(\epsilon) \quad\Longleftrightarrow\quad \Delta(\lambda,k;\epsilon)=0 \quad\text{for some }k\in[-\kappa/2,\kappa/2).
\end{equation}
\end{lemma}

\begin{proof}
By the reduction theorem, bounded solutions of the full spatial system with \(\lambda\in U\) are in one-to-one correspondence with bounded solutions of \eqref{eqn:center-system}.  Floquet theory for the latter finite-dimensional periodic system says that a nonzero bounded solution exists precisely when its monodromy matrix has a multiplier of modulus one, equivalently an eigenvalue \(e^{ikT}\) with real \(k\).  This is exactly the vanishing condition in \eqref{eqn:Evans-spectrum}.  Finally use the Bloch characterization \eqref{eqn:spec}.
\end{proof}

\section{The Benjamin--Feir instability}
\label{sec:BF}

In this section we study the spectrum near the origin.  We retain the notation of Sections~\ref{sec:linearize} and~\ref{sec:proj}, and put $\sigma=0$ throughout.  Recall that
\begin{equation*} \mathbf L(0)\Upphi_1(0)=-i\kappa\Upphi_1(0),\qquad \mathbf L(0)\Upphi_2(0)=i\kappa\Upphi_2(0),
\end{equation*}
whereas $\mathbf L(0)\Upphi_3(0)=\Upphi_4(0)$ and $\mathbf L(0)\Upphi_4(0)=0$.  Let $\mathbf A(0)$ be the matrix of $\mathbf L(0)$ on the generalized kernel in this ordered basis.  Thus
\begin{equation*} \mathbf A(0)=
 \begin{pmatrix}
 -i\kappa&0&0&0\\
 0&i\kappa&0&0\\
 0&0&0&0\\
 0&0&1&0
 \end{pmatrix},
 \qquad \mathbf X(T;0,0,0)=e^{\mathbf A(0)T}
 =\begin{pmatrix}
 1&0&0&0\\
 0&1&0&0\\
 0&0&1&0\\
 0&0&T&1
 \end{pmatrix},
 \qquad T=\frac{2\pi}{\kappa}.
\end{equation*}
Thus the Floquet multiplier one has algebraic multiplicity four at $(\delta,\gamma,\epsilon)=(0,0,0)$.

\subsection{Expansion of the monodromy matrix}
\label{sec:a0}

Throughout this subsection, $\sigma=0$.  For brevity, write
\begin{equation*} \mathbf B(x;\delta,\epsilon) :=\widetilde{\mathbf B}(x;0,\delta,\epsilon).
\end{equation*}
Let $v_k(x;\delta,\epsilon)$, $k=1,\ldots,4$, denote the four solutions of the reduced equation whose initial values are
\begin{equation*} v_k(0;\delta,\epsilon)=\Upphi_k(0).
\end{equation*}
The corresponding solutions of the full spatial system are
\begin{equation*} u_k(x;\delta,\epsilon) =v_k(x;\delta,\epsilon)+w(x,v_k;\delta,\epsilon),
\end{equation*}
where the first term belongs to $P_0Y$ and the second to $Q_0Y$.

We write
\begin{equation}\label{eqn:expansion-vk}
 v_k(x;\delta,\epsilon) =\sum_{j=1}^4\left(\sum_{m,n\geq0} a_{jk}^{(m,n)}(x)\delta^m\epsilon^n\right)\Upphi_j(0)
\end{equation}
and
\begin{equation*} w(x,v_k;\delta,\epsilon) =\sum_{m+n\geq1}w_k^{(m,n)}(x)\delta^m\epsilon^n.
\end{equation*}
Here and below, the superscript $(m,n)$ means the coefficient of $\delta^m\epsilon^n$, without factorials.  At zero amplitude and $\delta=0$, a straightforward calculation gives
\begin{equation}\label{eqn:a00}
 \big(a_{jk}^{(0,0)}(x)\big)_{j,k=1}^4 =e^{\mathbf A(0)x}
 =\begin{pmatrix}
 e^{-i\kappa x}&0&0&0\\
 0&e^{i\kappa x}&0&0\\
 0&0&1&0\\
 0&0&x&1
 \end{pmatrix},
 \qquad a_{jk}^{(m,n)}(0)=0\quad(m+n\geq1).
\end{equation}

We next describe how the remaining coefficients are calculated. Expand
\begin{equation}\label{eqn:B-expansion-low}
 \mathbf B(x;\delta,\epsilon) =\sum_{m+n\geq1}\delta^m\epsilon^n \mathbf B^{(m,n)}(x).
\end{equation}
Suppose that all coefficients of total degree less than $m+n$ have already been found.  At the order $\delta^m\epsilon^n$, put
\begin{equation}\label{eqn:forcing-low}
 \mathcal F_k^{(m,n)}(x) =\sum_{\substack{0\leq p\leq m,\ 0\leq q\leq n\\p+q\geq1}}
 \mathbf B^{(p,q)}(x) \left\{w_k^{(m-p,n-q)}(x) +\sum_{\ell=1}^4a_{\ell k}^{(m-p,n-q)}(x) \Upphi_\ell(0)\right\},
\end{equation}
where $w_k^{(0,0)}=0$.  Since $\mathbf B^{(0,0)}=0$, the right side of \eqref{eqn:forcing-low} contains no $w_k^{(m,n)}$ and is made up entirely of coefficients already determined.

Inserting \eqref{eqn:expansion-vk}--\eqref{eqn:B-expansion-low} into the two equations for $v$ and $w$ in \eqref{eqn:center-evolution}--\eqref{eqn:complement-evolution}, and comparing the coefficient of $\delta^m\epsilon^n$, gives
\begin{align}
 \frac{d}{dx}a_k^{(m,n)}
 &=\mathbf A(0)a_k^{(m,n)}+F_k^{(m,n)},
 \label{eqn:coefficient-equation-low}\\
 \frac{d}{dx}w_k^{(m,n)}
 &=\mathbf L_h w_k^{(m,n)}+Q_0\mathcal F_k^{(m,n)}, \qquad P_0w_k^{(m,n)}=0,
 \label{eqn:complementary-coefficient-low}
\end{align}
where
\begin{equation*} a_k^{(m,n)}=
 \begin{pmatrix}a_{1k}^{(m,n)}\\a_{2k}^{(m,n)}\\
 a_{3k}^{(m,n)}\\a_{4k}^{(m,n)}\end{pmatrix},
 \qquad F_k^{(m,n)}=
 \begin{pmatrix}
 \langle\mathcal F_k^{(m,n)},\uppsi_1(0)\rangle_Y\\
 \langle\mathcal F_k^{(m,n)},\uppsi_2(0)\rangle_Y\\
 \langle\mathcal F_k^{(m,n)},\uppsi_3(0)\rangle_Y\\
 \langle\mathcal F_k^{(m,n)},\uppsi_4(0)\rangle_Y
 \end{pmatrix}.
\end{equation*}
The first equation follows simply by pairing the center equation with the four dual vectors.

For clarity, we also display the second equation in components.  Write
\begin{equation*}
 w_k^{(m,n)}=(\phi,\chi,\eta)^T, \qquad Q_0\mathcal F_k^{(m,n)}=(f_1,f_2,f_3)^T,
\end{equation*}
and put
\begin{equation*} \upsilon(x,y)= -\rho e^{\rho y}\int_0^y e^{-\rho s}\chi(x,s)\,ds, \qquad \rho=\omega\mu_0.
\end{equation*}
Then \eqref{eqn:complementary-coefficient-low} reads
\begin{align}
 \phi_x&=\rho\upsilon-\rho\chi-\chi_y+f_1,
 &&0<y<1,\notag\\
 \chi_x&=-\rho\big(\phi-\phi(x,0)\big)+\phi_y+f_2,
 &&0<y<1,\notag\\
 \eta_x&=-\phi_y(x,1)+f_3,\notag\\
 \chi(x,0)&=0,\qquad \phi_y(x,0)=0, \qquad \chi_y(x,1)+\mu_0^2\eta(x)=0,
 \label{eqn:component-complementary-low}
\end{align}
together with $P_0w_k^{(m,n)}=0$.  Thus, once the lower order terms are known, \eqref{eqn:component-complementary-low} is a linear forced problem with explicitly known right side.  In the present calculation the forcing is a finite sum of exponential, trigonometric, and polynomial terms, and the solution is found by the method of undetermined coefficients.  The complementary condition removes the four center homogeneous terms.

After $w_k^{(m,n)}$ has been found, the first equation in \eqref{eqn:coefficient-equation-low}, subject to \eqref{eqn:a00}, gives
\begin{equation}\label{eqn:duhamel-low}
 a_k^{(m,n)}(x) =e^{\mathbf A(0)x}\int_0^x e^{-\mathbf A(0)s} F_k^{(m,n)}(s)\,ds.
\end{equation}
Equations \eqref{eqn:forcing-low}, \eqref{eqn:component-complementary-low}, and \eqref{eqn:duhamel-low} determine the coefficients successively in the total degree $m+n$.  This is the recursive calculation used below.

In what follows we abbreviate
\begin{equation*} a_{jk}^{(m,n)}=a_{jk}^{(m,n)}(T).
\end{equation*}
Translation invariance gives the following useful fact.

\begin{lemma}\label{lem:fourth-neutral-column}
For every $n\geq1$,
\begin{equation*} a_{j4}^{(0,n)}=0,\qquad j=1,\ldots,4.
\end{equation*}
\end{lemma}

\begin{proof}
The vector $\Upphi_4(0)=(1,0,0)^T$ is an exact solution when $\delta=\sigma=0$, and \( \mathbf B(x;0,\epsilon)\Upphi_4(0)=0. \) Hence $\mathbf B^{(0,n)}\Upphi_4(0)=0$ for every $n\geq1$. Starting at $n=1$, equations \eqref{eqn:forcing-low}-- \eqref{eqn:complementary-coefficient-low} give the assertion by induction.
\end{proof}

In particular, $\Delta(0,K\kappa;\epsilon)=0$ for every $K\in\mathbb Z$. A straightforward calculation using \eqref{eqn:duhamel-low} reveals that
\begin{align}
 \big(a_{jk}^{(1,0)}\big)
 &=\begin{pmatrix}
 a_{11}^{(1,0)}&0&a_{13}^{(1,0)}&0\\
 0&a_{11}^{(1,0)}&-a_{13}^{(1,0)}&0\\
 0&0&a_{33}^{(1,0)}&0\\
 a_{41}^{(1,0)}&a_{41}^{(1,0)}&a_{43}^{(1,0)}&a_{44}^{(1,0)}
 \end{pmatrix},
 \label{eqn:a10}\\
 \big(a_{jk}^{(0,1)}\big)
 &=\begin{pmatrix}
 0&0&a_{13}^{(0,1)}&0\\
 0&0&a_{13}^{(0,1)}&0\\
 0&0&0&0\\
 a_{41}^{(0,1)}&-a_{41}^{(0,1)}&a_{43}^{(0,1)}&0
 \end{pmatrix}.
 \label{eqn:a01}
\end{align}
At the next order we shall need only
\begin{equation}\label{eqn:second-order-coefficients-needed}
 a_{11}^{(2,0)},\quad a_{34}^{(2,0)},\quad a_{14}^{(1,1)},\quad a_{31}^{(1,1)},\quad a_{11}^{(0,2)}.
\end{equation}
The remaining entries either vanish or cancel from the invariant combination below.  Useful checks on the calculation are
\begin{align}
 a_{12}^{(2,0)}&=\frac{a_{13}^{(1,0)}a_{41}^{(1,0)}}{T},
 & a_{31}^{(2,0)}&=-\frac{a_{41}^{(1,0)} (a_{11}^{(1,0)}-a_{33}^{(1,0)})}{T}, \notag\\
 a_{11}^{(1,1)}&=\frac{a_{13}^{(0,1)}a_{41}^{(1,0)} +a_{13}^{(1,0)}a_{41}^{(0,1)}}{T},
 & a_{12}^{(1,1)}&=\frac{a_{13}^{(0,1)}a_{41}^{(1,0)} -a_{13}^{(1,0)}a_{41}^{(0,1)}}{T}.
 \label{eqn:relations-amn}
\end{align}
Explicit formulas for the first-order coefficients and the five entries in \eqref{eqn:second-order-coefficients-needed} are given in Appendix~\ref{app:all-explicit-low-order-coefficients}.  In particular, $a_{11}^{(2,0)}$ enters the determinant through $F_{\mathrm{low}}$ in \eqref{eqn:def-FQ-low}.

We shall use no operator coefficient beyond
\begin{equation*} \mathbf B^{(1,0)},\quad \mathbf B^{(0,1)},\quad \mathbf B^{(2,0)},\quad \mathbf B^{(1,1)},\quad \mathbf B^{(0,2)}.
\end{equation*}
Only these five coefficients enter the leading modulational index.

\subsection{The modulational instability index}\label{sec:index-low}

Let
\begin{equation*} D(\delta,\gamma,\epsilon) :=\Delta(\delta,\gamma;\epsilon) =\det\big(e^{i\gamma T}I-\mathbf X(T;0,\delta,\epsilon)\big).
\end{equation*}
The function $D$ is analytic near the origin. Our purpose is to determine the first term which moves the carrier roots away from the imaginary axis, using the monodromy coefficients calculated in the preceding subsection.

Recall that
\begin{equation}\label{eqn:rho-t-mu-low}
 \mu_0=\frac{\sqrt{\omega^2+4\kappa\coth\kappa}-\omega}{2}>0, \qquad \mu_0^2+\omega\mu_0=\kappa\coth\kappa.
\end{equation}
At zero amplitude, \eqref{eqn:sigma-pm} and \eqref{eqn:sigma-sym} give the two carrier roots
\begin{equation}\label{eqn:flat-carrier-roots-low}
 \delta_+(\gamma,0) =i\sigma_+(-\kappa+\gamma) =-i\sigma_-(\kappa-\gamma), \qquad \delta_-(\gamma,0)=i\sigma_-(\kappa+\gamma).
\end{equation}
Here $\kappa$, $\omega$, and hence $\mu_0$ are held fixed as the Fourier wave number varies. The dispersion relation gives
\[
 \mu_0^2\kappa\tanh\kappa +\frac{\omega^2\mu_0^2\tanh^2\kappa}{4} =\frac{(\kappa+\mu_0^2\tanh\kappa)^2}{4}>0.
\]
Thus $\sigma_-$ is analytic near $k=\kappa$, where $\sigma_-(\kappa)=0$. Taylor expansion of \eqref{eqn:flat-carrier-roots-low} yields
\begin{equation*} \delta_\pm(\gamma,0) =ic_g\gamma\pm ic_2\gamma^2+O(\gamma^3),
\end{equation*}
where differentiation of \eqref{eqn:sigma-pm} and use of \eqref{eqn:rho-t-mu-low} give
\begin{equation}\label{eqn:group-velocity-low}
 c_g:=\sigma_-'(\kappa) =\frac{\kappa}{\kappa+\mu_0^2\tanh\kappa} \left(1-\frac{2\kappa}{\sinh(2\kappa)}\right)
\end{equation}
and
\begin{equation*} c_2:=-\frac12\sigma_-''(\kappa) =-\frac{\kappa^2\sech^2\kappa} {\kappa+\mu_0^2\tanh\kappa} -\frac{\mu_0^4(\tanh\kappa-\kappa\sech^2\kappa)^2} {(\kappa+\mu_0^2\tanh\kappa)^3}.
\end{equation*}
In particular, $c_2<0$. Set
\begin{equation}\label{eqn:DeltaL}
 \Delta_L=c_g^2+(\omega\mu_0-2)c_g+1-\kappa\coth\kappa.
\end{equation}
The normalizing factors used below have fixed signs.

\begin{lemma}\label{lem:automatic-low-nonresonance}
For every $\kappa>0$ and $\omega\in\mathbb R$,
\begin{equation*} 0<c_g<1-\frac{2\kappa}{\sinh(2\kappa)}<1, \qquad \Delta_L<-\kappa^2\sech^2\kappa<0.
\end{equation*}
In particular, $c_g\Delta_L\neq0$.
\end{lemma}

\begin{proof}
Since $\mu_0>0$ and $\sinh(2\kappa)>2\kappa$ for $\kappa>0$, the bounds for $c_g$ follow directly from \eqref{eqn:group-velocity-low}. Moreover, the same formula gives
\[
 (\kappa\coth\kappa+\mu_0^2)c_g =\kappa\coth\kappa-\frac{\kappa^2}{\sinh^2\kappa}.
\]
Using this identity and $\omega\mu_0=\kappa\coth\kappa-\mu_0^2$ in \eqref{eqn:DeltaL}, we obtain
\begin{align*}
 \Delta_L
 &=c_g^2+2(\kappa\coth\kappa-1)c_g +1-2\kappa\coth\kappa+\frac{\kappa^2}{\sinh^2\kappa}\\
 &=\big(c_g+\kappa\coth\kappa-1\big)^2-\kappa^2.
\end{align*}
Since $\kappa\coth\kappa>1$, the bounds for $c_g$ imply
\[
 0<c_g+\kappa\coth\kappa-1 <\kappa\coth\kappa-\frac{2\kappa}{\sinh(2\kappa)} =\kappa\tanh\kappa.
\]
Consequently,
\[
 \Delta_L <\kappa^2\tanh^2\kappa-\kappa^2 =-\kappa^2\sech^2\kappa<0.
\]
\end{proof}

We first note the amplitude symmetry
\begin{equation}\label{eqn:Evans-amplitude-parity-low}
 D(\delta,\gamma,-\epsilon)=D(\delta,\gamma,\epsilon).
\end{equation}
Indeed, the exact amplitude normalization and uniqueness in Proposition~\ref{prop:Stokes} identify the wave of amplitude $-\epsilon$ with the half-period translate of the wave of amplitude $\epsilon$.  The corresponding reduced systems are related by translation and a periodic change of coordinates.  Their monodromy matrices are consequently conjugate, which proves \eqref{eqn:Evans-amplitude-parity-low}.

For brevity, write
\begin{equation*} \mathbf a^{(m,n)}=(a_{jk}^{(m,n)})_{j,k=1}^4, \qquad A=a_{11}^{(1,0)}=\frac{T}{c_g}.
\end{equation*}
The last identity follows from the flat dispersion relation, or from the first-order calculation above.  In the notation of \eqref{eqn:expansion-vk},
\begin{align}
 \mathbf X(T;0,\delta,\epsilon) ={}&I+TE_{43} +\delta\mathbf a^{(1,0)}+\epsilon\mathbf a^{(0,1)} \notag\\
 &+\delta^2\mathbf a^{(2,0)} +\delta\epsilon\mathbf a^{(1,1)} +\epsilon^2\mathbf a^{(0,2)} +O\big((|\delta|+|\epsilon|)^3\big),
 \label{eqn:monodromy-series-low}
\end{align}
where $E_{43}$ has its only nonzero entry, equal to one, in position $(4,3)$.  The remainder is analytic.  The following calculation explains why the displayed terms suffice.

\begin{lemma}\label{lem:second-order-sufficiency-low}
Introduce the variables
\begin{equation}\label{eqn:carrier-scaling}
 \gamma=\tau m,\qquad \epsilon=\tau\alpha,\qquad \delta=ic_g\tau m+\tau^2\zeta.
\end{equation}
For bounded $(m,\alpha,\zeta)$ and sufficiently small $\tau$, one has the analytic identity
\begin{equation}\label{eqn:scaled-Evans-sixth-order-low}
 D(ic_g\tau m+\tau^2\zeta,\tau m,\tau\alpha) =\tau^6\bigl(T\det S_2(m,\alpha,\zeta) +\tau\mathcal R(m,\alpha,\zeta,\tau)\bigr),
\end{equation}
where $\mathcal R$ is analytic and $S_2$ depends only on the five monodromy matrices of total degree one and two in \eqref{eqn:monodromy-series-low}.
\end{lemma}

\begin{proof}
Set $M=e^{i\gamma T}I-\mathbf X(T;0,\delta,\epsilon)$ and make the substitution \eqref{eqn:carrier-scaling}.  Then
\begin{equation*} M=-TE_{43}+\tau M_1+\tau^2M_2+O(\tau^3),
\end{equation*}
where
\begin{align}
 M_1={}&iTmI-ic_gm\mathbf a^{(1,0)} -\alpha\mathbf a^{(0,1)}, \notag\\
 M_2={}&-\zeta\mathbf a^{(1,0)}-\frac{T^2m^2}{2}I +c_g^2m^2\mathbf a^{(2,0)} \notag\\
 &-ic_gm\alpha\mathbf a^{(1,1)} -\alpha^2\mathbf a^{(0,2)}.
 \label{eqn:M1-M2-low}
\end{align}
Use the ordered row and column lists
\begin{equation*} I_0=(1,2,3),\qquad J_0=(1,2,4).
\end{equation*}
The first-order matrices \eqref{eqn:a10}--\eqref{eqn:a01}, together with $c_gA=T$, give
\begin{equation*} (M_1)_{I_0,J_0}=0.
\end{equation*}
After ordering the columns of $M$ as $(1,2,4,3)$, we obtain
\begin{equation*} \widetilde M=
 \begin{pmatrix}
 H&b\\ c&\mathfrak p
 \end{pmatrix},
 \quad
 \begin{aligned}
 H&=\tau^2H_2+O(\tau^3),& b&=\tau b_1+O(\tau^2),\\
 c&=\tau c_1+O(\tau^2),& \mathfrak p&=-T+O(\tau),
 \end{aligned}
\end{equation*}
where
\begin{equation*} H_2=(M_2)_{I_0,J_0},\qquad b_1=(M_1)_{I_0,3},\qquad c_1=(M_1)_{4,J_0}.
\end{equation*}
Notice that $\mathfrak p$ is nonzero for small $\tau$.  Subtracting $b\mathfrak p^{-1}$ times the last row from the first three rows gives
\begin{equation}\label{eqn:exact-elimination-low}
 \det M=-\mathfrak p\det(H-b\mathfrak p^{-1}c), \qquad H-b\mathfrak p^{-1}c=\tau^2S_2+O(\tau^3),
\end{equation}
with
\begin{equation}\label{eqn:effective-second-order-matrix-low}
 S_2=H_2+\frac{b_1c_1}{T}.
\end{equation}
The minus sign in \eqref{eqn:exact-elimination-low} accounts for the single column exchange.  All three columns of the remaining matrix are divisible by $\tau^2$, analytically.  Its determinant is therefore divisible by $\tau^6$.  A term containing one $O(\tau^3)$ column has order at least $3+2+2=7$.  Since $-\mathfrak p=T+O(\tau)$, this proves \eqref{eqn:scaled-Evans-sixth-order-low}, including the analyticity of the remainder.  In particular, no monodromy coefficient of total degree three or higher contributes to its sixth-order term.
\end{proof}

We now compute the determinant in \eqref{eqn:scaled-Evans-sixth-order-low}.  Substituting \eqref{eqn:M1-M2-low} into \eqref{eqn:effective-second-order-matrix-low}, and using \eqref{eqn:relations-amn} and the coefficient identities in Appendix~\ref{app:remaining-low-coefficient-relations}, gives the following entries of $S_2$:
\begin{subequations}\label{eqn:explicit-S2-low}
\begin{align}
 (S_2)_{11}={}&-a_{11}^{(1,0)}\zeta +c_g^2m^2\left[ a_{11}^{(2,0)}-\frac{(a_{11}^{(1,0)})^2}{2} -\frac{a_{13}^{(1,0)}a_{41}^{(1,0)}}{T}\right] \notag\\
 &+\alpha^2\left[ \frac{a_{13}^{(0,1)}a_{41}^{(0,1)}}{T}-a_{11}^{(0,2)}\right], \\
 (S_2)_{22}={}&-a_{11}^{(1,0)}\zeta -c_g^2m^2\left[ a_{11}^{(2,0)}-\frac{(a_{11}^{(1,0)})^2}{2} -\frac{a_{13}^{(1,0)}a_{41}^{(1,0)}}{T}\right] \notag\\
 &-\alpha^2\left[ \frac{a_{13}^{(0,1)}a_{41}^{(0,1)}}{T}-a_{11}^{(0,2)}\right], \\
 (S_2)_{12}={}&-(S_2)_{21} =\alpha^2\left[ a_{11}^{(0,2)}-\frac{a_{13}^{(0,1)}a_{41}^{(0,1)}}{T}\right], \\
 (S_2)_{13}={}&(S_2)_{23} =-ic_gm\alpha\left[ a_{14}^{(1,1)}+ \frac{a_{13}^{(0,1)}(a_{11}^{(1,0)}-a_{44}^{(1,0)})}{T} \right], \\
 (S_2)_{31}={}&-ic_gm\alpha\left[ a_{31}^{(1,1)}+ \frac{a_{41}^{(0,1)}(a_{11}^{(1,0)}-a_{33}^{(1,0)})}{T} \right]-\alpha^2a_{31}^{(0,2)}, \\
 (S_2)_{32}={}&ic_gm\alpha\left[ a_{31}^{(1,1)}+ \frac{a_{41}^{(0,1)}(a_{11}^{(1,0)}-a_{33}^{(1,0)})}{T} \right]-\alpha^2a_{32}^{(0,2)}, \\
 (S_2)_{33}={}&\frac{c_g^2m^2}{T}\left[ Ta_{34}^{(2,0)} -(a_{11}^{(1,0)}-a_{33}^{(1,0)}) (a_{11}^{(1,0)}-a_{44}^{(1,0)})\right].
\end{align}
\end{subequations}

Three combinations of coefficients occur in the determinant.  We write
\begin{align}
 F_{\mathrm{low}}
 &:=(a_{11}^{(1,0)})^2-2a_{11}^{(2,0)} +\frac{2a_{13}^{(1,0)}a_{41}^{(1,0)}}{T}, \notag\\
 Q_L
 &:=Ta_{34}^{(2,0)} -(a_{11}^{(1,0)}-a_{33}^{(1,0)}) (a_{11}^{(1,0)}-a_{44}^{(1,0)}),
 \label{eqn:def-FQ-low}
\end{align}
and
\begin{align}
 G_{\mathrm{low}}:={}&-a_{11}^{(0,2)}Q_L+a_{31}^{(1,1)} \left[Ta_{14}^{(1,1)}+a_{13}^{(0,1)} (a_{11}^{(1,0)}-a_{44}^{(1,0)})\right] \notag\\
 &+a_{41}^{(0,1)} \left[a_{14}^{(1,1)}(a_{11}^{(1,0)}-a_{33}^{(1,0)}) +a_{13}^{(0,1)}a_{34}^{(2,0)}\right].
 \label{eqn:def-G-low}
\end{align}
By \eqref{eqn:Evans-amplitude-parity-low} and Lemma~\ref{lem:second-order-sufficiency-low}, $T\det S_2$ is even in $\alpha$.  We may therefore discard its odd part when expanding \eqref{eqn:explicit-S2-low}.  The terms containing $a_{31}^{(0,2)}$ or $a_{32}^{(0,2)}$ occur only in this odd part. The remaining terms give
\begin{align}
 T\det S_2={}&c_g^2Q_L(a_{11}^{(1,0)})^2m^2\zeta^2 -\frac{c_g^6Q_LF_{\mathrm{low}}^2}{4}m^6+c_g^4F_{\mathrm{low}}G_{\mathrm{low}}m^4\alpha^2.
 \label{eqn:S2-determinant-low}
\end{align}
Thus only the coefficients listed in \eqref{eqn:second-order-coefficients-needed}, together with the first-order coefficients, enter the determinant.

Their explicit expressions are recorded in Appendix~\ref{app:all-explicit-low-order-coefficients}. Substituting them into \eqref{eqn:def-FQ-low} and using $\mu_0^2+\omega\mu_0=\kappa\coth\kappa$, we find
\begin{equation}\label{eqn:simplified-FQ-low}
 a_{11}^{(1,0)}=\frac{T}{c_g}, \qquad Q_L=\frac{T^2\Delta_L}{c_g^2(\kappa\coth\kappa-1)}, \qquad F_{\mathrm{low}}=-\frac{2iTc_2}{c_g^3}.
\end{equation}
In particular, $Q_L<0$ by Lemma~\ref{lem:automatic-low-nonresonance} and $\kappa\coth\kappa>1$.  Inserting \eqref{eqn:simplified-FQ-low} into \eqref{eqn:S2-determinant-low} yields
\begin{equation*} T\det S_2 =T^2Q_Lm^2(\zeta^2+c_2^2m^4) +c_g^4F_{\mathrm{low}}G_{\mathrm{low}}m^4\alpha^2.
\end{equation*}
We define
\begin{equation}\label{eqn:index-from-amn}
 \mathrm{ind}_{\mathrm{low}}(\omega,\kappa) :=-\frac{c_g^4F_{\mathrm{low}}G_{\mathrm{low}}}{T^2Q_L} =-\frac{T^2F_{\mathrm{low}}G_{\mathrm{low}}} {(a_{11}^{(1,0)})^4Q_L}.
\end{equation}
Then
\begin{equation}\label{eqn:carrier-quadratic-leading-low}
 T\det S_2 =T^2Q_Lm^2 \bigl(\zeta^2+c_2^2m^4 -\mathrm{ind}_{\mathrm{low}}m^2\alpha^2\bigr).
\end{equation}
The index is real: the coefficient formulas show that $F_{\mathrm{low}}$ and $G_{\mathrm{low}}$ are purely imaginary, whereas $Q_L$ is real.  Its closed expression in $\omega$ and $\kappa$ is given in Appendix~\ref{index:low}. Together with Lemma~\ref{lem:second-order-sufficiency-low}, \eqref{eqn:carrier-quadratic-leading-low} determines the leading term of the exact Evans function under \eqref{eqn:carrier-scaling}.

We use this quadratic to locate the unstable spectral roots.

\begin{theorem}[Benjamin--Feir instability]\label{thm:BF}
Let $\omega\neq0$ and $\kappa>0$.  If
\begin{equation*} \mathrm{ind}_{\mathrm{low}}(\omega,\kappa)>0,
\end{equation*}
then every sufficiently small, nonzero Stokes wave is spectrally unstable to long-wave perturbations.  More precisely, fix a nonempty compact interval $J\subset(0,\infty)$ such that
\begin{equation}\label{eqn:unstable-ray-interval-low}
 \inf_{m\in J} \bigl(\mathrm{ind}_{\mathrm{low}}-c_2^2m^2\bigr)>0.
\end{equation}
There is $\epsilon_J>0$ such that, for $m\in J$ and $0<|\epsilon|<\epsilon_J$, the spectrum at Bloch exponent $\gamma=m\epsilon$ contains two points satisfying
\begin{equation}\label{eqn:carrier-roots-low}
 \delta_\pm(m\epsilon,\epsilon) =ic_gm\epsilon \pm m\epsilon^2 \sqrt{\mathrm{ind}_{\mathrm{low}}-c_2^2m^2} +O_J(|\epsilon|^3).
\end{equation}
Their real parts have opposite signs.  In particular, these unstable spectral points tend to the origin as $\epsilon\to0$.
\end{theorem}

\begin{proof}
Set $\alpha=1$ and $\tau=\epsilon$ in \eqref{eqn:carrier-scaling}.  For $m$ in a neighborhood of $J$, define
\begin{equation*} \mathcal E(\zeta,m,\epsilon) =\frac{D(ic_gm\epsilon+\epsilon^2\zeta,m\epsilon,\epsilon)} {\epsilon^6T^2Q_Lm^2}, \qquad \epsilon\neq0.
\end{equation*}
The analytic factorization \eqref{eqn:scaled-Evans-sixth-order-low} shows that this function extends analytically to $\epsilon=0$.  Its value there is
\[
 \mathcal E(\zeta,m,0) =\zeta^2-m^2 \bigl(\mathrm{ind}_{\mathrm{low}}-c_2^2m^2\bigr).
\]
For each $m\in J$, its two roots are $\zeta_\pm^0(m)=\pm m\sqrt{\mathrm{ind}_{\mathrm{low}}-c_2^2m^2}$. They are real, nonzero, and simple, since $\partial_\zeta\mathcal E(\zeta_\pm^0(m),m,0)=2\zeta_\pm^0(m)\neq0$. The analytic implicit function theorem therefore gives unique local roots $\zeta_\pm(m,\epsilon)$ near $\zeta_\pm^0(m)$, analytic in $(m,\epsilon)$.  The local roots agree on overlaps by uniqueness.  Compactness of $J$ and \eqref{eqn:unstable-ray-interval-low} give a common neighborhood and the uniform estimate $\zeta_\pm(m,\epsilon)=\zeta_\pm^0(m)+O_J(|\epsilon|)$. Returning to $\delta$ proves \eqref{eqn:carrier-roots-low}. For real parameters, $ic_gm\epsilon$ is purely imaginary, while the real parts of $\zeta_+$ and $\zeta_-$ remain respectively positive and negative for small $|\epsilon|$.  These are zeros of the exact Evans function.  Since $\gamma=m\epsilon$ is real, Lemma~\ref{lem:Evans-spectrum} identifies them as temporal spectral points. This proves the assertion.
\end{proof}

For $\epsilon>0$, the leading terms in \eqref{eqn:carrier-roots-low}, together with their complex conjugates, give the curve
\begin{equation*} \left(\pm m\sqrt{\mathrm{ind}_{\mathrm{low}}-c_2^2m^2},\ c_gm\right), \qquad |m|\leq\frac{\sqrt{\mathrm{ind}_{\mathrm{low}}}}{|c_2|},
\end{equation*}
in the coordinates $(\operatorname{Re}\lambda/\epsilon^2, \operatorname{Im}\lambda/\epsilon)$. This is the figure-eight curve in Figure~\ref{fig:leading-spectra}(a). Theorem~\ref{thm:BF} justifies the expansion away from its crossing and endpoints, where the leading roots cease to be simple.

\section{The spectrum away from the origin}
\label{sec:high-freq}

We study the spectrum near the order-two resonance \eqref{eqn:resonance}.  Throughout this section, $\kappa>0$ and $\omega\in\mathbb R\setminus\{0\}$ are fixed, and $T=2\pi/\kappa$.  The Stokes amplitude has the normalization used in Section~\ref{sec:Stokes}.  The two resonant modes have the same Floquet multiplier. We calculate how this double multiplier splits when the amplitude becomes nonzero.

\subsection{The resonant frequency}
\label{sec:high-resonance}

Recall that $k_2(\sigma)>\kappa$ and $k_4(\sigma)>0$ are determined by
\[
 \sigma_-(k_2(\sigma))=\sigma, \qquad \sigma_+(k_4(\sigma))=\sigma, \qquad \sigma>0.
\]
Their continuous extensions at zero satisfy $k_2(0)=\kappa$ and $k_4(0)=0$.

\begin{lemma}\label{lem:high-resonance}
There is a unique $\sigma_H>0$ for which $k_2(\sigma_H)-k_4(\sigma_H)=2\kappa$. Moreover, $\sigma_H>\sigma_c$, and
\begin{equation}\label{eqn:high-slope-order}
 0<\sigma_-'(k_2(\sigma_H))<1 <\sigma_+'(k_4(\sigma_H)).
\end{equation}
\end{lemma}

\begin{proof}
Recall that $\rho=\omega\mu_0$.  For $k>0$, let
\[
 Q(k)=\mu_0^2k\tanh k+\frac{\rho^2\tanh^2k}{4}.
\]
Then $Q'>0$, and \eqref{eqn:sigmaplus-positive-identity} gives
\[
 (Q')^2-\rho^2\sech^4k\,Q =\mu_0^4(\tanh k+k\sech^2k)^2 +\mu_0^2\rho^2\tanh^2k\,\sech^2k>0.
\]
It follows that
\[
 \sigma_+'(k)=1-\frac{\rho\sech^2k}{2} +\frac{Q'}{2\sqrt Q}>1.
\]
By Lemma~\ref{lem:sigma-geometry}, $\sigma_-$ is strictly convex, vanishes at $0$ and $\kappa$, and is increasing for $k\geq\kappa$. Its explicit formula gives $\sigma_-'(k)\to1$ as $k\to\infty$. It follows that $0<\sigma_-'(k)<1$ for $k\geq\kappa$.  Thus
\begin{equation*} \frac{d}{d\sigma}(k_2-k_4) =\frac1{\sigma_-'(k_2)}-\frac1{\sigma_+'(k_4)}>0.
\end{equation*}
The gap equals $\kappa$ at $\sigma=0$.  As $k\to\infty$,
\[
 \sigma_\pm(k)=k\pm\mu_0\sqrt k-\frac{\rho}{2} +O(k^{-1/2}),
\]
and hence
\[
 k_2(\sigma)-k_4(\sigma)=2\mu_0\sqrt\sigma+O(1) \longrightarrow\infty.
\]
This proves existence and uniqueness of $\sigma_H$.

Let $r_c\in(0,\kappa)$ be the minimum point of $\sigma_-$, so that $\sigma_c=-\sigma_-(r_c)$.  Since $\sigma_-'$ is strictly increasing,
\[
 \sigma_c=\int_{r_c}^{\kappa}\sigma_-'(s)\,ds <\int_{\kappa}^{2\kappa-r_c}\sigma_-'(s)\,ds =\sigma_-(2\kappa-r_c).
\]
At the resonance, $k_2=2\kappa+k_4>2\kappa$.  Therefore $\sigma_H=\sigma_-(k_2)>\sigma_-(2\kappa-r_c)>\sigma_c$. The slope inequalities were established above.
\end{proof}

In the remainder of this section, the argument $\sigma_H$ is suppressed in $k_j$, $\Upphi_j$, and $\uppsi_j$.  In particular,
\begin{equation*} \rho_0=e^{ik_2T}=e^{ik_4T},\qquad |\rho_0|=1.
\end{equation*}
Lemma~\ref{lem:high-resonance} places the resonance in the region where $\mathbf L(i\sigma_H)$ has exactly two eigenvalues on the imaginary axis.  The reduction of Section~\ref{sec:spec} therefore has dimension two, with ordered basis $(\Upphi_2,\Upphi_4)$.

\subsection{Expansion of the monodromy matrix}
\label{sec:aH}

Let $\delta=\lambda-i\sigma_H$, and write
\[
 \mathbf L_H=\mathbf L(i\sigma_H),\qquad \mathbf B_H(x;\delta,\epsilon) =\widetilde{\mathbf B}(x;\sigma_H,\delta,\epsilon).
\]
The coefficients $\mathbf B_H^{(m,n)}$ are defined by
\[
 \mathbf B_H(x;\delta,\epsilon) =\sum_{m+n\geq1}\delta^m\epsilon^n \mathbf B_H^{(m,n)}(x).
\]
As in Section~\ref{sec:a0}, superscripts denote Taylor coefficients without factorials.  Every forcing term $f\in Y$ is projected by
\begin{equation*} [P_Hf]_{(\Upphi_2,\Upphi_4)}
 =\begin{pmatrix}
 \langle f,\uppsi_2\rangle_Y\\
 \langle f,\uppsi_4\rangle_Y
 \end{pmatrix}.
\end{equation*}

Let $\mathbf a(x;\delta,\epsilon)$ be the principal fundamental matrix of the reduced system in this basis, so that $\mathbf a(0;\delta,\epsilon)=I_2$. We use the analytic expansion
\begin{equation}\label{eqn:aexp}
 \mathbf a(x;\delta,\epsilon) =\sum_{m,n\geq0}\mathbf a^{(m,n)}(x)\delta^m\epsilon^n, \qquad \mathbf a^{(0,0)}(x)=\operatorname{diag}(e^{ik_2x},e^{ik_4x}).
\end{equation}
Below, $\mathbf a^{(m,n)}$ without an argument means $\mathbf a^{(m,n)}(T)$.

\begin{lemma}\label{lem:high-monodromy-structure}
The first coefficients in \eqref{eqn:aexp} satisfy
\[
 \mathbf a^{(0,0)}=\rho_0I_2,\qquad \mathbf a^{(1,0)}=\rho_0
 \begin{pmatrix}\alpha_1&0\\0&\alpha_2\end{pmatrix},
 \qquad \mathbf a^{(0,1)}=0,
\]
where
\begin{equation}\label{eqn:high-alpha}
 \alpha_1=\frac{T}{\sigma_-'(k_2)},\qquad \alpha_2=\frac{T}{\sigma_+'(k_4)},\qquad \alpha_1>T>\alpha_2>0.
\end{equation}
Moreover,
\begin{equation}\label{eqn:high-beta}
 \mathbf a^{(0,2)}=i\rho_0\boldsymbol\beta, \qquad \boldsymbol\beta=(\beta_{j\ell})_{j,\ell=1}^2 \in\mathbb R^{2\times2}.
\end{equation}
In particular, $a_{11}^{(1,0)}a_{22}^{(1,0)}\ne0$ and $a_{11}^{(1,0)}-a_{22}^{(1,0)}\ne0$.
\end{lemma}

\begin{proof}
At zero amplitude, the two simple spatial eigenvalues continue analytically as $ik_j(\sigma_H-i\delta)$, $j=2,4$.  Their multipliers have derivatives
\[
 \left.\frac{d}{d\delta}e^{ik_j(\sigma_H-i\delta)T} \right|_{\delta=0} =\rho_0T k_j'(\sigma_H).
\]
To express this in the fixed center basis, diagonalize the flat reduced matrix with an analytic change of basis $S(\delta)$, normalized by $S(0)=I_2$.  At $\delta=0$ both multipliers equal $\rho_0$; the two terms containing $S'(0)$ in the derivative of the conjugated monodromy therefore cancel.  This proves the formula for $\mathbf a^{(1,0)}$. The inequalities in \eqref{eqn:high-alpha} follow from \eqref{eqn:high-slope-order}.

The first-order Stokes profiles and the absence of a linear correction to $\mu$ imply
\begin{equation*} \mathbf B_H^{(0,1)}(x) =\mathbf B_{H,1}^{(0,1)}e^{i\kappa x} +\mathbf B_{H,-1}^{(0,1)}e^{-i\kappa x}.
\end{equation*}
The first-order variation of constants gives
\[
 a_{\ell j}^{(0,1)}(T) =\rho_0\int_0^T e^{i(k_{2j}-k_{2\ell})x} \langle\mathbf B_H^{(0,1)}(x)\Upphi_{2j}, \uppsi_{2\ell}\rangle_Y\,dx.
\]
Here $k_{2j}-k_{2\ell}\in\{0,2\kappa,-2\kappa\}$, so every exponential in the integral has a nonzero integer multiple of $\kappa$ as its frequency.  Thus $\mathbf a^{(0,1)}(T)=0$.

It remains to prove \eqref{eqn:high-beta}.  Define the antilinear map
\[
 \mathscr S(\phi,\chi,\eta) =(\overline\phi,-\overline\chi,-\overline\eta).
\]
The reality and reflection identities used in the proof of Lemma~\ref{lem:sym} imply that, at $\delta=0$ and real $\epsilon$, $u(x)\mapsto\mathscr S u(-x)$ preserves solutions of the full spatial system.  Each $\Upphi_{2j}$ is fixed by $\mathscr S$.  The spectral projection $P_H$ commutes with this map, and uniqueness of the periodic invariant graph in Theorem~\ref{thm:Mielke-specialized} implies that the graph also respects this symmetry.  Hence, in the chosen center coordinates, $v(x)\mapsto\overline{v(-x)}$ preserves the reduced equation.  Uniqueness of its principal fundamental matrix gives
\[
 \overline{\mathbf a(-x;0,\epsilon)}=\mathbf a(x;0,\epsilon).
\]
Periodicity gives $\mathbf a(-T;0,\epsilon)=\mathbf a(T;0,\epsilon)^{-1}$, and therefore
\begin{equation*} \overline{\mathbf a(T;0,\epsilon)} =\mathbf a(T;0,\epsilon)^{-1}.
\end{equation*}
Using $\mathbf a^{(0,1)}=0$ and comparing the coefficients of $\epsilon^2$, we obtain
\[
 \overline{\mathbf a^{(0,2)}}=-\rho_0^{-2}\mathbf a^{(0,2)}.
\]
Since $|\rho_0|=1$, this is precisely \eqref{eqn:high-beta}.
\end{proof}

We next give a formula for all four second-order coefficients.  For $n=1,2$, let
\begin{equation*} \mathbf B_{H,r}^{(0,n)} =\frac1T\int_0^T e^{-ir\kappa x}\mathbf B_H^{(0,n)}(x)\,dx, \qquad r\in\mathbb Z.
\end{equation*}
Only $r=\pm1$ occur when $n=1$, and only $r=0,\pm2$ occur when $n=2$. For $j=1,2$ and $s\in\{-1,1\}$, define
\begin{equation*} R_{j,s}=\bigl(i(k_{2j}+s\kappa)I-\mathbf L_H\bigr)^{-1}.
\end{equation*}
These inverses exist: the only imaginary eigenvalues of $\mathbf L_H$ are $ik_2$ and $ik_4$, whose wave numbers differ by $2\kappa$. They act from $Y$ to $\mathcal D(\mathbf L_H)$, so the operator products below are well defined.

\begin{lemma}\label{lem:high-second-coefficients}
For $j,\ell\in\{1,2\}$,
\begin{align}
 a_{\ell j}^{(0,2)} =\rho_0T\Bigg\langle
 &\mathbf B_{H,2(j-\ell)}^{(0,2)}\Upphi_{2j} +\sum_{s=\pm1}\mathbf B_{H,2(j-\ell)-s}^{(0,1)} R_{j,s}\mathbf B_{H,s}^{(0,1)}\Upphi_{2j}, \uppsi_{2\ell}\Bigg\rangle_Y.
 \label{eqn:high-second-coefficients}
\end{align}
The first-order Fourier coefficients outside $\{-1,1\}$ are zero. In particular, with
\begin{equation}\label{eqn:high-middle-resolvent}
 q_* = k_4+\kappa=k_2-\kappa,\qquad R_*=(iq_*I-\mathbf L_H)^{-1},
\end{equation}
the two coupling coefficients are
\begin{align}
 a_{12}^{(0,2)}
 &=\rho_0T\Big\langle \bigl(\mathbf B_{H,2}^{(0,2)} +\mathbf B_{H,1}^{(0,1)}R_*\mathbf B_{H,1}^{(0,1)}\bigr) \Upphi_4,\uppsi_2\Big\rangle_Y,
 \label{eqn:high-a12-explicit}\\
 a_{21}^{(0,2)}
 &=\rho_0T\Big\langle \bigl(\mathbf B_{H,-2}^{(0,2)} +\mathbf B_{H,-1}^{(0,1)}R_*\mathbf B_{H,-1}^{(0,1)}\bigr) \Upphi_2,\uppsi_4\Big\rangle_Y.
 \label{eqn:high-a21-explicit}
\end{align}
\end{lemma}

\begin{proof}
Let $u_j(x;0,\epsilon)$ be the full solution on the invariant graph whose center coordinates at $x=0$ are the $j$th coordinate vector. Its zeroth-order term is $e^{ik_{2j}x}\Upphi_{2j}$.  The first-order term solves
\[
 (\partial_x-\mathbf L_H)u_j^{(0,1)} =\sum_{s=\pm1}e^{i(k_{2j}+s\kappa)x} \mathbf B_{H,s}^{(0,1)}\Upphi_{2j}, \qquad P_Hu_j^{(0,1)}(0)=0.
\]
A bounded solution satisfying this initial condition is
\[
 u_j^{(0,1)}(x) =\sum_{s=\pm1}e^{i(k_{2j}+s\kappa)x} R_{j,s}\mathbf B_{H,s}^{(0,1)}\Upphi_{2j}-\sum_{h=1}^2e^{ik_{2h}x}\Upphi_{2h} \sum_{s=\pm1} \big\langle R_{j,s}\mathbf B_{H,s}^{(0,1)}\Upphi_{2j}, \uppsi_{2h}\big\rangle_Y.
\]
The complementary resolvent estimate and the center initial condition give uniqueness in the class of bounded solutions, so this is the coefficient selected by the invariant graph.

At second order, variation of constants in the center equation gives
\[
 a_{\ell j}^{(0,2)}(T)=\rho_0\int_0^T e^{-ik_{2\ell} x} \Big\langle\mathbf B_H^{(0,2)}(x)e^{ik_{2j}x}\Upphi_{2j}+\mathbf B_H^{(0,1)}(x)u_j^{(0,1)}(x), \uppsi_{2\ell}\Big\rangle_Y\,dx.
\]
The second line in the formula for $u_j^{(0,1)}$ makes no contribution: after multiplication by $\mathbf B_H^{(0,1)}$, each such term has frequency $k_{2h}-k_{2\ell}\pm\kappa\ne0$ in the integral.  The remaining terms survive precisely when their total frequency is zero.  This gives \eqref{eqn:high-second-coefficients}.  For the two off-diagonal entries, the Fourier shifts are $2\kappa$ and $-2\kappa$, respectively. Thus only the common intermediate wave number $q_*$ occurs.
\end{proof}

There are two contributions to each coupling coefficient.  The second harmonic of $\mathbf B_H^{(0,2)}$ connects the two modes in one step.  The first harmonic of $\mathbf B_H^{(0,1)}$ connects them in two steps, through $q_*=k_4+\kappa=k_2-\kappa$. The resolvent $R_*$ calculates this intermediate response. Thus the second-order Stokes profiles and the first-order correction already determine the leading coupling.  These are the same terms obtained by the recursive calculation in Section~\ref{sec:a0}.

\subsection{The high-frequency instability index}
\label{sec:ind2}

For a real Floquet displacement $\gamma$, define
\begin{equation*} D_{\mathrm{high}}(\delta,\gamma,\epsilon) =\det\bigl(\mathbf a(T;\delta,\epsilon) -\rho_0e^{i\gamma T}I_2\bigr).
\end{equation*}
The Floquet exponent is $k_4+\gamma$, understood modulo $\kappa$. This is the Evans function of Definition~\ref{def:Delta} at $\lambda=i\sigma_H+\delta$.

The product of the two coupling coefficients determines the leading splitting.  We define
\begin{equation}\label{def:ind2}
 \mathrm{ind}_{\mathrm{high}}(\omega,\kappa) =\frac{a_{12}^{(0,2)}a_{21}^{(0,2)}} {a_{11}^{(1,0)}a_{22}^{(1,0)}} =-\frac{\beta_{12}\beta_{21}}{\alpha_1\alpha_2}.
\end{equation}
The coefficients in this formula are given by \eqref{eqn:high-alpha} and \eqref{eqn:high-second-coefficients}. The index is real by Lemma~\ref{lem:high-monodromy-structure}. In fact, the two couplings have opposite signs whenever they are nonzero.

\begin{proposition}\label{prop:high-square}
For every $\kappa>0$ and $\omega\in\mathbb R\setminus\{0\}$,
\begin{equation}\label{eqn:high-beta-reciprocity}
 \beta_{21}=-R_H\beta_{12},
\end{equation}
where
\begin{equation}
 R_H=\frac{\alpha_2}{\alpha_1} \frac{k_2\sinh(2k_2)}{k_4\sinh(2k_4)} \frac{(k_4-\sigma_H)^2}{(k_2-\sigma_H)^2} \frac{2(k_2-\sigma_H)-\rho\tanh k_2} {\rho\tanh k_4-2(k_4-\sigma_H)}>0.
 \label{eqn:high-reciprocity-factor}
\end{equation}
Consequently,
\begin{equation}\label{eqn:high-square}
 \mathrm{ind}_{\mathrm{high}}(\omega,\kappa) =\frac{R_H}{\alpha_1\alpha_2}\,\beta_{12}^2\geq0, \qquad \mathrm{ind}_{\mathrm{high}}\ne0 \quad\Longleftrightarrow\quad\beta_{12}\ne0.
\end{equation}
\end{proposition}

\begin{proof}
Substitution of the Stokes coefficients into \eqref{eqn:high-a12-explicit}--\eqref{eqn:high-a21-explicit}, followed by the endpoint dispersion relations, gives \eqref{eqn:high-beta-reciprocity}.  The algebra is given in Appendix~\ref{A:high-index}; no division by either coupling is needed.

We check the sign of the factor in \eqref{eqn:high-reciprocity-factor}.  Both wave numbers and both $\alpha_j$ are positive.  Also, $k_j-\sigma_H\ne0$ for $j=2,4$: otherwise the dispersion relation would imply $\mu_0^2k_j\tanh k_j=0$.  Finally, $\sigma_H=\sigma_-(k_2)=\sigma_+(k_4)$ gives
\begin{align*}
 2(k_2-\sigma_H)-\rho\tanh k_2
 &=\sqrt{\rho^2\tanh^2k_2+4\mu_0^2k_2\tanh k_2}>0,\\
 \rho\tanh k_4-2(k_4-\sigma_H)
 &=\sqrt{\rho^2\tanh^2k_4+4\mu_0^2k_4\tanh k_4}>0.
\end{align*}
Every factor in $R_H$ is therefore positive.  Substitution into \eqref{def:ind2} proves \eqref{eqn:high-square}.
\end{proof}

To describe the spectral curves, let
\[
 m_c=\frac{\alpha_2\beta_{11}-\alpha_1\beta_{22}} {T(\alpha_2-\alpha_1)},\qquad v_H(m)=\frac12\left( \frac{Tm-\beta_{11}}{\alpha_1} +\frac{Tm-\beta_{22}}{\alpha_2}\right).
\]
The value $m_c$ is where the two diagonal approximations coincide; $v_H(m)$ is their average.  These two quantities give the center of the splitting in the Floquet and spectral variables, respectively.

\begin{theorem}\label{thm:high}
Suppose that $\beta_{12}\ne0$, and let $J$ be a compact subinterval of
\begin{equation}\label{eqn:high-bubble-interval}
 \left( m_c-\frac{2\alpha_1\alpha_2\sqrt{\mathrm{ind}_{\mathrm{high}}}} {T(\alpha_1-\alpha_2)},\quad m_c+\frac{2\alpha_1\alpha_2\sqrt{\mathrm{ind}_{\mathrm{high}}}} {T(\alpha_1-\alpha_2)} \right).
\end{equation}
For all sufficiently small nonzero real $\epsilon$ and all $m\in J$, the temporal spectrum at Floquet displacement $\gamma=\epsilon^2m$ contains two roots satisfying
\begin{align}
 \lambda_\pm(m,\epsilon) ={}&i\sigma_H+i\epsilon^2v_H(m) \notag\\
 &{}\pm\epsilon^2\sqrt{ \mathrm{ind}_{\mathrm{high}} -\frac{T^2}{4}\left(\frac1{\alpha_1}-\frac1{\alpha_2}\right)^2 (m-m_c)^2} +O(|\epsilon|^3).
 \label{eqn:high-root-expansion}
\end{align}
The remainder is uniform for $m\in J$.  In particular, every sufficiently small nonzero Stokes wave is spectrally unstable when $\beta_{12}\ne0$.
\end{theorem}

\begin{proof}
Lemma~\ref{lem:high-monodromy-structure} gives
\[
 \mathbf a(T;\delta,\epsilon) =\rho_0\bigl(I_2+\delta\operatorname{diag}(\alpha_1,\alpha_2) +i\epsilon^2\boldsymbol\beta\bigr) +O(|\delta|^2+|\delta\epsilon|+|\epsilon|^3).
\]
The linear terms in $\delta,\gamma$ balance the quadratic term in $\epsilon$ when
\begin{equation*} \delta=\epsilon^2z,\qquad \gamma=\epsilon^2m.
\end{equation*}
Under this scaling,
\[\rho_0^{-1}\bigl(\mathbf a(T;\epsilon^2z,\epsilon) -\rho_0e^{iT\epsilon^2m}I_2\bigr)=\epsilon^2\left[
 \begin{pmatrix}
 \alpha_1z+i(\beta_{11}-Tm)&i\beta_{12}\\
 i\beta_{21}&\alpha_2z+i(\beta_{22}-Tm)
 \end{pmatrix}
 +\epsilon\mathbf E(z,m,\epsilon)\right],
 \]
where $\mathbf E$ is analytic near every finite $(z,m,0)$. Each column is divisible by $\epsilon^2$.  Hence the determinant is divisible by $\epsilon^4$, and its quotient extends analytically to $\epsilon=0$.  Taking the determinant and completing the square gives
\[ \frac{D_{\mathrm{high}}(\epsilon^2z,\epsilon^2m,\epsilon)} {\rho_0^2\alpha_1\alpha_2\epsilon^4} =(z-iv_H(m))^2-\mathrm{ind}_{\mathrm{high}}+\frac{T^2}{4} \left(\frac1{\alpha_1}-\frac1{\alpha_2}\right)^2(m-m_c)^2 +O(\epsilon).
 \]

For $m$ in \eqref{eqn:high-bubble-interval}, the polynomial at $\epsilon=0$ has two distinct roots.  Its $z$-derivative is nonzero at both roots, so the analytic implicit-function theorem continues each root to small $\epsilon$, with an $O(\epsilon)$ correction. Uniqueness identifies these local branches on overlapping neighborhoods.  Compactness of $J$ gives a common smallness bound and a uniform remainder.  Returning to $\lambda=i\sigma_H+\epsilon^2z$ proves \eqref{eqn:high-root-expansion}.

The square root in \eqref{eqn:high-root-expansion} is bounded away from zero on $J$.  Thus one root has positive real part and the other negative real part for small nonzero $\epsilon$.  The Floquet exponent $k_4+\epsilon^2m$ is real, so Lemma~\ref{lem:Evans-spectrum} places both roots in the temporal spectrum.
\end{proof}

At $m=m_c$ the real splitting is largest, and
\begin{equation*} \lambda_\pm(m_c,\epsilon) =i\sigma_H+i\epsilon^2 \frac{\beta_{11}-\beta_{22}}{\alpha_2-\alpha_1} \pm\epsilon^2\sqrt{\mathrm{ind}_{\mathrm{high}}} +O(|\epsilon|^3).
\end{equation*}
This explains the normalization in \eqref{def:ind2}: the square root of the index gives the leading growth rate at the center.

The same expansion describes the shape of the spectral bubble. Since
\[
 v_H(m)-v_H(m_c) =\frac T2\left(\frac1{\alpha_1}+\frac1{\alpha_2}\right)(m-m_c),
\]
the limiting curves for $z=(\lambda-i\sigma_H)/\epsilon^2$ lie on
\begin{equation}\label{eqn:high-leading-ellipse}
 (\operatorname{Re}z)^2 +\left(\frac{\alpha_1-\alpha_2}{\alpha_1+\alpha_2}\right)^2 \bigl(\operatorname{Im}z-v_H(m_c)\bigr)^2 =\mathrm{ind}_{\mathrm{high}}.
\end{equation}
Figure~\ref{fig:leading-spectra}(b) illustrates this ellipse. The root expansion is uniform on compact subsets of the open interval \eqref{eqn:high-bubble-interval}.  At the endpoints the limiting roots coalesce, and the simple-root argument does not apply.

\begin{remark}\label{rem:high-index-degeneracy}
The two coupling coefficients vanish simultaneously by \eqref{eqn:high-beta-reciprocity}.  At such a parameter pair, $\mathrm{ind}_{\mathrm{high}}=0$, and the limiting quadratic has a double root at $m=m_c$.  Higher-order terms are then needed to decide whether this resonance produces instability.  Either coupling can be used to locate these exceptional parameters.  Its sign may change across a zero, while the index stays positive wherever the coupling is nonzero. Thus the zero set of the coupling alone does not establish a boundary between stability and instability.
\end{remark}

\section{Illustration of the instability criteria}
\label{sec:comparison-indices}

\begin{figure}[!htbp]
 \centering
 \includegraphics[width=\textwidth]{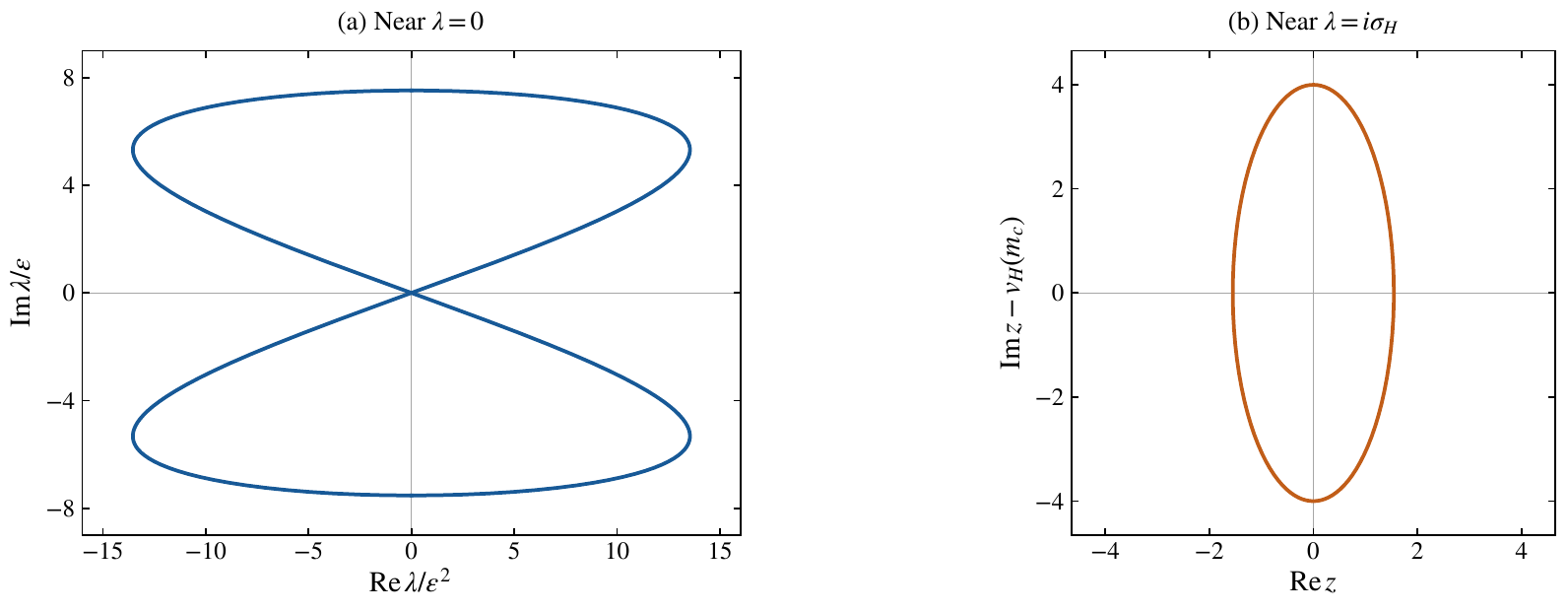}
\caption{Leading spectral curves for $(\omega,\kappa)=(1/2,2)$. (a) The low-frequency figure-eight in the coordinates $(\operatorname{Re}\lambda/\epsilon^2, \operatorname{Im}\lambda/\epsilon)$, for $\epsilon>0$. (b) The high-frequency ellipse in the coordinates $(\operatorname{Re}z,\operatorname{Im}z-v_H(m_c))$, where $z=(\lambda-i\sigma_H)/\epsilon^2$ and $\sigma_H\approx2.061971$. These are the curves determined by the leading terms of the expansions, rather than a numerical computation of the full spectrum at a fixed nonzero amplitude.}
 \label{fig:leading-spectra}
\end{figure}

Figure~\ref{fig:leading-spectra} shows the leading spectral curves at $(\omega,\kappa)=(1/2,2)$, where the low-frequency index is positive and the high-frequency index is nonzero.  The figure-eight curve is obtained from \eqref{eqn:carrier-roots-low}, and the ellipse from \eqref{eqn:high-leading-ellipse}.  The latter is translated to its leading center.  The scales on the axes reflect the different orders of the real and imaginary parts in the two expansions.

\begin{figure}[htbp]
 \centering
 \includegraphics[width=\textwidth]{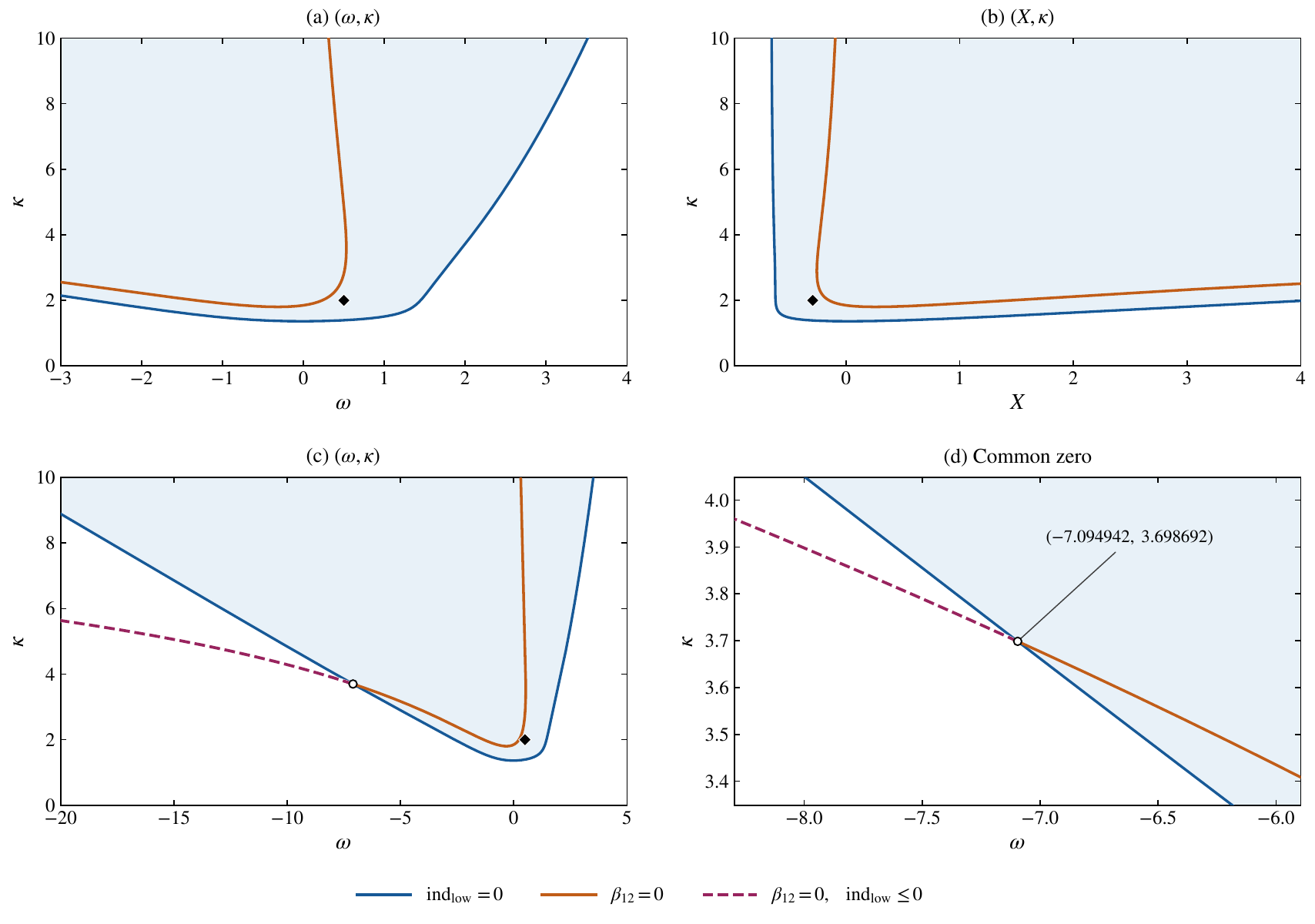}
\caption{Comparison of the two instability criteria. Blue shading denotes $\mathrm{ind}_{\mathrm{low}}>0$, with its zero set in blue.  The orange and magenta curves together give $\beta_{12}=0$; the magenta dashed part satisfies \eqref{eqn:undecided-leading-indices}. Panels~(a)--(c) show two parameter ranges and the coordinate $X$; panel~(d) enlarges the common zero (open circle). Diamonds mark the parameters in Figure~\ref{fig:leading-spectra}. The formulas are continued through $\omega=0$ for plotting. The unshaded region is not a region of spectral stability.}
 \label{fig:low-high-stability}
\end{figure}

Theorems~\ref{thm:BF} and \ref{thm:high} establish instability at sufficiently small nonzero amplitude whenever
\begin{equation*} \mathrm{ind}_{\mathrm{low}}(\omega,\kappa)>0 \qquad\text{or}\qquad \beta_{12}(\omega,\kappa)\ne0.
\end{equation*}
Figure~\ref{fig:low-high-stability} compares these conditions. Since the high-frequency index is nonnegative, we plot the zero set of its coupling coefficient.  The two criteria leave stability undecided when
\begin{equation}\label{eqn:undecided-leading-indices}
 \beta_{12}(\omega,\kappa)=0, \qquad \mathrm{ind}_{\mathrm{low}}(\omega,\kappa)\leq0.
\end{equation}

The curves are evaluated from the formulas in Appendices~\ref{index:low} and \ref{A:high-index}.  Panel~(b) uses the shear coordinate of \cite{TKM12}, which in our normalization is
\begin{equation*} X=-\frac{\omega}{\mu_0+\omega} =-\frac{\rho\tanh\kappa}{\kappa}>-1.
\end{equation*}
The low-frequency boundary and the coupling zero set have a numerically located intersection at
\begin{equation*} (\omega,\kappa)\approx(-7.094942,\,3.698692).
\end{equation*}

\clearpage

\appendix
\section{Details of the calculations}
\label{A:details}

We collect functional-analytic details and algebraic calculations. The subsections are arranged in the order in which the results are first used in Sections~\ref{sec:Stokes}, \ref{sec:spec}, and \ref{sec:BF}.

\subsection{The order-three Stokes coefficient}\label{app:Stokes-order3}

With $s=\sinh\kappa$ and $c=\cosh\kappa$, the polynomials in \eqref{eqn:coefmu2} are
\begin{align*}
 N_1={}&\frac98\kappa^4c-\frac18\kappa^4c^5 -\frac58\kappa^4c^3-\frac14 \mu_0^6\omega^2c^3 +\frac18 \mu_0^6\omega^2c^5+\frac18 \mu_0^6\omega^2c\\
 &+\frac12\kappa^2 \mu_0^3\omega c +\frac12\kappa^2 \mu_0^2\omega^2c -\frac12\kappa^2 \mu_0^3\omega c^5 -\kappa^2 \mu_0^2\omega^2c^3 +\frac12\kappa^2\mu_0^2\omega^2c^5,\\
 N_2={}&\frac14\kappa^3\mu_0\omega+\frac18\kappa \mu_0^5\omega +\frac38\kappa^3\mu_0^2c^2-\frac98\kappa^3\mu_0^2c^4 -\frac18\kappa \mu_0^4\omega^2\\
 &-\frac12\kappa \mu_0^5\omega c^2 +\frac38\kappa \mu_0^5\omega c^4 +\frac12\kappa \mu_0^4\omega^2c^2 -\frac38\kappa \mu_0^4\omega^2c^4 .
\end{align*}

\subsection{The good unknown}
\label{A:good-unknown-proof}

We give the details needed to justify the transformations in Section~\ref{sec:linearize}.  Differentiation of the flattened equations at \((\phi_\epsilon,\psi_\epsilon,\eta_\epsilon)\) produces terms in which a derivative falls on the moving-boundary correction.  On substituting
\begin{equation*}
 \phi^{\rm old}=\phi+ \frac{y\phi_{\epsilon,y}}{1+\eta_\epsilon}\eta, \qquad \psi^{\rm old}=\psi+ \frac{y\psi_{\epsilon,y}}{1+\eta_\epsilon}\eta,
\end{equation*}
these highest-order terms cancel.  The remaining interior equations are the first two equations of \eqref{eqn:linearized}.  At \(y=1\), use the traveling-wave equations for the basic state before collecting the coefficients of \(\eta_x,\eta,\phi_y,\psi_y\).  This gives exactly \(A_1,A_2,A_3,A_4\) in \eqref{coefficients:A}.  This calculation is an identity between differential expressions and therefore holds in the Sobolev trace sense by density.

Differentiating \eqref{eqn:barpsi} with respect to \(y\) yields
\begin{equation}\label{eqn:good-unknown}
 A_1\psi_y+\omega\mu_\epsilon\psi =\frac{A_3}{\mu_0^2}\chi_y-\lambda\phi-A_4\phi_y.
\end{equation}
At \(y=1\), the fourth equation of \eqref{eqn:linearized} and \eqref{eqn:good-unknown} give \(\chi_y(1)+\mu_0^2\eta=0\); at \(y=0\), the bottom conditions give \(\chi(0)=0\).  Conversely, solving the scalar initial-value problem \eqref{eqn:good-unknown}, \(\psi(0)=0\), gives \eqref{Psi:recover}.  Hence the transformation is one-to-one.

Let \(1\le r\le2\).  Multiplication by a \(C^2\)-coefficient is bounded on \(H^r(0,1)\), and the map
\begin{equation*}
 f\longmapsto e^{-cy}\int_0^y e^{cs}f(s)\,ds
\end{equation*}
is bounded from \(H^{r-1}(0,1)\) to \(H^r(0,1)\), locally uniformly in \(c\in\mathbb C\).  Together with \eqref{eqn:A13}, this proves that both \eqref{eqn:barpsi} and \eqref{Psi:recover} are bounded on the spaces used in Section~4, with bounds locally uniform in \((\lambda,\epsilon)\).  It also proves analytic dependence on these parameters.  No step requires \(\omega\ne0\).

Differentiating \eqref{eqn:barpsi} in \(x\), inserting the first two equations of \eqref{eqn:linearized}, and using \eqref{eqn:good-unknown} gives
\begin{align}
\chi_x={}&\frac{\eta_{\epsilon,x}}{1+\eta_\epsilon}(y\chi_y-\chi) -\frac{A_{3,x}}{A_3}\chi \notag\\
&+\frac{\mu_0^2}{A_3}\Bigg[ \left(A_{1,x}+\frac{\eta_{\epsilon,x}}{1+\eta_\epsilon}A_1\right)\psi +\left(A_{4,x}+\frac{\eta_{\epsilon,x}}{1+\eta_\epsilon}A_4\right) (\phi-\phi(0)) \notag\\
&\qquad+\frac1{1+\eta_\epsilon} \left\{\left(\lambda-\frac{\omega\mu_\epsilon A_4}{A_1}\right)\psi -\left(\omega\mu_\epsilon+\frac{\lambda A_4}{A_1}\right) (\phi-\phi(0)) \right.\notag\\
&\hspace{4.2cm}\left. +\frac{A_3A_4}{\mu_0^2A_1}(\chi_y-\chi_y(0)) -\frac{A_1^2+A_4^2}{A_1}\phi_y\right\}\Bigg].
\label{eqn:chi-evolution}
\end{align}
Equations \eqref{Psi:recover} and \eqref{eqn:chi-evolution}, together with the surface equation
\begin{equation}\label{eqn:eta-x-explicit}
 \eta_x=-\frac1{A_1}\left( \frac{\eta_{\epsilon,x}}{(1+\eta_\epsilon)^2}\psi_y(1) -\frac1{(1+\eta_\epsilon)^2}\phi_y(1)+A_2\eta\right),
\end{equation}
give \eqref{eqn:spatial-system}.  More explicitly, substitute \eqref{Psi:recover} into the first equation of \eqref{eqn:linearized}, into \eqref{eqn:chi-evolution}, and into \eqref{eqn:eta-x-explicit}; subtract the three components of \(\mathbf L(\lambda)\) in \eqref{eqn:operatorL}.  The remaining three-component expression is, by definition, \(\mathbf B(x;\lambda,\epsilon)u\).  Every coefficient in this remainder vanishes at \(\epsilon=0\), and the trace maps \(H^2(0,1)\to\mathbb C\) are bounded.  This proves \eqref{eqn:B}, including its analytic parameter dependence.

\subsection{Closedness and compact resolvent}
\label{A:operator-properties}

\begin{lemma}\label{lem:L-closed-compact}
For every \(\lambda\in\mathbb C\), \(\mathbf L(\lambda):\mathcal D\subset Y\to Y\) is densely defined and closed.  Its resolvent is compact whenever it exists, and its resolvent set is nonempty.  On compact \(\lambda\)-sets there is a constant \(C\) such that
\begin{equation}\label{eqn:graph-norm-equivalence}
 \|u\|_{\mathcal D} \le C\big(\|u\|_Y+\|\mathbf L(\lambda)u\|_Y\big), \qquad u\in\mathcal D.
\end{equation}
\end{lemma}

\begin{proof}
To see density, start with arbitrary \((\phi,\chi,\eta)\in Y\).  Smooth approximants can be chosen with \(\chi(0)=0\).  The derivative conditions at \(y=0,1\) may then be imposed by adding smooth boundary-layer functions whose prescribed endpoint derivative is fixed while their \(H^1\)-norm tends to zero. Thus one obtains \(\phi_n\to\phi\) and \(\chi_n\to\chi\) in \(H^1\), with \(\phi_{n,y}(0)=0\) and \(\chi_{n,y}(1)=-\mu_0^2\eta_n\) for some \(\eta_n\to\eta\). Consequently \((\phi_n,\chi_n,\eta_n)\in\mathcal D\) and \(\mathcal D\) is dense in \(Y\). The Volterra operator \(\mathcal V_\lambda\) in \eqref{def:V} is bounded from \(H^s\times H^s_D\) to \(H^{s+1}_D\), \(s=0,1\).  Consequently the lower-order terms in \(\mathbf L(\lambda)\) are bounded from \(Y\) to itself, locally uniformly in \(\lambda\).

Writing out the first two components of \(\mathbf L(\lambda)u\), and then differentiating once in \(y\), gives directly
\begin{equation*} \|\phi\|_{H^2}+\|\chi\|_{H^2} \le C_K\big(\|\mathbf L(\lambda)u\|_Y+\|u\|_Y\big), \qquad \lambda\in K,
\end{equation*}
for every compact \(K\subset\mathbb C\). Here the differentiated Volterra term is estimated by its defining first-order equation \(\upsilon_y=\omega\mu_0\upsilon-\omega\mu_0\chi+\lambda\phi\). The boundary relation in \(\mathcal D\) controls no additional infinite-dimensional component, and \(|\eta|\) is already part of the \(Y\)-norm.  This proves \eqref{eqn:graph-norm-equivalence}.

Let \((u_n)\subset\mathcal D\) satisfy \(u_n\to u\) in \(Y\) and \(\mathbf L(\lambda)u_n\to f\) in \(Y\).  The first two components of \eqref{eqn:operatorL} express \(\chi_{n,y}\) and \(\phi_{n,y}\) as the sum of a convergent \(H^1\)-sequence and lower-order Volterra terms.  Hence \((\phi_n,\chi_n)\) is Cauchy in \(H^2\times H^2\).  The trace theorem passes \(\phi_y(0)=0\), \(\chi(0)=0\), and \(\chi_y(1)+\mu_0^2\eta=0\) to the limit.  Thus \(u\in\mathcal D\) and \(\mathbf L(\lambda)u=f\), proving closedness.

We next establish the Fredholm property without assuming that the resolvent is nonempty.  Fix \(z\in\mathbb C\) and write \((\mathbf L(\lambda)-z)u=f\).  Put \(a_0=\omega\mu_0\).  The first two equations, together with the defining equation for \(\upsilon\), form the first-order system
\begin{equation}\label{eqn:resolvent-first-order-system}
\begin{aligned}
 \phi_y&=f_2-\lambda\upsilon+(\lambda+z)\chi +a_0(\phi-\phi(0)),\\
 \chi_y&=a_0\upsilon-a_0\chi+(\lambda-z)\phi-f_1,\\
 \upsilon_y&=a_0\upsilon-a_0\chi+\lambda\phi,
\end{aligned}
\qquad \chi(0)=\upsilon(0)=0.
\end{equation}
For a prescribed scalar \(c=\phi(0)\), variation of constants gives a unique solution in \(H^2\times H^2\times H^2\), depending boundedly and affinely on \((c,f_1,f_2)\).  The domain condition gives \(\eta=-\chi_y(1)/\mu_0^2\), and the third resolvent equation is the single scalar boundary condition
\begin{equation}\label{eqn:scalar-resolvent-boundary}
 -\frac{\lambda-z}{\mu_0^2}\chi_y(1)-\phi_y(1)=f_3.
\end{equation}
For the homogeneous problem the left-hand side of \eqref{eqn:scalar-resolvent-boundary} is \(M(z,\lambda)c\), where \(M\) is analytic and not identically zero. Thus \(\mathbf L(\lambda)-z\) is Fredholm of index zero: when \(M\ne0\) it is bijective, while at a zero of \(M\) its kernel and cokernel are both one dimensional; the order of that zero records the algebraic multiplicity.  Equivalently, this is the one-dimensional Schur complement of the ODE solution operator in \eqref{eqn:resolvent-first-order-system}.

The large-imaginary-parameter estimate proved directly in Appendix~\ref{A:resolvent-estimate} implies that \(\mathbf L(\lambda)-ia\) is injective for all sufficiently large \(|a|\).  The index-zero property then makes it surjective.  Hence the resolvent set is nonempty.  Estimate \eqref{eqn:graph-norm-equivalence} shows that the resolvent maps \(Y\) boundedly into \(\mathcal D\).  Since \(\mathcal D\hookrightarrow Y\) is compact by Rellich's theorem and the trace component is finite dimensional, the resolvent is compact.
\end{proof}

The same graph-norm argument, now in the \(x\)-variable, shows that \(\mathcal L_\theta(\lambda,\epsilon)\) in \eqref{eqn:Bloch} is closed on the parameter-independent domain \(\mathcal D_{\rm per}\).  Since \(\mathbf B\) is bounded from \(\mathcal D\) to \(Y\), it is a graph-bounded perturbation and does not change this conclusion.  In particular, with the graph norm on \(\mathcal D_{\rm per}\), the map \(\theta\mapsto\mathcal L_\theta(\lambda,\epsilon)\) is continuous, and so is its inverse wherever the latter exists.

\subsection{Bloch decomposition}
\label{A:bloch-proof}

Let \(\mathscr B\) denote the unitary Bloch transform
\begin{equation*}
 L^2(\mathbb R;Y)\longrightarrow \int_{[-\kappa/2,\kappa/2)}^\oplus X_{\rm per}\,d\theta.
\end{equation*}
On smooth compactly supported functions it may be normalized by
\begin{equation*}
 (\mathscr Bu)(\theta,x) =\left(\frac{T}{2\pi}\right)^{1/2} \sum_{n\in\mathbb Z}e^{-i\theta(x+nT)}u(x+nT), \qquad 0<x<T,
\end{equation*}
and extends by density to a unitary map.  The same formula and Plancherel's theorem identify the whole-line domain with the direct integral of the common fiber domain \(\mathcal D_{\rm per}\). On the common domain it satisfies
\begin{equation}\label{eqn:Bloch-direct-integral}
 \mathscr B\mathcal L(\lambda,\epsilon)\mathscr B^{-1} =\int_{[-\kappa/2,\kappa/2)}^\oplus \mathcal L_\theta(\lambda,\epsilon)\,d\theta.
\end{equation}
Suppose first that every fiber is invertible.  By the common-domain graph estimate following Lemma~\ref{lem:L-closed-compact}, the map \(\theta\mapsto\mathcal L_\theta\) is continuous as an element of \(\mathcal L(\mathcal D_{\rm per},X_{\rm per})\), where the domain has its graph norm.  The inverse map is therefore locally norm-continuous. Compactness of \([-\kappa/2,\kappa/2]\) yields
\begin{equation*}
 \sup_\theta\|\mathcal L_\theta^{-1}\|<\infty,
\end{equation*}
and \eqref{eqn:Bloch-direct-integral} gives a bounded inverse of the whole-line pencil.  Conversely, if a fiber at \(\theta_0\) is not invertible, the direct-integral criterion gives noninvertibility of the whole-line pencil.  Explicitly, approximate null vectors for the fiber (or for its adjoint if surjectivity fails), multiplied by \(L^2\)-normalized cutoffs supported in shrinking \(\theta\)-neighborhoods of \(\theta_0\), form an approximate null sequence for the direct integral.  This proves the equivalence in Definition~\ref{def:temporal-spectrum} without any unproved assertion about the point spectrum.

\subsection{Proof of Lemma~\ref{lem:sym}}
\label{A:symm}

All coefficients of the linearized problem are real.  Consequently
\begin{equation*}
 \mathcal L_\theta(\lambda,\epsilon)u=0 \quad\Longrightarrow\quad \mathcal L_{-\theta}(\overline\lambda,\epsilon)\overline u=0.
\end{equation*}
With the crest at \(x=0\), the basic Stokes wave has the usual reversible parity.  Direct substitution into \eqref{eqn:linearized} gives
\begin{equation*}
 \mathcal R(\phi,\chi,\eta)(x,y) =(-\phi(-x,y),\chi(-x,y),\eta(-x))
\end{equation*}
and
\begin{equation*}
 \mathcal L_\theta(\lambda,\epsilon)u=0 \quad\Longrightarrow\quad \mathcal L_{-\theta}(-\lambda,\epsilon)\mathcal Ru=0.
\end{equation*}
The same identities hold for approximate null sequences and hence for noninvertibility, not merely for eigenfunctions.  They prove \eqref{eqn:symm}.

\subsection{The boundary determinant at zero amplitude}
\label{A:eps=0}

Let \(k\in\mathbb R\setminus\{0\}\) and suppose \(\mathbf L(\lambda)U=ikU\).  At zero amplitude the first two untransformed interior equations are the Cauchy--Riemann system
\begin{equation}\label{eqn:flat-CR}
 ik\phi=\psi_y, \qquad ik\psi=-\phi_y.
\end{equation}
Thus \(\phi_{yy}=k^2\phi\), and the bottom condition gives \(\phi=C\cosh(ky)\).  If \(C=0\), then \eqref{eqn:flat-CR} and the boundary conditions imply \(\psi=\eta=0\); hence every eigenfunction has \(C\ne0\), and we normalize \(C=1\).  Substitution into the surface conditions gives the scalar determinant
\begin{equation*}
 d(k,\lambda)= \lambda^2+i(\omega\mu_0\tanh k-2k)\lambda -k^2+(\mu_0^2+\omega\mu_0)k\tanh k.
\end{equation*}
Solving \(d(k,\lambda)=0\) gives
\begin{equation*} \lambda=-i\left(\frac{\omega\mu_0\tanh k}{2}-k \pm\sqrt{\mu_0^2k\tanh k+ \frac{\omega^2\mu_0^2\tanh^2k}{4}}\right).
\end{equation*}
This derivation does not divide by \(\omega\), and hence includes the irrotational case.

In the \((\phi,\chi,\eta)\)-variables, direct substitution gives
\begin{equation*} U(k,\sigma)=
\begin{pmatrix}
\cosh(ky)\\[1mm]
\dfrac1{ik}\big((k-\sigma)\sinh(ky) -\omega\mu_0\cosh(ky)+\omega\mu_0\big)\\[1mm]
\dfrac{i}{\mu_0^2}\big((k-\sigma)\cosh k -\omega\mu_0\sinh k\big)
\end{pmatrix}.
\end{equation*}
Indeed, the first two components of \((\mathbf L(i\sigma)-ik)U\) vanish identically in \(y\), the bottom conditions hold, and its third-component residual is exactly
\begin{equation*}
 \big[(\mathbf L(i\sigma)-ik)U\big]_3 =-\frac{\cosh k}{\mu_0^2}\,d(k,i\sigma).
\end{equation*}

At \(k=0\), solve the equations before division by \(k\).  For \(\sigma=0\), the kernel is \(\operatorname{span}\{(1,0,0)^T\}\), and the vector
\begin{equation*}
 \left(0,\frac{\omega\mu_0}{2}y^2-y, \frac{1-\omega\mu_0}{\mu_0^2}\right)^T
\end{equation*}
is mapped to \((1,0,0)^T\).  This proves the identities in \eqref{eqn:Jordan0}, again for all \(\omega\in\mathbb R\).

\subsection{The dispersion curves}
\label{A:disp}

Put
\begin{equation*}
 \rho=\omega\mu_0, \qquad t=\tanh k, \qquad u=\sech^2k, \qquad Q(k)=\mu_0^2kt+\frac{a^2t^2}{4}.
\end{equation*}
For \(k>0\), \(Q>0\) and
\begin{equation*}
 Q'=\mu_0^2(t+ku)+\frac{\rho^2}{2}tu>0.
\end{equation*}
If \(\rho\le0\), the formula
\begin{equation*}
 \sigma_+'=1-\frac{\rho u}{2}+\frac{Q'}{2\sqrt Q}
\end{equation*}
is immediately positive.  If \(\rho>0\), the exact identity
\begin{equation}\label{eqn:sigmaplus-positive-identity}
 (Q')^2-\rho^2u^2Q =\mu_0^4(t+ku)^2+\mu_0^2\rho^2t^2u>0
\end{equation}
implies \(Q'/\sqrt Q>\rho u\), and therefore \(\sigma_+'(k)>1\).  This proves part (i) of Lemma~\ref{lem:sigma-geometry}.

To prove strict convexity, write
\begin{equation*}
 R(k)=kt+\frac{\omega^2t^2}{4}, \qquad Q=\mu_0^2R.
\end{equation*}
A direct differentiation, with \(t'=u\) and \(u'=-2tu\), gives
\begin{equation}\label{eqn:sigmaminus-second}
 \sigma_-''(k)=\mu_0\left( \omega tu+\frac{S}{4R^{3/2}}\right),
\end{equation}
where
\begin{equation}\label{eqn:S-correct}
 S=(t-ku)^2+4k^2t^2u+3\omega^2kt^3u +\frac{\omega^4}{2}t^4u.
\end{equation}
Expanding the first term in \eqref{eqn:S-correct} gives \((t-ku)^2=t^2-2ktu+k^2u^2\).

If \(\omega\ge0\), positivity in \eqref{eqn:sigmaminus-second} is immediate.  Suppose \(\omega<0\), and set \(x=k\), \(z=\omega^2t/4\).  Dropping the first positive term in \eqref{eqn:S-correct} gives
\begin{align*}
 S&>4t^2u(x+z)(x+2z),\\
 4|\omega|R^{3/2}tu
 &=8t^2u\sqrt z\,(x+z)^{3/2}.
\end{align*}
Here \(t-ku>0\) for \(k>0\), since this function vanishes at zero and has derivative \(2ktu>0\). But
\begin{equation*}
 x+2z>2\sqrt{z(x+z)},
\end{equation*}
because the difference of the squares is \(x^2>0\).  Hence \(S>4|\omega|R^{3/2}tu\), and \eqref{eqn:sigmaminus-second} is strictly positive also when \(\omega<0\).

Finally, \eqref{eqn:mu0} implies
\begin{equation*}
 \sigma_-(\kappa)=0, \qquad \sigma_-(0)=0.
\end{equation*}
Strict convexity gives \(\sigma_-<0\) on \((0,\kappa)\), a unique critical point \(r_c\in(0,\kappa)\), and strict increase on \((\kappa,\infty)\).  The symmetry \(\sigma_+(-k)=-\sigma_-(k)\) completes the proof of Lemma~\ref{lem:sigma-geometry} and all stated root counts.

\subsection{Algebraic multiplicity. Jordan chains}
\label{A:Jordan}

Fix \(\sigma\in\mathbb R\).  The boundary determinant for the spatial eigenvalue \(z=ik\) is, up to a nonvanishing analytic factor,
\begin{equation}\label{eqn:Dksigma}
 D(k,\sigma)=d(k,i\sigma) =-\sigma^2+(2k-\omega\mu_0\tanh k)\sigma -k^2+(\mu_0^2+\omega\mu_0)k\tanh k.
\end{equation}
The kernel of the interior boundary-value problem is one dimensional for \(k\ne0\).  Analytic Fredholm reduction therefore identifies the algebraic multiplicity of \(ik\) as an eigenvalue of \(\mathbf L(i\sigma)\) with the order of \(k\) as a zero of \(D(\cdot,\sigma)\).  Thus the simple-root assertions in Lemma~\ref{lem:flat-eigenvalues} follow from Lemma~\ref{lem:sigma-geometry}.

At \(\sigma=0\),
\begin{equation*} D(k,0)=k\big((\mu_0^2+\omega\mu_0)\tanh k-k\big).
\end{equation*}
Since \(\mu_0^2+\omega\mu_0=\kappa\coth\kappa>1\),
\begin{equation*}
 D(k,0)= (\kappa\coth\kappa-1)k^2+O(k^4) \quad\text{as }k\to0.
\end{equation*}
Hence \(k=0\) has order exactly two.  At \(k=\pm\kappa\),
\begin{equation*}
 \partial_k\big((\mu_0^2+\omega\mu_0)\tanh k-k\big) =\frac{2\kappa}{\sinh(2\kappa)}-1\ne0,
\end{equation*}
so those roots are simple.  The direct identities in \eqref{eqn:Jordan0} exhibit the length-two chain at zero; the order-two determinant excludes a longer chain.

At \((k_c,\sigma_c)\), factorization of \eqref{eqn:Dksigma} gives \(D=-(\sigma-\sigma_+(k))(\sigma-\sigma_-(k))\).  Hence \(D_\sigma(k_c,\sigma_c)=-2\sqrt{Q(k_c)}\ne0\), and strict convexity, together with \(\sigma_+(-k)=-\sigma_-(k)\), gives
\begin{equation*}
 D(k_c,\sigma_c)=0, \qquad D_k(k_c,\sigma_c)=0, \qquad D_{kk}(k_c,\sigma_c)\ne0.
\end{equation*}
Differentiate
\begin{equation*}
 (\mathbf L(i\sigma_c)-ik)\Upphi(k,\sigma_c)=0
\end{equation*}
with respect to \(k\).  The differentiated surface condition is satisfied precisely because \(D_k(k_c,\sigma_c)=0\), and one obtains
\begin{equation*}
 (\mathbf L(i\sigma_c)-ik_c) \big(-i\partial_k\Upphi(k_c,\sigma_c)\big) =\Upphi(k_c,\sigma_c).
\end{equation*}
This proves \eqref{eqn:Jordanc}; the nonzero second derivative excludes a chain of length three.

\subsection{The adjoint problem and biorthonormal coordinates}
\label{A:adjoint-coordinate-formulas}

We give a self-contained calculation of the dual vectors used in \eqref{eqn:biorthogonal-projection-zero} and \eqref{eqn:biorthogonal-projection-high}.  In particular, all the quantities entering the two projections are displayed explicitly below in terms of elementary functions.  Write \(\rho=\omega\mu_0\), let \(\lambda\in\mathbb C\), and, for \(V=(\Phi,\mathsf X,H)^T\), set
\begin{equation*} \nu(y)=\int_0^y e^{\rho s} \big(\rho\Phi(s)+\overline\lambda\mathsf X(s)\big)\,ds.
\end{equation*}
Define the common expression
\begin{align*}
 \mathcal R_V(y)={}&e^{-\rho y}\left\{ \nu(1)-\nu(y)+\int_y^1\rho e^{\rho s} \big(\rho\Phi_y(s)+\overline\lambda\mathsf X_y(s)\big)\,ds \right\} \\
 &+\rho\Phi_y(y)+\overline\lambda\mathsf X_y(y).
\end{align*}
Let \(p_\Phi,p_{\mathsf X}\in H^2(0,1)\) be the unique solutions of
\begin{align}
 p_\Phi-p_{\Phi,yy}&=\overline\lambda\,\mathcal R_V,
 &p_{\Phi,y}(0)&=-\rho\int_0^1\mathsf X(s)\,ds,
 &p_{\Phi,y}(1)&=\mathsf X(1),
 \label{eqn:pPhi-BVP}\\
 p_{\mathsf X}-p_{{\mathsf X},yy}&=-\rho\,\mathcal R_V,
 &p_{\mathsf X}(0)&=\rho\Phi(0)+\overline\lambda\mathsf X(0)-\Phi_y(0),
 &p_{{\mathsf X},y}(1)&=-\Phi(1).
 \label{eqn:pChi-BVP}
\end{align}
The one-dimensional operators in \eqref{eqn:pPhi-BVP} and \eqref{eqn:pChi-BVP} are isomorphisms with the displayed separated boundary conditions.  Hence these equations define bounded linear functions of \(V\), with no hidden solvability condition.

For later reference, these two solutions can be written without any implicit differential equation.  If \(f_\Phi=\overline\lambda\mathcal R_V\), \(f_{\mathsf X}=-\rho\mathcal R_V\), and
\begin{equation*}
 a_\Phi=-\rho\int_0^1\mathsf X(s)\,ds, \qquad a_{\mathsf X}=\rho\Phi(0)+\overline\lambda\mathsf X(0)-\Phi_y(0),
\end{equation*}
then variation of parameters gives
\begin{align*}
 p_\Phi(y)
 &=C_\Phi\cosh y+a_\Phi\sinh y -\int_0^y\sinh(y-s)f_\Phi(s)\,ds,\\
 p_{\mathsf X}(y)
 &=a_{\mathsf X}\cosh y+C_{\mathsf X}\sinh y -\int_0^y\sinh(y-s)f_{\mathsf X}(s)\,ds,
\end{align*}
where
\begin{align*}
 C_\Phi
 &=\frac{\mathsf X(1)-a_\Phi\cosh1 +\int_0^1\cosh(1-s)f_\Phi(s)\,ds}{\sinh1}, \\
 C_{\mathsf X}
 &=\frac{-\Phi(1)-a_{\mathsf X}\sinh1 +\int_0^1\cosh(1-s)f_{\mathsf X}(s)\,ds}{\cosh1}.
\end{align*}
Thus every integral occurring below is an elementary integral of exponentials.

\begin{lemma}[Adjoint boundary-value problem]
\label{lem:adjoint-BVP}
With respect to \eqref{eqn:Y-product},
\begin{align}
 \mathcal D(\mathbf L(\lambda)^*) =\big\{(\Phi,\mathsf X,H)\in H^2(0,1)\times H_D^2(0,1) \times\mathbb C:\ \Phi_y(0)=0,\ \mathsf X_y(1)=H\big\},
\label{eqn:adjoint-domain}
\end{align}
and
\begin{equation}\label{eqn:adjoint-operator-visible}
 \mathbf L(\lambda)^*V=
 \begin{pmatrix}
 \overline\lambda\Phi-\rho\mathsf X-\mathsf X_y+p_\Phi\\
 -\rho\Phi-\overline\lambda\mathsf X+\Phi_y+p_{\mathsf X}\\
 \overline\lambda H+\mu_0^2\Phi_y(1)
 \end{pmatrix}.
\end{equation}
\end{lemma}

\begin{proof}
Insert \eqref{eqn:operatorL} and \eqref{def:V} in \(\langle\mathbf L(\lambda)U,V\rangle_Y\), change the order in each Volterra integral, and integrate the two derivative terms by parts. The coefficients of the free traces \(\phi(0)\), \(\phi_y(1)\), and \(\chi(1)\) vanish precisely when \eqref{eqn:adjoint-domain} and the boundary conditions in \eqref{eqn:pPhi-BVP}--\eqref{eqn:pChi-BVP} hold.  The remaining volume terms are \(\langle U,\mathbf L(\lambda)^*V\rangle_Y\) with \eqref{eqn:adjoint-operator-visible}.  Conversely, if one of these trace conditions fails, a sequence of smooth functions with bounded \(Y\)-norm and the corresponding unbounded trace shows that the functional is not continuous on \(Y\).  This proves both the formula and the asserted adjoint domain.
\end{proof}

We now record the dual vectors at the two spectral anchors in a form which will be useful below.  Since the second entry of \eqref{eqn:Y-product} is conjugated, it is convenient to display \(\overline{\uppsi_j}\), the quantity which occurs directly in the biorthogonal pairings.  Throughout this calculation, put
\begin{equation*} \rho=\omega\mu_0.
\end{equation*}

We begin at \(\sigma=0\).  For \(j=1,2\), let \(k=k_j(0)\), where \(k_1(0)=-\kappa\) and \(k_2(0)=\kappa\).  Notice that
\begin{equation}\label{eqn:low-root-simplification}
 k\cosh k-\rho\sinh k=\mu_0^2\sinh k.
\end{equation}
A straightforward calculation, together with \eqref{eqn:low-root-simplification}, gives
\begin{equation}\label{eqn:adjoint-simple-zero}
\overline{\uppsi_j(0)}(y)=p_j^0
\begin{pmatrix}
\displaystyle -\frac{2i(k-\rho)}{k} \left(\cosh(ky)-\frac{\sinh k}{k\sinh1}\cosh y\right)\\[3mm]
\displaystyle e^{-ky}-\frac{k-\rho}{k+\rho}e^{ky} +p_{1,j}^0e^{-\rho y}+p_{2,j}^0e^{-y}+p_{3,j}^0e^y\\[3mm]
\displaystyle -ke^{-k}-\frac{k(k-\rho)}{k+\rho}e^k -\rho p_{1,j}^0e^{-\rho} -p_{2,j}^0e^{-1}+p_{3,j}^0e
\end{pmatrix}.
\end{equation}
The first two coefficients in \eqref{eqn:adjoint-simple-zero} are
\begin{align*}
 p_j^0
 &=\frac{ik^2} {2(k-\rho)(k^2-1)(\sinh k\cosh k-k)}, \\
 p_{1,j}^0
 &=-\frac{2\rho^3\mu_0^2e^\rho(k^2-1)\sinh k} {k^3(k+\rho)(\rho^2-1)}.
\end{align*}
We keep the remaining coefficients in the original exponential basis by writing
\begin{equation*} p_{2,j}^0=\frac{C_j^0-B_j^0}{2}, \qquad p_{3,j}^0=\frac{C_j^0+B_j^0}{2}.
\end{equation*}
Here, upon setting
\begin{equation*} H_0(k)=k^3(k+\rho)(\rho^2-1),
\end{equation*}
a straightforward simplification gives
\begin{align*}
 B_j^0
 &=\frac{2}{H_0(k)}\Big[ k^3\rho(\rho^2-1)\tanh1 -\mu_0^2\sinh k\bigl\{ (k^2-\rho^2)\sech1 +e^\rho\rho^3(k^2-1)\tanh1\bigr\}\Big], \\
 C_j^0
 &=\frac{2}{H_0(k)}\Big[ e^\rho\mu_0^2\rho^3(k^2-1)\sinh k -k^3\rho(\rho^2-1)\Big].
\end{align*}
Equivalently, \(p_{2,j}^0e^{-y}+p_{3,j}^0e^y =C_j^0\cosh y+B_j^0\sinh y\); this identity is used only to shorten the displayed coefficient formulas.  The scalar \(p_j^0\) is chosen so that the two simple modes are biorthonormal.  Since their eigenvalues are distinct, the cross pairings vanish by adjointness, and hence
\begin{equation*} \langle\Upphi_\ell(0),\uppsi_j(0)\rangle_Y=\delta_{\ell j}, \qquad \ell,j\in\{1,2\}.
\end{equation*}

The zero Jordan block admits a shorter formula.  Put
\begin{equation*} \rho=\omega\mu_0,\qquad q=\mu_0^2+\rho=\kappa\coth\kappa,\qquad A=q-1,\qquad C_0=\frac{\mu_0^2}{A(1-\rho^2)}.
\end{equation*}
Since \(A>0\), direct substitution into the adjoint boundary problem and the chain relation gives
\begin{align}
 \overline{\uppsi_3(0)}(y)
 &=\begin{pmatrix}
 0\\[1mm]
 C_0\left\{ \rho e^\rho\left(\dfrac{\cosh(1-y)}{\cosh1}-e^{-\rho y}\right) -\dfrac{\sinh y}{\cosh1}\right\}\\[2mm]
 -\dfrac{\mu_0^2}{A}
 \end{pmatrix},
 \label{eqn:dual-zero-column-three}\\
 \overline{\uppsi_4(0)}(y)
 &=\begin{pmatrix}
 \dfrac{q-\cosh y/\sinh1}{A}\\[2mm]0\\0
 \end{pmatrix}.
 \label{eqn:dual-zero-column-four}
\end{align}
In the coefficient notation of the unreduced formulas, the index \(j\) is spurious and the four coefficients become
\begin{align*}
 p_4^0&=\frac{\mu_0^2\rho e^\rho}{A(1-\rho^2)},
 & p_5^0&=-\frac{\mu_0^2 (\sech1+\rho e^\rho\tanh1)}{A(1-\rho^2)}, \\
 p_6^0&=\frac qA,
 & p_7^0&=-\frac1{A\sinh1}.
\end{align*}
Indeed, the common denominator in the first two coefficients factors as
\begin{equation*} -\omega^2\mu_0^4-\omega^3\mu_0^3+\omega^2\mu_0^2 +\mu_0^2+\omega\mu_0-1=A(1-\rho^2).
\end{equation*}
The apparently singular cases \(\rho=\pm1\) in \eqref{eqn:dual-zero-column-three} are understood by analytic continuation.  Equivalently, one replaces the colliding exponential by the corresponding polynomial multiple.  The vector itself remains regular.

At the high-frequency anchor, write
\begin{equation*} r=k_4(\sigma_H),\qquad k_2(\sigma_H)=r+2\kappa.
\end{equation*}
For each \(j\in\{2,4\}\), we abbreviate \(k=k_j(\sigma_H)\).  Thus \(k=r+2\kappa\) for \(j=2\), whereas \(k=r\) for \(j=4\).  A straightforward calculation reveals that
\begin{equation}\label{eqn:adjoint-simple-high}
\overline{\uppsi_j(\sigma_H)}(y)=p_j
\begin{pmatrix}
\begin{aligned}
 &p_{1,j}\frac{i(k+\sigma_H+\rho)}{k}e^{ky} -\frac{i(k+\sigma_H-\rho)}{k}e^{-ky} -p_{2,j}\frac{i\sigma_H}{\rho}e^{-\rho y}\\
 &\hspace{28mm} +p_{3,j}e^y+p_{4,j}e^{-y}
\end{aligned}\\[4mm]
\begin{aligned}
 &p_{1,j}e^{ky}+e^{-ky}+p_{2,j}e^{-\rho y}\\
 &\hspace{28mm} +p_{5,j}e^y+p_{6,j}e^{-y}
\end{aligned}\\[4mm]
\displaystyle p_{1,j}ke^k-ke^{-k}-\rho p_{2,j}e^{-\rho} +p_{5,j}e-p_{6,j}e^{-1}
\end{pmatrix}.
\end{equation}
The two coefficients belonging to the characteristic exponents \(k\) and \(-\rho\) take the particularly simple form
\begin{align*}
 p_{1,j}
 &=-\frac{k-\rho}{k+\rho}, \\
 p_{2,j}
 &=\frac{2\rho^2e^\rho(k^2-1)} {k^2(k+\rho)(\rho^2-1)} \left\{ \frac{\rho^2-k\sigma_H}{k-\sigma_H}\sinh k -\rho\cosh k \right\}.
\end{align*}
We retain the four boundary-correction coefficients in the same exponential notation by putting
\begin{equation*} \begin{aligned}
 p_{3,j}&=\frac{B_{\Phi,j}+C_{\Phi,j}}{2},& p_{4,j}&=\frac{B_{\Phi,j}-C_{\Phi,j}}{2},\\
 p_{5,j}&=\frac{B_{\mathsf X,j}+C_{\mathsf X,j}}{2},& p_{6,j}&=\frac{B_{\mathsf X,j}-C_{\mathsf X,j}}{2}.
 \end{aligned}
\end{equation*}
It remains to give \(B_{\Phi,j}\), \(C_{\Phi,j}\), \(B_{\mathsf X,j}\), and \(C_{\mathsf X,j}\).  For brevity, let
\begin{equation*} S_a=\sinh a,\qquad C_a=\cosh a,\qquad E_\rho=e^\rho,
\end{equation*}
and define
\begin{align*}
 \mathcal N_\Phi(k,\sigma) ={}& \rho\sigma(k-\sigma) \bigl\{E_\rho\rho^2(k^2-1)C_1-(k^2-\rho^2)\bigr\}C_k \\
 &-\Big\{ E_\rho C_1\rho^2(k^2-1)\sigma(\rho^2-k\sigma) +(k^2-\rho^2) (k\rho^2+k\sigma^2-k-\rho^2\sigma) \Big\}S_k \\
 &-k^2\rho\sigma(k-\sigma)(\rho^2-1)C_1,\\
 \mathcal N_{\mathsf X}(k,\sigma) ={}& (\sigma-k) \bigl\{E_\rho\rho^3(k^2-1)S_1+(k^2-\rho^2)\bigr\}C_k \\
 &+\rho\Big\{ E_\rho S_1\rho(k^2-1)(\rho^2-k\sigma) +(k^2-\rho^2)(1-k\sigma) \Big\}S_k \\
 &-k^2\rho(\sigma-k)(\rho^2-1)S_1.
\end{align*}
The desired coefficients are
\begin{align*}
 B_{\Phi,j}
 &=-\frac{2i\,\mathcal N_\Phi(k,\sigma_H)} {k^2(k+\rho)(\rho^2-1)(k-\sigma_H)S_1}, \\
 C_{\Phi,j}
 &=-i\big\{ p_{1,j}(k+\sigma_H+\rho) +(k+\sigma_H-\rho)+\sigma_Hp_{2,j} \big\}, \\
 B_{\mathsf X,j}
 &=-p_{1,j}-1-p_{2,j}, \\
 C_{\mathsf X,j}
 &=\frac{2\,\mathcal N_{\mathsf X}(k,\sigma_H)} {k^2(k+\rho)(\rho^2-1)(k-\sigma_H)C_1}.
\end{align*}
These formulas determine, without ambiguity, all the coefficients appearing in \eqref{eqn:adjoint-simple-high}.  Thus the notation used in the calculation below is unchanged; only the common combinations in the coefficient formulas have been collected.

To state the normalization, recall from \eqref{eqn:Dksigma} that
\begin{equation*} D(k,\sigma)=-(k-\sigma)^2+ \big\{\mu_0^2k+\rho(k-\sigma)\big\}\tanh k.
\end{equation*}
Using \(D(k,\sigma_H)=0\), the scalar factor reduces to
\begin{equation*} p_j= \frac{ik^2\mu_0^2(k-\sigma_H)} {2(k-\rho)(k^2-1)\cosh^2k \{\rho\tanh k-(k-\sigma_H)\} \partial_kD(k,\sigma_H)}.
\end{equation*}
Consequently,
\begin{equation*} \langle\Upphi_\ell(\sigma_H),\uppsi_j(\sigma_H)\rangle_Y =\delta_{\ell j},\qquad \ell,j\in\{2,4\}.
\end{equation*}
The formulas above, together with \eqref{eqn:dual-zero-column-three}-- \eqref{eqn:dual-zero-column-four}, give all the dual vectors entering \eqref{eqn:biorthogonal-projection-zero} and \eqref{eqn:biorthogonal-projection-high}.

\subsection{The resolvent estimate}
\label{A:resolvent-estimate}

We prove Lemma~\ref{lem:resolvent}.  Write
\[
 L_*=\mathbf L(i\sigma_*),\qquad P\in\{P_0,P_H\},\qquad Q=I-P,\qquad L_h=L_*|_{QY}.
\]
By Lemma~\ref{lem:biorthogonal-projections}, \(P\) reduces \(L_*\). Moreover, its range contains all generalized eigenspaces of \(L_*\) on \(i\mathbb R\), by Lemma~\ref{lem:flat-eigenvalues}.  Since \(L_*\) has compact resolvent, the restriction to \(QY\) therefore satisfies
\begin{equation}\label{eqn:no-imag-spectrum-Q}
 i\mathbb R\subset\rho(L_h).
\end{equation}
It remains to establish the decay of the resolvent norm as \(|a|\to\infty\).

Put \(\rho=\omega\mu_0\) and \(\upsilon=\mathcal V_{i\sigma_*}(\phi,\chi)\).  The equation
\[
 \upsilon_y=\rho\upsilon-\rho\chi+i\sigma_*\phi,\qquad \upsilon(0)=0,
\]
implies, for \(j=0,1\),
\begin{equation}\label{eqn:Volterra-Hj-estimate}
 \|\upsilon\|_{H^{j+1}} \le C\big(\|\phi\|_{H^j}+\|\chi\|_{H^j}\big).
\end{equation}
Let \(u=(\phi,\chi,\eta)\in\mathcal D\) solve
\begin{equation}\label{eqn:resolvent-equation-A8}
 (L_*-ia)u=f=(f_1,f_2,f_3)\in Y.
\end{equation}
The first two components are equivalently
\begin{equation}\label{eqn:first-order-resolvent-pair}
 -\chi_y-ia\phi=F_1,\qquad \phi_y-ia\chi=F_2,
\end{equation}
where
\begin{align*}
 F_1&=f_1-i\sigma_*\phi+\rho\chi-\rho\upsilon,\\
 F_2&=f_2+i\sigma_*\chi+\rho(\phi-\phi(0)) -i\sigma_*\upsilon.
\end{align*}
Consequently \eqref{eqn:Volterra-Hj-estimate} and the trace inequality at \(y=0\) give
\begin{equation}\label{eqn:F-H1-bound}
 \|F_1\|_{H^1}+\|F_2\|_{H^1} \le C\big(\|f_1\|_{H^1}+\|f_2\|_{H^1} +\|\phi\|_{H^1}+\|\chi\|_{H^1}\big).
\end{equation}

Taking the \(L^2(0,1)\)-norms in \eqref{eqn:first-order-resolvent-pair}, expanding the squares, and integrating the cross terms by parts yields the exact identity
\begin{align}
 &\|\chi_y\|_2^2+\|\phi_y\|_2^2 +a^2(\|\phi\|_2^2+\|\chi\|_2^2) +2a\,\operatorname{Im}\! \big(\chi(1)\overline{\phi(1)}\big) \notag\\
 &\hspace{5cm}=\|F_1\|_2^2+\|F_2\|_2^2.
\label{eqn:zeroth-energy-identity}
\end{align}
Here the bottom contribution vanishes because \(\chi(0)=0\). Differentiating \eqref{eqn:first-order-resolvent-pair}, repeating the same computation, and using \(\phi_y(0)=0\), gives
\begin{align*}
 &\|\chi_{yy}\|_2^2+\|\phi_{yy}\|_2^2 +a^2(\|\phi_y\|_2^2+\|\chi_y\|_2^2) +2a\,\operatorname{Im}\! \big(\chi_y(1)\overline{\phi_y(1)}\big) \\
 &\hspace{5cm}=\|F_{1,y}\|_2^2+\|F_{2,y}\|_2^2.
\end{align*}
The boundary term in \eqref{eqn:zeroth-energy-identity} is lower order: the one-dimensional trace inequality gives
\begin{equation}\label{eqn:lower-boundary-term}
 2|a|\,|\chi(1)\phi(1)| \le C|a|\big(\|\phi\|_{H^1}^2+\|\chi\|_{H^1}^2\big).
\end{equation}
Combining \eqref{eqn:F-H1-bound}--\eqref{eqn:lower-boundary-term} and dropping the nonnegative \(H^2\)-terms gives
\begin{align}
 a^2S^2 \le C\|f\|_Y^2+C(1+|a|)S^2 +C|a|\,|\chi_y(1)\phi_y(1)|, \qquad S^2:=\|\phi\|_{H^1}^2+\|\chi\|_{H^1}^2.
\label{eqn:pre-trace-high-frequency}
\end{align}

It remains to control the top trace without losing a power of \(|a|\).  The domain and surface conditions in \eqref{eqn:resolvent-equation-A8} are
\begin{equation}\label{eqn:two-surface-relations-A8}
 \chi_y(1)=-\mu_0^2\eta,\qquad i(\sigma_*-a)\eta-\phi_y(1)=f_3.
\end{equation}
Hence, for \(|a|\ge2|\sigma_*|+1\),
\begin{equation}\label{eqn:boundary-product-first-bound}
 |a|\,|\chi_y(1)\phi_y(1)| \le C\big(|f_3|^2+|\phi_y(1)|^2\big).
\end{equation}
To estimate \(\phi_y(1)\), integrate the differentiated second equation in \eqref{eqn:first-order-resolvent-pair} against \(y^N\):
\begin{equation*} \phi_y(1) =N\int_0^1y^{N-1}\phi_y\,dy +ia\int_0^1y^N\chi_y\,dy +\int_0^1y^NF_{2,y}\,dy .
\end{equation*}
Cauchy--Schwarz gives
\begin{align}
 |\phi_y(1)|^2 \le{}& \frac{3N^2}{2N-1}\|\phi_y\|_2^2 +\frac{3a^2}{2N+1}\|\chi_y\|_2^2 +\frac{3}{2N+1}\|F_{2,y}\|_2^2 \notag\\
 \le{}& C_N S^2+\frac{3a^2}{2N+1}\|\chi_y\|_2^2 +\frac{C}{2N+1}\big(\|f\|_Y^2+S^2\big),
\label{eqn:weighted-trace-bound}
\end{align}
where \eqref{eqn:F-H1-bound} was used in the last line.

Insert \eqref{eqn:weighted-trace-bound} into \eqref{eqn:boundary-product-first-bound} and then into \eqref{eqn:pre-trace-high-frequency}.  First choose \(N\) so large that the coefficient of \(a^2\|\chi_y\|_2^2\) contributed by the trace is at most one quarter of the left-hand side.  With this \(N\) fixed, choose \(a_1\) so large that the terms \(C_N S^2+C(1+|a|)S^2\) are also absorbed.  We obtain
\begin{equation}\label{eqn:H1-high-frequency-final}
 |a|\,\big(\|\phi\|_{H^1}+\|\chi\|_{H^1}\big) \le C\|f\|_Y,\qquad |a|\ge a_1.
\end{equation}
Equations \eqref{eqn:weighted-trace-bound} and \eqref{eqn:H1-high-frequency-final} also give \(|\phi_y(1)|\le C\|f\|_Y\).  The second relation in \eqref{eqn:two-surface-relations-A8} therefore yields
\begin{equation*} |a|\,|\eta|\le C\|f\|_Y,\qquad |a|\ge a_1.
\end{equation*}
Together,
\begin{equation}\label{eqn:full-high-resolvent}
 \|u\|_Y\le\frac{C}{|a|}\|f\|_Y,\qquad |a|\ge a_1.
\end{equation}

The derivation of \eqref{eqn:full-high-resolvent} is an a priori estimate and does not assume invertibility.  For \(f=0\) it proves injectivity of \(L_*-ia\) for \(|a|\ge a_1\).  The index-zero Fredholm property established in Lemma~\ref{lem:L-closed-compact} then gives surjectivity, so \(ia\in\rho(L_*)\) and \eqref{eqn:full-high-resolvent} is the full resolvent estimate there. Since \(Q\) reduces \(L_*\),
\[
 (L_h-ia)^{-1}=Q(L_*-ia)^{-1}|_{QY}
\]
for these large values of \(|a|\).  On the compact interval \([-a_1,a_1]\), \eqref{eqn:no-imag-spectrum-Q} and norm-continuity of the resolvent give a finite bound.  Enlarging the constant combines the compact and large-\(|a|\) regions and proves \eqref{eqn:resolvent-estimate}.

\section{Explicit low-order monodromy coefficients}
\label{app:all-explicit-low-order-coefficients}
\begingroup \allowdisplaybreaks[3]

Throughout this appendix, $\kappa>0$ and $\omega\in\mathbb R$. The coefficients use the basis and amplitude normalization of the low-frequency calculation; the superscripts denote Taylor coefficients without factorials, evaluated at $x=2\pi/\kappa$. We display every entry of the two first-order matrices and the five second-order entries required in the modulational calculation. The explicit coefficient formulas are written in terms of $\kappa$ and $\omega$. At removable singularities the displayed quotients are understood by continuity; this convention includes $\omega=0$.

\subsection{Coefficients of order (1,0)}

\begin{align*}
a_{11}^{(1,0)}=a_{22}^{(1,0)}&={}\frac{\pi\,\sinh\kappa}{2\,\kappa^{2}\,\left(\sinh\kappa\cosh\kappa-\kappa\right)}\left(4\,\kappa\,\cosh\kappa+\left(\sqrt{\omega^2+4\kappa\coth\kappa}-\omega\right)^{2}\,\sinh\kappa\right),
\end{align*}

\begin{align*}
a_{13}^{(1,0)}=-a_{23}^{(1,0)}&={}\frac{\mathrm{i}\,\pi\,\omega\,\left(\sqrt{\omega^2+4\kappa\coth\kappa}-\omega\right)}{\kappa^{2}\,\left(\sinh\kappa\cosh\kappa-\kappa\right)} \left[4\,\kappa^{2}-\omega^{2}\,\left(\sqrt{\omega^2+4\kappa\coth\kappa}-\omega\right)^{2}\right]^{-1}\\
&\qquad{}\times \left[4\,\kappa-\left(\sqrt{\omega^2+4\kappa\coth\kappa}-\omega\right)^{2}\,\exp\!\left[\frac{\omega}{2}\left(\sqrt{\omega^2+4\kappa\coth\kappa}-\omega\right)\right]\,\sinh\kappa\right],
\end{align*}

\begin{align*}
a_{33}^{(1,0)}&={}\frac{\pi}{\kappa\,\omega\,\left(\kappa\coth\kappa-1\right)}\Big[\left(\sqrt{\omega^2+4\kappa\coth\kappa}-\omega\right)\,\left(1+\omega^{2}\right)\\&\qquad\qquad\qquad\qquad\qquad-\left(\sqrt{\omega^2+4\kappa\coth\kappa}-\omega\right)\,\exp\!\left[\frac{\omega}{2}\left(\sqrt{\omega^2+4\kappa\coth\kappa}-\omega\right)\right]-2\,\omega\Big],
\end{align*}

\begin{align*}
a_{41}^{(1,0)}=a_{42}^{(1,0)}&={}\frac{\pi\,\left(\sqrt{\omega^2+4\kappa\coth\kappa}-\omega\right)^{2}}{4\,\kappa^{2}\,\left(\kappa\coth\kappa-1\right)\left[4\,\kappa^{2}-\omega^{2}\,\left(\sqrt{\omega^2+4\kappa\coth\kappa}-\omega\right)^{2}\right]}\\
&\qquad{}\times \Big[\left(16\,\kappa^{2}+\omega\,\left(\sqrt{\omega^2+4\kappa\coth\kappa}-\omega\right)^{3}\right)\,\sinh\kappa\\*[1mm]
&\qquad\qquad-4\,\kappa\,\omega\,\left(\sqrt{\omega^2+4\kappa\coth\kappa}-\omega\right)\,\exp\!\left[\frac{\omega}{2}\left(\sqrt{\omega^2+4\kappa\coth\kappa}-\omega\right)\right]\Big],
\end{align*}

\begin{align*}
a_{43}^{(1,0)}={}\frac{\pi^{2}}{\kappa^{2}\,\left(\kappa\coth\kappa-1\right)} \left[\omega\,\left(\sqrt{\omega^2+4\kappa\coth\kappa}-\omega\right)-4\right],
\end{align*}

\begin{align*}
a_{44}^{(1,0)}&={}\frac{\pi}{\kappa\,\omega\,\left(\kappa\coth\kappa-1\right)} \Big[\left(\sqrt{\omega^2+4\kappa\coth\kappa}-\omega\right)\,\exp\!\left[\frac{\omega}{2}\left(\sqrt{\omega^2+4\kappa\coth\kappa}-\omega\right)\right]\\*[1mm]
&{}\qquad\qquad\qquad\qquad\qquad-\left(\sqrt{\omega^2+4\kappa\coth\kappa}-\omega\right)-2\,\omega\Big].
\end{align*}

The remaining first-order entries are
\begin{equation*}
 a_{12}^{(1,0)}=a_{14}^{(1,0)}=a_{21}^{(1,0)}=a_{24}^{(1,0)} =a_{31}^{(1,0)}=a_{32}^{(1,0)}=a_{34}^{(1,0)}=0.
\end{equation*}

\subsection{Coefficients of order (0,1)}

\begin{align*}
a_{13}^{(0,1)}=a_{23}^{(0,1)}&={}\frac{\pi\,\sinh^{2}\kappa}{8\,\kappa\,\left(\sqrt{\omega^2+4\kappa\coth\kappa}-\omega\right)^{2}\,\left(\sinh\kappa\cosh\kappa-\kappa\right)}\\[2mm]
&\qquad{}\times \Big[16\,\kappa^{2}\,\left(\omega\,\left(\sqrt{\omega^2+4\kappa\coth\kappa}-\omega\right)-2\right)\\*[1mm]
&\qquad\qquad+\left(6-\omega\,\left(\sqrt{\omega^2+4\kappa\coth\kappa}-\omega\right)\right)\,\left(\sqrt{\omega^2+4\kappa\coth\kappa}-\omega\right)^{4}\\*[1mm]
&\qquad\qquad+4\,\omega\,\left(\sqrt{\omega^2+4\kappa\coth\kappa}-\omega\right)^{3}\Big],
\end{align*}

\begin{align*}
a_{41}^{(0,1)}=-a_{42}^{(0,1)}&={}-\frac{\mathrm{i}\,\pi}{32\,\kappa\,\sinh^{2}\kappa\,\left(\kappa\coth\kappa-1\right)}\left[4\,\kappa^{2}-\omega^{2}\,\left(\sqrt{\omega^2+4\kappa\coth\kappa}-\omega\right)^{2}\right]^{-1}\\[2mm]
&\qquad{}\times \Big[128\,\kappa^{4}\,\left(1+2\,\sinh^{2}\kappa\right)\\*[1mm]
&\qquad\qquad+64\,\kappa^{2}\,\left(\sqrt{\omega^2+4\kappa\coth\kappa}-\omega\right)^{2}\,\sinh^{2}\kappa\,\left(1+2\,\sinh^{2}\kappa\right)\\*[1mm]
&\qquad\qquad-8\,\kappa^{2}\,\left(\sqrt{\omega^2+4\kappa\coth\kappa}-\omega\right)^{4}\,\sinh^{2}\kappa\\*[1mm]
&\qquad\qquad-16\,\kappa^{2}\,\omega\,\left(\sqrt{\omega^2+4\kappa\coth\kappa}-\omega\right)^{3}\,\sinh^{2}\kappa\,\left(1+4\,\sinh^{2}\kappa\right)\\*[1mm]
&\qquad\qquad+\left(\sqrt{\omega^2+4\kappa\coth\kappa}-\omega\right)^{6}\,\sinh^{4}\kappa\,\left(\omega\,\left(\sqrt{\omega^2+4\kappa\coth\kappa}-\omega\right)-4\right)\\*[1mm]
&\qquad\qquad-8\,\omega\,\left(\sqrt{\omega^2+4\kappa\coth\kappa}-\omega\right)^{5}\,\sinh^{4}\kappa\\*[1mm]
&\qquad\qquad-8\,\kappa\,\omega\,\left(\sqrt{\omega^2+4\kappa\coth\kappa}-\omega\right)^{5}\,\exp\!\left[\frac{\omega}{2}\left(\sqrt{\omega^2+4\kappa\coth\kappa}-\omega\right)\right]\,\sinh^{3}\kappa\\*[1mm]
&\qquad\qquad-32\,\kappa\,\omega\,\left(\sqrt{\omega^2+4\kappa\coth\kappa}-\omega\right)^{3}\,\exp\!\left[\frac{\omega}{2}\left(\sqrt{\omega^2+4\kappa\coth\kappa}-\omega\right)\right]\,\sinh^{3}\kappa\Big],
\end{align*}

\begin{align*}
a_{43}^{(0,1)}&={}\frac{\pi\,\sinh\kappa}{4\,\kappa\,\omega\,\left(\kappa\coth\kappa-1\right)} \Big[\left(\sqrt{\omega^2+4\kappa\coth\kappa}-\omega\right)^{3}\,\exp\!\left[\frac{\omega}{2}\left(\sqrt{\omega^2+4\kappa\coth\kappa}-\omega\right)\right]\\*[1mm]
&\qquad\qquad\qquad\qquad\qquad+4\,\left(\sqrt{\omega^2+4\kappa\coth\kappa}-\omega\right)\,\exp\!\left[\frac{\omega}{2}\left(\sqrt{\omega^2+4\kappa\coth\kappa}-\omega\right)\right]\\*[1mm]
&\qquad\qquad\qquad\qquad\qquad-\left(1+\omega^{2}\right)\,\left(\sqrt{\omega^2+4\kappa\coth\kappa}-\omega\right)^{3}\\*[1mm]
&\qquad\qquad\qquad\qquad\qquad-4\,\left(\sqrt{\omega^2+4\kappa\coth\kappa}-\omega\right)-8\,\omega\,\kappa^{2}+2\,\omega\,\left(\sqrt{\omega^2+4\kappa\coth\kappa}-\omega\right)^{2}\Big],
\end{align*}

All other entries of this matrix vanish:
\begin{align*}
 a_{11}^{(0,1)}=a_{12}^{(0,1)}=a_{14}^{(0,1)} =a_{21}^{(0,1)}=a_{22}^{(0,1)}=a_{24}^{(0,1)}&=0,\\
 a_{31}^{(0,1)}=a_{32}^{(0,1)}=a_{33}^{(0,1)} =a_{34}^{(0,1)}=a_{44}^{(0,1)}&=0.
\end{align*}

\subsection{The five required second-order coefficients}

\begin{align*}
a_{11}^{(2,0)}&={}\frac{\pi^{2}\,\sinh^{2}\kappa}{8\,\kappa^{4}\,\left(\sinh\kappa\cosh\kappa-\kappa\right)^{2}} \left(4\,\kappa\,\cosh\kappa+\left(\sqrt{\omega^2+4\kappa\coth\kappa}-\omega\right)^{2}\,\sinh\kappa\right)^{2}\\[2mm]
&{}\quad+\frac{\mathrm{i}\,\pi\,\omega\,\left(\sqrt{\omega^2+4\kappa\coth\kappa}-\omega\right)^{3}}{8\,\kappa^{3}\,\left(\kappa\coth\kappa-1\right)\,\left(\sinh\kappa\cosh\kappa-\kappa\right)}\\[2mm]
&\qquad{}\times \left[4\,\kappa^{2}-\omega^{2}\,\left(\sqrt{\omega^2+4\kappa\coth\kappa}-\omega\right)^{2}\right]^{-2}\\[2mm]
&\qquad{}\times \left[4\,\kappa-\left(\sqrt{\omega^2+4\kappa\coth\kappa}-\omega\right)^{2}\,\exp\!\left[\frac{\omega}{2}\left(\sqrt{\omega^2+4\kappa\coth\kappa}-\omega\right)\right]\,\sinh\kappa\right]\\[2mm]
&\qquad{}\times \Big[\left(16\,\kappa^{2}+\omega\,\left(\sqrt{\omega^2+4\kappa\coth\kappa}-\omega\right)^{3}\right)\,\sinh\kappa\\[1mm]
&\qquad\qquad-4\,\kappa\,\omega\,\left(\sqrt{\omega^2+4\kappa\coth\kappa}-\omega\right)\,\exp\!\left[\frac{\omega}{2}\left(\sqrt{\omega^2+4\kappa\coth\kappa}-\omega\right)\right]\Big]\\[2mm]
&{}\quad-\frac{\mathrm{i}\,\pi\,\sinh^{3}\kappa}{8\,\cosh\kappa\,\kappa^{4}\,\left(\sinh\kappa\cosh\kappa-\kappa\right)^{3}}\\[2mm]
&\qquad{}\times \Big[\kappa^{2}\,\left(4\,\kappa\,\cosh\kappa+\left(\sqrt{\omega^2+4\kappa\coth\kappa}-\omega\right)^{2}\,\sinh\kappa\right)^{2}\\*[1mm]
&\qquad\qquad\qquad+\left(\sqrt{\omega^2+4\kappa\coth\kappa}-\omega\right)^{4}\,\left(\sinh\kappa\cosh\kappa-\kappa\right)^{2}\Big],
\end{align*}

\begin{align*}
a_{34}^{(2,0)}&={}\frac{\pi\,\left(\sqrt{\omega^2+4\kappa\coth\kappa}-\omega\right)^{2}}{2\,\kappa\,\omega^{2}\,\left(\kappa\coth\kappa-1\right)^{2}} \Big[\omega^{2}\,\exp\!\left[\frac{\omega}{2}\left(\sqrt{\omega^2+4\kappa\coth\kappa}-\omega\right)\right]\\*[1mm]
&\qquad\qquad\qquad\qquad\qquad\qquad\qquad-\left(\exp\!\left[\frac{\omega}{2}\left(\sqrt{\omega^2+4\kappa\coth\kappa}-\omega\right)\right]-1\right)^{2}\Big],
\end{align*}

\begin{align*}
a_{14}^{(1,1)}&={}\frac{\pi\,\sinh^{2}\kappa}{64\,\left(\sqrt{\omega^2+4\kappa\coth\kappa}-\omega\right)^{2}\,\kappa^{2}\,\left(\sinh\kappa\cosh\kappa-\kappa\right)^{2}}\\[2mm]
&\qquad{}\times \Big[256\,\kappa^{4}\,\sinh^{2}\kappa-32\,\kappa^{2}\,\left(\sqrt{\omega^2+4\kappa\coth\kappa}-\omega\right)^{4}\,\sinh^{2}\kappa\\*[1mm]
&\qquad\qquad+\left(\sqrt{\omega^2+4\kappa\coth\kappa}-\omega\right)^{8}\,\sinh^{2}\kappa+8\,\left(\sqrt{\omega^2+4\kappa\coth\kappa}-\omega\right)^{6}\,\sinh^{2}\kappa\\*[1mm]
&\qquad\qquad+64\,\kappa^{2}\,\left(\sqrt{\omega^2+4\kappa\coth\kappa}-\omega\right)^{2}+32\,\kappa\,\left(\sqrt{\omega^2+4\kappa\coth\kappa}-\omega\right)^{4}\,\sinh\kappa\,\cosh\kappa\Big]\\[2mm]
&{}\quad-\frac{\pi\,\sinh^{2}\kappa}{8\,\kappa\,\left(\sqrt{\omega^2+4\kappa\coth\kappa}-\omega\right)^{2}\,\left(\sinh\kappa\cosh\kappa-\kappa\right)}\\[2mm]
&\qquad{}\times \Big[16\,\kappa^{2}\,\left(\omega\,\left(\sqrt{\omega^2+4\kappa\coth\kappa}-\omega\right)-2\right)\\*[1mm]
&\qquad\qquad+\left(6-\omega\,\left(\sqrt{\omega^2+4\kappa\coth\kappa}-\omega\right)\right)\,\left(\sqrt{\omega^2+4\kappa\coth\kappa}-\omega\right)^{4}\\*[1mm]
&\qquad\qquad+4\,\omega\,\left(\sqrt{\omega^2+4\kappa\coth\kappa}-\omega\right)^{3}\Big]\\[2mm]
&\qquad{}\times \Big[\frac{\sinh\kappa\,\cosh\kappa}{\left(\sinh\kappa\cosh\kappa-\kappa\right)}+\frac{\left(\sqrt{\omega^2+4\kappa\coth\kappa}-\omega\right)^{2}\,\sinh^{2}\kappa}{4\,\kappa\,\left(\sinh\kappa\cosh\kappa-\kappa\right)}\\*[1mm]
&\qquad\qquad-\frac{\left(\sqrt{\omega^2+4\kappa\coth\kappa}-\omega\right)\,\exp\!\left[\frac{\omega}{2}\left(\sqrt{\omega^2+4\kappa\coth\kappa}-\omega\right)\right]}{2\,\omega\,\left(\kappa\coth\kappa-1\right)}\\*[1mm]
&\qquad\qquad+\frac{\left(\sqrt{\omega^2+4\kappa\coth\kappa}-\omega\right)}{2\,\omega\,\left(\kappa\coth\kappa-1\right)}+\frac{1}{\left(\kappa\coth\kappa-1\right)}\Big],
\end{align*}

\begin{align*}
a_{31}^{(1,1)}&={}\frac{\mathrm{i}\,\pi}{64\,\kappa^{2}\,\left(\kappa\coth\kappa-1\right)\,\left(\sinh\kappa\cosh\kappa-\kappa\right)}\\[2mm]
&\qquad{}\times \Big[64\,\kappa^{4}+16\,\kappa^{2}\,\left(\sqrt{\omega^2+4\kappa\coth\kappa}-\omega\right)^{2}\,\sinh^{2}\kappa\,\left(\sinh^{2}\kappa+2\right)\\*[1mm]
&\qquad\qquad+8\,\kappa\,\left(\sqrt{\omega^2+4\kappa\coth\kappa}-\omega\right)^{4}\,\sinh^{3}\kappa\,\cosh\kappa+\left(\sqrt{\omega^2+4\kappa\coth\kappa}-\omega\right)^{6}\,\sinh^{4}\kappa\Big]\\[2mm]
&{}+\frac{\mathrm{i}\,\pi}{32\,\kappa\,\sinh^{2}\kappa\,\left(\kappa\coth\kappa-1\right)}\left[4\,\kappa^{2}-\omega^{2}\,\left(\sqrt{\omega^2+4\kappa\coth\kappa}-\omega\right)^{2}\right]^{-1}\\[2mm]
&\qquad{}\times \Big[128\,\kappa^{4}\,\left(1+2\,\sinh^{2}\kappa\right)\\*[1mm]
&\qquad\qquad+64\,\kappa^{2}\,\left(\sqrt{\omega^2+4\kappa\coth\kappa}-\omega\right)^{2}\,\sinh^{2}\kappa\,\left(1+2\,\sinh^{2}\kappa\right)\\*[1mm]
&\qquad\qquad-8\,\kappa^{2}\,\left(\sqrt{\omega^2+4\kappa\coth\kappa}-\omega\right)^{4}\,\sinh^{2}\kappa\\*[1mm]
&\qquad\qquad-16\,\kappa^{2}\,\omega\,\left(\sqrt{\omega^2+4\kappa\coth\kappa}-\omega\right)^{3}\,\sinh^{2}\kappa\,\left(1+4\,\sinh^{2}\kappa\right)\\*[1mm]
&\qquad\qquad+\left(\sqrt{\omega^2+4\kappa\coth\kappa}-\omega\right)^{6}\,\sinh^{4}\kappa\,\left(\omega\,\left(\sqrt{\omega^2+4\kappa\coth\kappa}-\omega\right)-4\right)\\*[1mm]
&\qquad\qquad-8\,\omega\,\left(\sqrt{\omega^2+4\kappa\coth\kappa}-\omega\right)^{5}\,\sinh^{4}\kappa\\*[1mm]
&\qquad\qquad-8\,\kappa\,\omega\,\left(\sqrt{\omega^2+4\kappa\coth\kappa}-\omega\right)^{5}\,\exp\!\left[\frac{\omega}{2}\left(\sqrt{\omega^2+4\kappa\coth\kappa}-\omega\right)\right]\,\sinh^{3}\kappa\\*[1mm]
&\qquad\qquad-32\,\kappa\,\omega\,\left(\sqrt{\omega^2+4\kappa\coth\kappa}-\omega\right)^{3}\,\exp\!\left[\frac{\omega}{2}\left(\sqrt{\omega^2+4\kappa\coth\kappa}-\omega\right)\right]\,\sinh^{3}\kappa\Big]\\[2mm]
&\qquad{}\times \Big[\frac{\sinh\kappa\,\cosh\kappa}{\left(\sinh\kappa\cosh\kappa-\kappa\right)}+\frac{\left(\sqrt{\omega^2+4\kappa\coth\kappa}-\omega\right)^{2}\,\sinh^{2}\kappa}{4\,\kappa\,\left(\sinh\kappa\cosh\kappa-\kappa\right)}-\frac{\left(\sqrt{\omega^2+4\kappa\coth\kappa}-\omega\right)\,\left(1+\omega^{2}\right)}{2\,\omega\,\left(\kappa\coth\kappa-1\right)}\\*[1mm]
&\qquad\qquad+\frac{\left(\sqrt{\omega^2+4\kappa\coth\kappa}-\omega\right)\,\exp\!\left[\frac{\omega}{2}\left(\sqrt{\omega^2+4\kappa\coth\kappa}-\omega\right)\right]}{2\,\omega\,\left(\kappa\coth\kappa-1\right)}+\frac{1}{\left(\kappa\coth\kappa-1\right)}\Big],
\end{align*}

\begin{align*}
a_{11}^{(0,2)}&={}-\frac{\mathrm{i}\,\pi}{512\,\kappa\,\left(\sqrt{\omega^2+4\kappa\coth\kappa}-\omega\right)^{2}\,\left(\sinh\kappa\cosh\kappa-\kappa\right)\,\left(\kappa\coth\kappa-1\right)}\\[2mm]
&\qquad{}\times \left[4\,\kappa^{2}-\omega^{2}\,\left(\sqrt{\omega^2+4\kappa\coth\kappa}-\omega\right)^{2}\right]^{-1}\\[2mm]
&\qquad{}\times \Big[16\,\kappa^{2}\,\left(\omega\,\left(\sqrt{\omega^2+4\kappa\coth\kappa}-\omega\right)-2\right)\\*[1mm]
&\qquad\qquad+\left(6-\omega\,\left(\sqrt{\omega^2+4\kappa\coth\kappa}-\omega\right)\right)\,\left(\sqrt{\omega^2+4\kappa\coth\kappa}-\omega\right)^{4}\\[1mm]
&\qquad\qquad+4\,\omega\,\left(\sqrt{\omega^2+4\kappa\coth\kappa}-\omega\right)^{3}\Big]\\[2mm]
&\qquad{}\times \Big[128\,\kappa^{4}\,\left(1+2\,\sinh^{2}\kappa\right)\\[1mm]
&\qquad\qquad+64\,\kappa^{2}\,\left(\sqrt{\omega^2+4\kappa\coth\kappa}-\omega\right)^{2}\,\sinh^{2}\kappa\,\left(1+2\,\sinh^{2}\kappa\right)\\[1mm]
&\qquad\qquad-8\,\kappa^{2}\,\left(\sqrt{\omega^2+4\kappa\coth\kappa}-\omega\right)^{4}\,\sinh^{2}\kappa\\[1mm]
&\qquad\qquad-16\,\kappa^{2}\,\omega\,\left(\sqrt{\omega^2+4\kappa\coth\kappa}-\omega\right)^{3}\,\sinh^{2}\kappa\,\left(1+4\,\sinh^{2}\kappa\right)\\[1mm]
&\qquad\qquad+\left(\sqrt{\omega^2+4\kappa\coth\kappa}-\omega\right)^{6}\,\sinh^{4}\kappa\,\left(\omega\,\left(\sqrt{\omega^2+4\kappa\coth\kappa}-\omega\right)-4\right)\\[1mm]
&\qquad\qquad-8\,\omega\,\left(\sqrt{\omega^2+4\kappa\coth\kappa}-\omega\right)^{5}\,\sinh^{4}\kappa\\[1mm]
&\qquad\qquad-8\,\kappa\,\omega\,\left(\sqrt{\omega^2+4\kappa\coth\kappa}-\omega\right)^{5}\,\exp\!\left[\frac{\omega}{2}\left(\sqrt{\omega^2+4\kappa\coth\kappa}-\omega\right)\right]\,\sinh^{3}\kappa\\[1mm]
&\qquad\qquad-32\,\kappa\,\omega\,\left(\sqrt{\omega^2+4\kappa\coth\kappa}-\omega\right)^{3}\,\exp\!\left[\frac{\omega}{2}\left(\sqrt{\omega^2+4\kappa\coth\kappa}-\omega\right)\right]\,\sinh^{3}\kappa\Big]\\[2mm]
&{}-\frac{\mathrm{i}\,\pi}{1024\,\kappa^{2}\,\sinh\kappa\,\left(\sinh\kappa\cosh\kappa-\kappa\right)\,\left(\sqrt{\omega^2+4\kappa\coth\kappa}-\omega\right)^{2}}\\[2mm]
&\qquad{}\times \Big[2048\,\kappa^{5}\,\sinh^{3}\kappa+256\,\kappa^{4}\,\left(\sqrt{\omega^2+4\kappa\coth\kappa}-\omega\right)^{2}\,\cosh\kappa\,\left(1-3\,\sinh^{2}\kappa+3\,\sinh^{4}\kappa\right)\\*[1mm]
&\qquad\qquad+128\,\kappa^{3}\,\left(\sqrt{\omega^2+4\kappa\coth\kappa}-\omega\right)^{4}\,\sinh\kappa\,\left(1-\sinh^{2}\kappa\right)\\*[1mm]
&\qquad\qquad+48\,\kappa^{2}\,\left(\sqrt{\omega^2+4\kappa\coth\kappa}-\omega\right)^{6}\,\cosh\kappa\,\sinh^{2}\kappa\\*[1mm]
&\qquad\qquad+8\,\kappa\,\left(\sqrt{\omega^2+4\kappa\coth\kappa}-\omega\right)^{8}\,\sinh^{3}\kappa\,\left(1+2\,\sinh^{2}\kappa\right)\\*[1mm]
&\qquad\qquad+\left(\sqrt{\omega^2+4\kappa\coth\kappa}-\omega\right)^{10}\,\cosh\kappa\,\sinh^{4}\kappa\Big].
\end{align*}

\subsection{Relations between the remaining coefficients}
\label{app:remaining-low-coefficient-relations}

The substitution leading to \eqref{eqn:explicit-S2-low} also uses \eqref{eqn:relations-amn} and the following identities:
\begin{align*}
 a_{14}^{(2,0)}
 &=-\frac{a_{13}^{(1,0)}(a_{11}^{(1,0)}-a_{44}^{(1,0)})}{T},& a_{21}^{(2,0)}&=-a_{12}^{(2,0)},\\
 a_{22}^{(2,0)}&=(a_{11}^{(1,0)})^2-a_{11}^{(2,0)},& a_{32}^{(2,0)}&=a_{31}^{(2,0)},\\
 a_{24}^{(2,0)}&=-a_{14}^{(2,0)},& a_{21}^{(1,1)}&=a_{12}^{(1,1)},\\
 a_{22}^{(1,1)}&=a_{11}^{(1,1)},& a_{24}^{(1,1)}&=a_{14}^{(1,1)},\\
 a_{32}^{(1,1)}&=-a_{31}^{(1,1)},& a_{34}^{(1,1)}&=0,\\
 a_{12}^{(0,2)}&=-a_{11}^{(0,2)},& a_{21}^{(0,2)}&=a_{11}^{(0,2)},\\
 a_{22}^{(0,2)}&=-a_{11}^{(0,2)},& a_{j4}^{(0,2)}&=0\quad(j=1,\ldots,4).
\end{align*}
These identities follow from the Duhamel formula \eqref{eqn:duhamel-low}.  The last also follows from Lemma~\ref{lem:fourth-neutral-column}.

The two remaining entries in the last row of \eqref{eqn:explicit-S2-low} satisfy
\begin{equation}\label{eqn:pure-amplitude-third-row-low}
 a_{31}^{(0,2)}=a_{32}^{(0,2)}=0.
\end{equation}
These two zeros can also be deduced from amplitude parity.  The odd part in $\alpha$ of $T\det S_2$, computed from \eqref{eqn:explicit-S2-low}, is
\[
 iTc_gm\alpha^3 \left[a_{14}^{(1,1)}+ \frac{a_{13}^{(0,1)}(a_{11}^{(1,0)}-a_{44}^{(1,0)})}{T}\right]\left[ a_{11}^{(1,0)}\zeta(a_{31}^{(0,2)}+a_{32}^{(0,2)}) -\frac{c_g^2m^2F_{\mathrm{low}}}{2} (a_{31}^{(0,2)}-a_{32}^{(0,2)})\right].
\]
It vanishes by \eqref{eqn:Evans-amplitude-parity-low}. The coefficient formulas above give, with $\mu_0$ as in \eqref{eqn:rho-t-mu-low},
\[
 a_{14}^{(1,1)}+ \frac{a_{13}^{(0,1)}(a_{11}^{(1,0)}-a_{44}^{(1,0)})}{T} = \frac{\pi\sinh^2\kappa} {\mu_0^2\kappa^2(\sinh\kappa\cosh\kappa-\kappa)^2} \Bigl[\sinh^2\kappa\bigl((\kappa^2-\mu_0^4)^2+2\mu_0^6\bigr) +\kappa^2\mu_0^2 +2\kappa\mu_0^4\sinh\kappa\cosh\kappa\Bigr]>0.
\]
Moreover, $a_{11}^{(1,0)}=T/c_g\ne0$ and $F_{\mathrm{low}}=-2iTc_2/c_g^3\ne0$. Comparing the coefficients of $m\alpha^3\zeta$ and $m^3\alpha^3$ gives $a_{31}^{(0,2)}+a_{32}^{(0,2)}=0$ and $a_{31}^{(0,2)}-a_{32}^{(0,2)}=0$, proving \eqref{eqn:pure-amplitude-third-row-low}. \endgroup

\section{Explicit formula for the low-frequency index}
\label{index:low}

For the algebraic simplification in this appendix, use $t=\tanh\kappa$ and $\rho=\omega\mu_0$ as in \eqref{eqn:rho-t-mu-low}, and set
\begin{equation*} X=-\frac{\rho t}{\kappa}.
\end{equation*}
Define
\begin{align*}
 U={}&9-12t^2+13t^4-2t^6 +(27-18t^2+15t^4)X \\
 &+(33-3t^2+4t^4)X^2 +(21+5t^2)X^3+(7+2t^2)X^4+X^5,\\
 V={}&(1+X)^2(2-c_g-\rho)+1+X +(c_g+\rho-1)t^2,\\
 N_W={}&(1+X)(2+X)+(1-c_g)(1-t^2), \\
 D_W={}&t(1-c_g)(1-c_g-\rho)-\kappa(1+X), \qquad W=\frac{2t^3N_W}{D_W}.
\end{align*}
Substituting the coefficient formulas from Appendix~\ref{app:all-explicit-low-order-coefficients} into \eqref{eqn:def-G-low} and using the dispersion relation gives
\begin{equation}\label{eqn:simplified-G-low}
 G_{\mathrm{low}} =-\frac{i\pi^3\sinh^2\kappa\,\Delta_L} {c_g^3(1-\kappa\coth\kappa)(1+X)(2+X)t^4}\,(U+V W).
\end{equation}
These are exact algebraic identities.

Combining \eqref{eqn:index-from-amn}, \eqref{eqn:simplified-FQ-low}, and \eqref{eqn:simplified-G-low}, we obtain
\begin{equation}\label{eqn:ind1}
 \mathrm{ind}_{\mathrm{low}}(\omega,\kappa) =-\frac{\kappa^3c_2\sinh^2\kappa} {4(1+X)(2+X)t^4}\,(U+V W).
\end{equation}
Since $1+X=\mu_0^2t/\kappa>0$ and $c_2<0$, the factor multiplying $U+VW$ in \eqref{eqn:ind1} is positive.  Hence
\begin{equation*}
\operatorname{sign} (\mathrm{ind}_{\mathrm{low}}(\omega,\kappa)) = \operatorname{sign} (U+V W).
\end{equation*} 
The formula depends only on $\omega$ and $\kappa$ through the quantities defined above.  Substituting their definitions gives the following expression without auxiliary parameters:
\begingroup\allowdisplaybreaks[3]
\begin{align*}
\mathrm{ind}_{\mathrm{low}}(\omega,\kappa) ={}& \frac{\kappa^3 \bigl(\sqrt{\omega^2+4\kappa\coth\kappa}+\omega\bigr)^2} {64(\omega^2+4\kappa\coth\kappa)^2\tanh^3\!\kappa} \\[0.2mm]
&\quad\times \Biggl[ (\omega^2+4\kappa\coth\kappa) \bigl(\sqrt{\omega^2+4\kappa\coth\kappa}+\omega\bigr)^2\tanh^2\!\kappa +\frac{4\bigl[\tanh\kappa-\kappa(1-\tanh^2\!\kappa)\bigr]^2}{\tanh^2\!\kappa(1-\tanh^2\!\kappa)} \Biggr] \\[0.2mm]
&\quad\times\Biggl\{ \left( \frac{\sqrt{\omega^2+4\kappa\coth\kappa}-\omega} {\sqrt{\omega^2+4\kappa\coth\kappa}+\omega} \right)^5 \\[0.2mm]
&\qquad +2(1+\tanh^2\!\kappa) \left( \frac{\sqrt{\omega^2+4\kappa\coth\kappa}-\omega} {\sqrt{\omega^2+4\kappa\coth\kappa}+\omega} \right)^4  +3(1-\tanh^2\!\kappa) \left( \frac{\sqrt{\omega^2+4\kappa\coth\kappa}-\omega} {\sqrt{\omega^2+4\kappa\coth\kappa}+\omega} \right)^3 \\[0.2mm]
&\qquad +2(1-\tanh^2\!\kappa)(1-2\tanh^2\!\kappa) \left( \frac{\sqrt{\omega^2+4\kappa\coth\kappa}-\omega} {\sqrt{\omega^2+4\kappa\coth\kappa}+\omega} \right)^2 \\[0.2mm]
&\qquad +(1-5\tanh^2\!\kappa+7\tanh^4\!\kappa) \left( \frac{\sqrt{\omega^2+4\kappa\coth\kappa}-\omega} {\sqrt{\omega^2+4\kappa\coth\kappa}+\omega} \right) +2\tanh^4\!\kappa(1-\tanh^2\!\kappa) \\[0.2mm]
&\qquad +2\tanh^2\!\kappa\Biggl[ \tanh\kappa(1-\tanh^2\!\kappa) \left( \frac{\sqrt{\omega^2+4\kappa\coth\kappa}-\omega} {\sqrt{\omega^2+4\kappa\coth\kappa}+\omega} \right) \\[0.2mm]
&\qquad\quad +2\tanh\kappa \left( \frac{\sqrt{\omega^2+4\kappa\coth\kappa}-\omega} {\sqrt{\omega^2+4\kappa\coth\kappa}+\omega} \right)^2 +2\tanh\kappa \left( \frac{\sqrt{\omega^2+4\kappa\coth\kappa}-\omega}
     {\sqrt{\omega^2+4\kappa\coth\kappa}+\omega} \right)^3 \\[0.2mm]
&\qquad\quad +\kappa\left( \left( \frac{\sqrt{\omega^2+4\kappa\coth\kappa}-\omega} {\sqrt{\omega^2+4\kappa\coth\kappa}+\omega} \right)^2-\tanh^2\!\kappa \right)^2 \Biggr] \\[0.2mm]
&\qquad\quad\times \Biggl[ \tanh\kappa(2-\tanh^2\!\kappa) \left( \frac{\sqrt{\omega^2+4\kappa\coth\kappa}-\omega} {\sqrt{\omega^2+4\kappa\coth\kappa}+\omega} \right) \\[0.2mm]
&\qquad\qquad +2\tanh\kappa \left( \frac{\sqrt{\omega^2+4\kappa\coth\kappa}-\omega} {\sqrt{\omega^2+4\kappa\coth\kappa}+\omega} \right)^2 +\tanh\kappa \left( \frac{\sqrt{\omega^2+4\kappa\coth\kappa}-\omega}
     {\sqrt{\omega^2+4\kappa\coth\kappa}+\omega} \right)^3 \\[0.2mm]
&\qquad\qquad{}+\kappa(1-\tanh^2\!\kappa)^2 \Biggr] \\[0.2mm]
&\qquad\quad\times \Biggl[ \bigl(\tanh^2\!\kappa-2\kappa \tanh\kappa+\kappa^2(1-\tanh^2\!\kappa)\bigr) \left( \frac{\sqrt{\omega^2+4\kappa\coth\kappa}-\omega} {\sqrt{\omega^2+4\kappa\coth\kappa}+\omega} \right)^2 \\[0.2mm]
&\qquad\qquad -2\kappa \tanh^3\!\kappa \left( \frac{\sqrt{\omega^2+4\kappa\coth\kappa}-\omega} {\sqrt{\omega^2+4\kappa\coth\kappa}+\omega} \right) -\kappa^2\tanh^2\!\kappa(1-\tanh^2\!\kappa) \Biggr]^{-1} \Biggr\}.
\end{align*}
\endgroup

\section{The relation between the high-frequency coupling coefficients}
\label{A:high-index}

We prove \eqref{eqn:high-beta-reciprocity} by comparing the two formulas in \eqref{eqn:high-a12-explicit}--\eqref{eqn:high-a21-explicit}. The calculation has two parts: evaluate the direct second-order term and the two first-order interactions, then use the dispersion relations at $k_2$ and $k_4$.

Throughout this appendix, $k_2-k_4=2\kappa$ and $q_*=k_2-\kappa=k_4+\kappa$, as in \eqref{eqn:high-middle-resolvent}.  We use the dispersion function $d(k,\lambda)$ from \eqref{eqn:disp}.  At $\lambda=i\sigma_H$,
\begin{equation*} d(\xi,i\sigma_H) =\bigl[\mu_0^2\xi+\rho(\xi-\sigma_H)\bigr]\tanh\xi -(\xi-\sigma_H)^2.
\end{equation*}
In particular, $d(k_2,i\sigma_H)=d(k_4,i\sigma_H)=0$. At the intermediate wave number,
\[
 d(q_*,i\sigma_H) =-(\sigma_H-\sigma_+(q_*))(\sigma_H-\sigma_-(q_*))>0.
\]
Indeed, $k_4<q_*<k_2$, and the branch ordering places $\sigma_H$ strictly between $\sigma_-(q_*)$ and $\sigma_+(q_*)$.

We first record the finite products obtained from the coefficient calculation.  For an interaction from wave number $\zeta$ to wave number $\xi$, the first-order terms give
\begin{equation}\label{eqn:high-algebra-N1}
 N_1(\xi,\zeta)=
 \begin{pmatrix}
 -\xi\zeta&-\mu_0^2\xi\\
 \zeta\bigl[\mu_0^2+\rho-(\xi-\sigma_H)\tanh\zeta\bigr]&-\kappa^2
 \end{pmatrix}.
\end{equation}
The direct second-order terms give a matrix $N_2(\xi,\zeta)$ with entries
\begin{align*}
 (N_2)_{11}&=-\xi\zeta\left(\frac{2E}{\sinh^2\kappa} +\frac{\zeta\tanh\zeta}{2}\right), 
 \qquad (N_2)_{12}
 =\frac{\xi}{2}\left(2\mu_0^2\rho+\rho^2 -\frac{4\mu_0^2E+\kappa^2}{\sinh^2\kappa}\right), \\
 (N_2)_{21}&=\frac{\zeta}{2}\Bigg\{ \frac{4E}{\sinh^2\kappa} \bigl[\mu_0^2+\rho-(\xi-\sigma_H)\tanh\zeta\bigr]-2\mu_0^2\rho-\rho^2 +\frac{\kappa^2}{\sinh^2\kappa} -\zeta(\xi-\sigma_H) +\bigl[\rho\zeta+\mu_0^2(\xi+\zeta)\bigr]\tanh\zeta \Bigg\}, \\
 (N_2)_{22}&=\kappa^2\left(3\mu_0^2-\frac{8E}{\sinh^2\kappa}\right).
\end{align*}
Here $E$ is the Stokes coefficient already given in \eqref{eqn:coefE}.  All other quantities have been defined above.

For $(\xi,\zeta)=(k_2,k_4)$ or $(k_4,k_2)$, let
\begin{align*}
 P_{\xi\zeta} ={}&(\mu_0^2,\sigma_H-\xi) \Bigg[-d(q_*,i\sigma_H)N_2(\xi,\zeta) \\
 &\quad{}-N_1(\xi,q_*)
 \begin{pmatrix}
 \mu_0^2&\sigma_H-q_*\\
 q_*-\sigma_H-\rho\tanh q_*&-q_*\tanh q_*
 \end{pmatrix}
 N_1(q_*,\zeta)\Bigg]
 \begin{pmatrix}
 \mu_0^2\\ \zeta-\sigma_H-\rho\tanh\zeta
 \end{pmatrix}.
\end{align*}
The first term contains the direct second-order interaction.  The second contains the two first-order interactions through $q_*$. The row and column vectors retain the normalization of the endpoint modes.  This is the finite-dimensional form of the two terms in \eqref{eqn:high-a12-explicit}--\eqref{eqn:high-a21-explicit}.

More precisely, substitution into those formulas gives
\begin{align}
 \frac{\beta_{12}}{\alpha_1}
 &=\frac{\sinh^2\kappa\cosh k_4} {4\mu_0^2\cosh k_2\,d(q_*,i\sigma_H)} \frac{P_{k_2k_4}} {2(k_2-\sigma_H)-\rho\tanh k_2},
 \label{eqn:high-algebra-coefficients}\\
 \frac{\beta_{21}}{\alpha_2}
 &=\frac{\sinh^2\kappa\cosh k_2} {4\mu_0^2\cosh k_4\,d(q_*,i\sigma_H)} \frac{P_{k_4k_2}} {2(k_4-\sigma_H)-\rho\tanh k_4}. \notag
\end{align}
To obtain \eqref{eqn:high-algebra-N1}--\eqref{eqn:high-algebra-coefficients}, write the hyperbolic functions as exponentials and integrate each term in $y$.  The integrals are of the form
\[
 \int_0^1 y^j e^{ay}\,dy =\partial_a^j\left(\frac{e^a-1}{a}\right),
\]
with the continuous value at $a=0$.  Collect the direct second-order term and the two first-order terms separately, and use $\mu_0^2+\rho=\kappa\coth\kappa$ and $k_2-k_4=2\kappa$. The ratios of hyperbolic cosines in \eqref{eqn:high-algebra-coefficients} come from the potential components of $\Upphi_2,\Upphi_4$.  Division by $\alpha_j$ normalizes the coefficient of $\delta$.

We now compare the two numerators.  In the finite products, temporarily regard $\tanh k_2$ and $\tanh k_4$ as independent variables. Then $P_{k_2k_4}$ is affine in $\tanh k_4$ and independent of $\tanh k_2$; the roles are reversed for $P_{k_4k_2}$. Write $[\tanh k_j]P$ for the coefficient of $\tanh k_j$ in such a polynomial.  Multiplication of the matrices gives
\begin{align*}
&\,k_2\bigl[\mu_0^2k_4+\rho(k_4-\sigma_H)\bigr]P_{k_2k_4} -k_4\bigl[\mu_0^2k_2+\rho(k_2-\sigma_H)\bigr]P_{k_4k_2} \\=&\,k_2\bigl([\tanh k_4]P_{k_2k_4}\bigr)d(k_4,i\sigma_H) -k_4\bigl([\tanh k_2]P_{k_4k_2}\bigr)d(k_2,i\sigma_H).
\end{align*}
This identity is checked by comparing the constant and linear terms in the two variables, using $q_*=(k_2+k_4)/2$ and $\kappa^2=(k_2-k_4)^2/4$.  The value of $E$ is not needed in this comparison.  At the resonant endpoints the right-hand side vanishes. Thus the dispersion relations give
\begin{equation*} \frac{(k_4-\sigma_H)^2}{k_4\tanh k_4}\,P_{k_2k_4} =\frac{(k_2-\sigma_H)^2}{k_2\tanh k_2}\,P_{k_4k_2}.
\end{equation*}
Insert this relation into \eqref{eqn:high-algebra-coefficients} and use
\[
 \frac{\cosh^2k_2}{\cosh^2k_4} \frac{k_2\tanh k_2}{k_4\tanh k_4} =\frac{k_2\sinh(2k_2)}{k_4\sinh(2k_4)}.
\]
The result is \eqref{eqn:high-beta-reciprocity} with the factor \eqref{eqn:high-reciprocity-factor}.  This argument uses no division by $P_{k_2k_4}$, $P_{k_4k_2}$, or either coupling coefficient, so it also covers their simultaneous vanishing.  The displayed denominators are nonzero by the branch ordering and $\mu_0>0$; removable factors in the intermediate dual coordinates are interpreted by analytic continuation.

\bibliographystyle{plain}
\bibliography{references.bib}
\addcontentsline{toc}{section}{References}

\end{document}